\documentclass{amsart}

\usepackage[english]{babel}
\usepackage[latin1]{inputenc}
\usepackage{amssymb}
\usepackage{amsmath}
\usepackage{amscd}
\usepackage{latexsym}
\usepackage{amsthm}
\usepackage{mathrsfs}
\usepackage{mathabx}
\usepackage{algorithm}
\usepackage{algorithmic}
\usepackage[all]{xy}
\usepackage{comment}
\usepackage{marginnote}
\theoremstyle{plain}
\newtheorem{thm}{Theorem}[section]
\newtheorem{theorem}[thm]{Theorem}
\usepackage[dvipsnames]{xcolor}
\usepackage{hyperref}

\newtheorem{cor}[thm]{Corollary}
\newtheorem{corollary}[thm]{Corollary}

\newtheorem{lem}[thm]{Lemma}
\newtheorem{lemma}[thm]{Lemma}

\newtheorem{prop}[thm]{Proposition}

\newtheorem{defn}[thm]{Definition}

\newtheorem{remark}[thm]{Remark}

\newtheorem{esem}[thm]{Example}
\newtheorem{example}[thm]{Example}

\numberwithin{equation}{section}
\newcommand{\rk}{{\rm rank}}

\newcommand{\ra}{\rightarrow}

\newcommand{\Z}{\mathbb{Z}}
\newcommand{\Q}{\mathbb{Q}}

\newcommand{\Aut}{\mbox{Aut}}

\usepackage{imakeidx}

\title{Singular symplectic surfaces}
\author{Alice Garbagnati, Matteo Penegini, Arvid Perego}

\begin{document}
\begin{abstract}
In this paper we classify the irreducible symplectic surfaces, i.e., compact, connected complex surfaces with canonical singularities that have a holomorphic symplectic form $\sigma$ on the smooth locus, and for which every finite quasi-\'etale covering has the algebra of reflexive forms spanned by the reflexive pull-back of $\sigma$. 
More precisely, we classify all singular symplectic surfaces distinguish them in primitive symplectic surfaces, irreducible symplectic surfaces and 2-dimensional irreducible symplectic orbifolds. Moreover, we prove that the Hilbert scheme of two points on such a surface $X$ is an irreducible symplectic variety, at least in the case where the smooth locus of $X$ is simply connected. 
\end{abstract}
\maketitle
\tableofcontents
\section{Introduction}
Irreducible symplectic manifolds are compact K\"ahler manifolds that are simply connected and carry a holomorphic symplectic form that spans the space of holomorphic $2-$forms. They appear as one of the three building blocks of compact K\"ahler manifolds with trivial first Chern class according to the Bogomolov Decomposition Theorem, the two others being complex tori and irreducible Calabi-Yau manifolds.

In dimension 2 the irreducible symplectic manifolds are exactly the K3 surfaces, and Hilbert schemes of points on K3 surfaces provide higher dimensional examples. Few other examples of irreducible symplectic manifolds have been described over the years, starting from moduli spaces of sheaves either on K3 surfaces or on Abelian surfaces.
	

The notion of irreducible symplectic manifold may be generalized in the singular setting in at least three non-equivalent ways. One of them is in the context of orbifolds, and was introduced by Campana in \cite{Cam}:

\begin{defn}\label{def: ISO}
An \textbf{irreducible symplectic orbifold} is a compact K\"ahler orbifold $X$ that verifies the three following conditions:
\begin{enumerate}
    \item the smooth locus $X^{s}$ of $X$ is simply connected;
    \item on $X^{s}$ there is a holomorphic symplectic form $\sigma$;
    \item the space $H^{0}(X^{s},\Omega_{X^{s}}^{2})$ is spanned by $\sigma$.
\end{enumerate}
\end{defn}

It is immediate to notice that an irreducible symplectic orbifold is smooth if and only if it is an irreducible symplectic manifold. Moreover, in \cite{Cam} a Bogomolov Decomposition Theorem for orbifolds is proved, and irreducible symplectic orbifolds appear as one of the three building blocks of compact K\"ahler orbifolds with trivial first Chern class. Examples of irreducible symplectic orbifolds appear in \cite{Men18}, \cite{Men14}, \cite{Men22}, \cite{BGMM} as terminalization of quotients of Hilbert schemes of points on K3 surfaces (or generalized Kummer varieties) by the action of symplectic automorphisms, or as quotients of products of K3 surfaces; in \cite{CGKK} a complete family of irreducible symplectic orbifolds is described, as deformation of a quotient of Hilbert schemes of two points on K3 surfaces by the action of a symplectic involution.

A second possible generalization to the singular setting of the notion of an irreducible symplectic manifold was proposed in \cite{GKP}:
 
\begin{defn}
\label{def: THE DEFINTION}
An \textbf{irreducible symplectic variety} is a normal compact K\"ahler space $X$ that verifies the three following conditions:
\begin{enumerate}
    \item it has canonical singularities,
    \item it admits a holomorphic symplectic form $\sigma$ on its smooth locus, that we will view as a global section of $\Omega_{X}^{[2]}:=(\Omega_{X}^{2})^{**}$,
    \item for every finite quasi-\'etale covering $f:Y\longrightarrow X$ we have that $$\bigoplus_{p}H^{0}(Y,\Omega_{Y}^{[p]})=\mathbb{C}[f^{[*]}\sigma]$$as a $\mathbb{C}-$algebra, where we recall that a finite quasi-\'etale morphism is a finite morphism which is \'etale in codimension 1.
\end{enumerate}
\end{defn} 

It is worthwhile to mention that as a consequence of the Bogomolov Decomposition Theorem for manifolds and of the fact that irreducible symplectic varieties are simply connected by \cite{GGK}, an irreducible symplectic variety is smooth if and only if it is an irreducible symplectic manifold. 

Moreover, in \cite{HP} and \cite{BGL} a generalization of the Bogomolov Decomposition Theorem for klt varities  is proved, and irreducible symplectic varieties are one of the three building blocks of normal compact K\"ahler spaces with klt singularities and trivial first Chern class. A consequence of this is that irreducible symplectic orbifolds are irreducible symplectic varieties (see Proposition 2.14 of \cite{Per}).  

Examples of irreducible symplectic varieties appear in \cite{Men18}, \cite{Men14}, \cite{Men22}, \cite{BGMM}, \cite{Per} as quotients of Hilbert schemes of points by the action of symplectic groups, in \cite{PR} as moduli spaces of sheaves on K3 or Abelian surfaces, or in \cite{BCGPSV} using relative Prym varieties.

A further generalization of the notion of an irreducible symplectic manifold to the singular setting was first considered in \cite{Fuj} in the context of orbifolds, and then in \cite{Nam}:
\begin{defn}
\label{defn:psv}
A \textbf{symplectic variety} is a normal variety which has a holomorphic symplectic form on its smooth locus, which extends to a holomorphic $2-$form on a (and hence any) resolution of the singularities. 

A \textbf{primitive symplectic variety} is a compact K\"ahler symplectic variety $X$ such that $h^{1}(X,\mathcal{O}_{X})=0$ and $H^{0}(X^{s},\Omega_{X^{s}}^{2})$ is spanned by the class of a holomorphic symplectic form $\sigma$.
\end{defn}

It is worthwhile to mention that irreducible symplectic varieties are primitive symplectic varieties (see Proposition 2.19 of \cite{Per}), and that since \cite{Nam} and \cite{BL} we know that on the free part of the second integral cohomology of a primitive symplectic variety there is an integral non-degenerate quadratic form and that a local and global Torelli Theorem hold. The symmetric product of a K3 surface is an example of primitive symplectic variety.

\begin{remark}
\label{rem:namikawa}
{\rm In the definition of a primitive symplectic variety $X$ we ask that $X$ is a symplectic variety. 
By Theorem 6 of \cite{Nam} a normal K\"ahler space $X$ having a holomorphic symplectic form on its smooth locus is a symplectic variety if and only if it has rational Gorenstein singularities.}  
\end{remark}

\begin{remark}
{\rm The notion of primitive symplectic variety appears in several papers and the definition is not always equivalent to Definition \ref{defn:psv}. An important example of this is given by the notion of \textit{primitive symplectic orbifold} that appears in \cite{FuMen}, which is more restrictive than that of primitive symplectic variety: following Definition \ref{defn:psv}, a primitive symplectic orbifold, as defined in \cite{FuMen}, is an orbifold that is a primitive symplectic variety with terminal singularities.} 
\end{remark}

The examples in \cite{PR} are all irreducible symplectic varieties which are not orbifolds, and hence not irreducible symplectic orbifolds. Moreover, the quotients of Hilbert schemes of 2 points on a K3 surface by a natural automorphism of prime order at least 3 is an orbifold and an irreducible symplectic variety (see Proposition 2.15 of \cite{Per}) but it is not an irreducible symplectic orbifold since its smooth locus is not simply connected. The symmetric products of K3 surfaces are examples of primitive symplectic varieties which are not irreducible symplectic.

Some examples of irreducible symplectic varieties  have a symplectic resolution of the singularities: in this case, the resolution is an irreducible symplectic manifold (see \cite{PR}, \cite{Per}). Other examples do not have any symplectic resolutions at all (see again \cite{PR}). There are moreover examples of singular symplectic varieties having a symplectic resolution of the singularities which is an irreducible symplectic manifolds and which are not irreducible symplectic varieties (e. g. the symmetric product of a K3 surface): they are however always primitive symplectic varieties (see Proposition 2.19 of \cite{Per}). 

A classification of primitive symplectic varieties, irreducible symplectic varieties and irreducible symplectic orbifolds is far from being obtained. One of the aims of the present paper is to study the case of the smallest possible dimension, i.e., the case of surfaces, and to construct examples of higher dimensional primitive symplectic varieties, irreducible symplectic varieties or irreducible symplectic orbifolds starting from this kind of singular surfaces. We will use the following:

\begin{defn}
Let $X$ be a normal, compact K\"ahler surface whose smooth locus has a holomorphic symplectic form.
\begin{enumerate}
    \item If $X$ is a  primitive symplectic variety, we will call it \textbf{primitive symplectic surface}.
    \item If $X$ is an  irreducible symplectic variety, we will call it \textbf{irreducible symplectic surface}.    \item If $X$ is an irreducible symplectic orbifold, we will call it \textbf{simple symplectic surface}.
\end{enumerate}
\end{defn}

The first aim of this paper is to present a classification of primitive symplectic surfaces, irreducible symplectic surfaces and simple symplectic surfaces. The smooth ones are the K3 surfaces, as well known, so we only need to classify the singular ones.

\begin{remark}
\label{remark:orbifold}
{\rm The following are known to hold: 
\begin{itemize} \item simple symplectic surfaces are irreducible symplectic, \item  irreducible symplectic surfaces are primitive symplectic,
\item primitive symplectic surfaces are compact K\"ahler holomorphic symplectic orbifolds,\item a primitive symplectic surface is simple if and only if its smooth locus is simply connected.\end{itemize}
The first two items follow by the discussions above.
As we will see their converses are false, since there are primitive symplectic surfaces which are not irreducible symplectic, and there are irreducible symplectic surfaces which are not simple symplectic.

 For the remaining items, we notice that a normal compact K\"ahler surface is symplectic if and only if it is a $2-$dimensional compact K\"ahler holomorphic symplectic orbifold, i.e., a compact K\"ahler orbifold of dimension 2 that carries a holomorphic symplectic form on its smooth locus.

 Indeed, if $X$ is symplectic then the singularities of $X$ are canonical. Since canonical singularities on a surfaces are ADE singularities, and these are all quotient singularities, it follows that $X$ is a $2-$dimensional orbifold. Conversely, a $2-$dimensional compact K\"ahler holomorphic symplectic orbifold is a normal, compact K\"ahler surface with rational Gorenstein singularities that has a holomorphic symplectic form on its smooth locus: by Theorem 6 of \cite{Nam} we get that $X$ is symplectic.} 
\end{remark}

The first result we will prove is the following, which provides a characterization of the surfaces we are willing to classify:

\begin{thm}
\label{thm:mainsurfaces}
Let $X$ be a normal compact K\"ahler surface whose smooth locus carries a holomorphic symplectic form.
\begin{enumerate}
    \item The surface $X$ is a symplectic variety if and only if it is primitive symplectic if and only if it is the contraction of an ADE configuration $B$ of rational curves on a K3 surface $S$.
    \item The surface $X$ is irreducible symplectic if and only if it is the contraction of an ADE configuration $B$ of rational curves on a K3 surface $S$, where $B$ is not a configuration of the list appearing in Theorem \ref{thm:fujiki}.
    \item The surface $X$ is simple symplectic if and only if it is the contraction of an ADE configuration $B$ of rational curves on a K3 surface $S$, where $B$ has a primitive embedding in $NS(S)$ and is not a configuration in the list appearing in Theorem \ref{thm:fujiki}, nor $B_{\mathfrak{A}_{5}},B_{\mathfrak{A}_{6}},B_{L_{2}(7)},B_{M_{20}}$ 
    (see Remark \ref{rem: subset S} for the definitions).
\end{enumerate}
\end{thm}

The proof of Theorem \ref{thm:mainsurfaces} is mainly contained in Sections \ref{covering of singular K3} and \ref{sec: Covers of K3}. More precisely Corollaries  \ref{cor:singprim} and \ref{cor:orbif} prove the first point, Theorem \ref{thm:mainfujiki} and Proposition \ref{prop: abelian covers} prove the second and Theorem \ref{thm:mainsimple} and its lattice theoretic version, Corollary \ref{cor: no lattices in S}, proves the last.

Once this characterization is obtained, we see that in order to get a complete classification of primitive symplectic, irreducible symplectic and simple symplectic surfaces we need to get a complete list of the ADE configurations that may be realized as ADE configurations of rational curves on K3 surfaces. This is the content of Sections \ref{sec: Covers of K3} and \ref{sec_ADEconfig}, and the summary of the results we get is the following: 
\begin{thm}
\label{thm:classification}
\begin{enumerate}
    \item There are exactly 5836 different ADE configurations on K3 surfaces whose contraction is a singular primitively symplectic surface.
    \item There are exactly 5826 different ADE configurations on K3 surfaces whose contraction is a singular irreducible symplectic surface.
    \item There are exactly 4697 different ADE configurations on K3 surfaces whose contraction is a singular simple symplectic surface.
    \item There are exactly 81 different ADE configurations on K3 surfaces whose contraction is a singular irreducible symplectic surface that admits a smooth finite quasi-\'etale cover.
 \end{enumerate}
\end{thm}
As contractions of different ADE configurations of rational curves on a K3 surface provide surfaces with different singularities, the previous result provides 5826 different locally trivial deformation classes of irreducible symplectic surfaces.

Theorem \ref{thm:classification} is a combination of some results that have been known since several years, see \cite{Fuj}, \cite{Xi} and \cite{Sh}, here reinterpreted in view of the definitions of singular symplectic surfaces we provided before. We add the complete computations of the ADE configurations and describe their geometric properties with respect to the existence of quasi-\'etale covers. 

The proof of Theorem \ref{thm:classification} is presented in Section \ref{subsec: proof} and is based on the content of Sections \ref{sec: Covers of K3} and \ref{sec_ADEconfig}.  
Besides giving explicit proofs of all the results, we will provide a computer program that allows to write explicitly the list of all the possible  ADE configurations of rational curves on K3 surfaces: this is described explicitly in Appendix \ref{sec: the algortihm GAP}.

Some of the ADE configurations mentioned in Theorem \ref{thm:classification} correspond to more than one family of singular symplectic surfaces. In Theorem \ref{theo: A1 singularities} we give a more precise classification of the irreducible symplectic surfaces whose singularities are all of type $A_1$: there are 15 possible ADE configurations, which give 20 different families; among them 11 corresponds to simple symplectic surfaces.

\begin{remark}{\rm  In dimension 2 smooth Calabi--Yau manifolds and irreducible symplectic manifolds coincide, and they are exactly the K3 surfaces. Definitions \ref{def: ISO} and \ref{def: THE DEFINTION} have an analogue which correspond to irreducible Calabi--Yau orbifolds and to irreducible Calabi--Yau varieties respectively, see \cite{Per} for precise definitions. Observe that even in the singular setting the definition of irreducible symplectic variety (resp. orbifold) coincide with the one of  irreducible Calabi--Yau variety (resp. orbifold). Hence the classification present in this paper is also a classification of the singular analogue of Calabi--Yau manifolds in dimension 2.}\end{remark}

The last part of the paper is about the construction of higher dimensional examples of singular symplectic varieties from singular symplectic surfaces. A classical construction of irreducible symplectic manifolds of higher dimension is given by the Hilbert scheme of $n$ points on K3 surfaces. It is natural to ask if the Hilbert scheme of $n$ points on an irreducible (resp. primitive, simple) symplectic surface is an irreducible symplectic variety (resp. primitive symplectic variety, irreducible symplectic orbifold) as well.

The results we prove in this direction are the following (see Proposition \ref{prop:genhilb2} and Corollary \ref{cor:simpleisv}): 

\begin{thm}
\label{thm:hilb2thm}
Let $X$ be a projective primitive symplectic surface.
\begin{enumerate}
    \item The Hilbert scheme $Hilb^{n}(X)$ is a $(2n)-$dimensional 
    primitive symplectic variety.
    \item If $X$ is simple, then $Hilb^{2}(X)$ is an irreducible symplectic 4-dimensional orbifold.
\end{enumerate}
\end{thm}

In the recent paper \cite[Theorem A]{LBP}, it is proved that if $X$ is a projective irreducible symplectic surface, then $Hilb^2(X)$ is an irreducible symplectic variety. 
Finally, we will calculate the second rational cohomology of all these examples: this is a very important locally trivial deformation invariant of a primitive symplectic variety by the Global Torelli Theorem proved in \cite{BL}. 

By a result of Fu and Menet \cite{FuMen} we know that if $X$ is a primitive symplectic variety of dimension 4 which is an orbifold with terminal singularities, then $3\leq b_{2}(X)\leq 23$. In \cite{FuMen}, \cite{Men18}, \cite{Men14}, \cite{Men22} and \cite{BGMM} there are examples of irreducible symplectic orbifolds of dimension 4 with terminal singularities and $b_{2}\in\{3,5,6,8,9,10,11,14,16,23\}$.

By calculating the second rational cohomology group of the irreducible symplectic orbifolds of Theorem \ref{thm:hilb2thm} we prove the following theorem, which stated that all the possible values of $b_2$ are realized if one removes the condition that all the singularities are terminal. It remains open the problem of realizing all the possible values of the $b_2(X)$ if one requires that the all the singularities are terminal (in particular of codimension at least 4).
\begin{thm}
\label{thm:b2qualunque}
For every integer $3\leq n\leq 23$ there exists an irreducible symplectic orbifold $X$ of dimension 4 with $b_2(X)=n$.
\end{thm}


\begin{center}{\bf Acknowledgments}\end{center}

The authors wish to thank F. Catanese, S. Kondo, M. Loenne, V. Bertini, L. Li Bassi and E. Romano for fruitful discussions, M. Mauri for pointing out a mistake in the previous formulation of Proposition 2.1 and for showing us a proof of the existence of a symplectic resolution of the Hilbert scheme of points on a symplectic surface, and A. Craw for his interesting comments and for bringing to our attention the reference \cite{CY}. Moreover, they are thankful to the anonymous referees for their careful reading of
the manuscript. The authors were partially supported by the Research Project PRIN 2020 - CuRVI, CUP J37G21000000001. Arvid Perego was partially supported by the Research Project PRIN 2022 PEKBY - \textit{Symplectic varieties: their interplay with Fano manifolds and derived categories}, CUP J53D23003840006. Matteo Penegini and Arvid Perego wish to thank the MIUR Excellence Department of Mathematics, University of Genoa, CUP D33C23001110001; Alice Garbagnati wishes to thank the MIUR Excellence Department of Mathematics, University of Milano. The authors are members of the INDAM-GNSAGA.

\section{Primitive symplectic surfaces as contractions of K3 surfaces}

A primitive symplectic surface $X$ is a normal singular surface, so there is a smooth surface $S$ such that $X$ is the contraction of a family of curves on $S$. The condition on $X$ of being a primitive symplectic surface imposes several conditions on the possible smooth minimal surfaces $S$ whose contraction is $X$. 

In this section we investigate the relations between these conditions and in particular we prove 
that $S$ is necessarily a K3 surface. 

We will use the following notation: if $X$ is a normal variety and $p\in\mathbb{N}$, we let $$h^{[p],0}(X):=\dim H^{0}(X,\Omega_{X}^{[p]}).$$

We recall that a resolution of the singularities of $X$ is a map $f:\tilde{X}\ra X$ from a smooth surface $\tilde{X}$ to $X$ which is an isomorphisms outside the singularities of $X$. It is called \textit{minimal} if for every resolution $g:\tilde{X}'\ra X$ there is a morphism $\psi:\tilde{X}'\ra \tilde{X}$ such that $g=f\circ \psi$. It follows that if $f:\tilde{X}\ra X$ is a minimal resolution and $p\in Sing(X)$, then $f^{-1}(p)$ does not contain any $(-1)$-curve.

\begin{prop}
\label{prop:contraction}
Let $X$ be normal compact complex surface with at worst rational singularities whose smooth locus carries a holomorphic symplectic form and let $S$ be its minimal model. Then:\begin{enumerate}
\item $S$ is either a K3 surface, a torus, or a primary Kodaira surface
\item if $X$ is K\"ahler, then $S$ is a K3 surface or a torus;
\item if $X$ is singular, then $S$ is a K3 surface; in particular $X$ and $S$ are K\"ahler.\end{enumerate}
\end{prop}

\proof Let $f:\tilde{X}\longrightarrow X$ be a minimal resolution of $X$ and denote $E=\bigcup_{i=1}^{n}E_i$ the exceptional locus of $f$, where $E_1,\cdots,E_{n}$ are the irreducible components of $E$. 

The smooth locus of $X$ is isomorphic to $\tilde{X}-E$, hence $\tilde{X}-E$ is endowed with a holomorphic symplectic form. It follows that the canonical bundle of $K_{\tilde{X}}$ restricted to $\tilde{X}-E$ is trivial, so that \begin{equation}\label{eq: KXtilde}K_{\tilde{X}}=\sum a_i E_i, \ \ \ a_i\geq 0.\end{equation}

Since $E$ is a contractible set, the intersection form on the lattice $\langle E_1,\cdots,E_{n}\rangle$ is negative definite. As a consequence we get $$K_{\tilde{X}}^2\leq 0\mbox{ and }K_{\tilde{X}}^2=0\mbox{ if and only if }a_i=0\ \forall i.$$In particular $K_{\tilde{X}}^2=0$ if and only if $K_{\tilde{X}}=0$.

By definition of minimal resolution we see that $\tilde{X}$ contains no $(-1)-$curves. Indeed, suppose that $R$ is a $(-1)-$curve on $\tilde{X}$. By the genus formula we get that $RK_{\tilde{X}}=-1$, so $R\sum a_iE_i=-1$ with $a_i\geq 0$ and $E_i$ irreducible curves. Hence there is $1\leq i\leq n$ such that $R=E_{i}$, which contradicts the hypothesis that $\tilde{X}$ is a minimal resolution. 

It follows that $\tilde{X}$ is a minimal surface, and hence $\tilde{X}=S$. In particular, we get that $K_{S}=K_{\tilde{X}}$. If  $K_S\neq 0$,  $S$ would be is a minimal surface such that $K_S^2<0$, hence $\kappa(S)=-\infty$. By \eqref{eq: KXtilde}, $K_{\widetilde{X}}=K_S$ is effective, which gives a contradiction. So $K_S= 0$ and hence $S$ is either a K3 surface or a complex torus or a primary Kodaira surface (see \cite[Table 10]{BHPV}).

If $X$ is  K\"ahler, then $S$ is K\"ahler, and this excludes the primary Kodaira surfaces.

If $X$ is singular, it is obtained as a contraction of rational curves on $S$, but complex tori and primary Kodaira surfaces do not contain rational curves.\endproof

\begin{remark}
{\rm The hypothesis of Proposition \ref{prop:contraction} are in particular satisfied by singular symplectic surfaces (see Remark \ref{rem:namikawa}). Moreover, by Remark \ref{remark:orbifold}, point (1) of Proposition \ref{prop:contraction} is equivalent to Lemma 6.6 in \cite{Fuj}.}
\end{remark}

\begin{corollary}
\label{cor:singprim}
A singular compact complex surface $X$ is a symplectic surface if and only if it is a primitive symplectic surface.
\end{corollary}

\proof
One implication is trivial in any dimension both for smooth and singular varieties. The other one depends on the classification of surfaces: if $X$ is a singular symplectic surface, it satisfies the assumption of Proposition \ref{prop:contraction}, case (3), so its minimal model $S$ is a K3 surface. By the Leray spectral sequence we have that $h^{1}(X,\mathcal{O}_{X})\simeq h^{1}(S,\mathcal{O}_{S})$, and by \cite[Theorem 1.5 and Remark 1.5.2]{GKKP} we have that $$h^{[2],0}(X)=h^{2,0}(S)=1,$$hence $H^{0}(X^{s},\Omega_{X^{s}}^{2})$ is spanned by the class of a holomorphic symplectic form on $X^{s}$. It follows that $X$ is primitive symplectic. 
\endproof

An immediate corollary of Proposition \ref{prop:contraction} is then the following, which is point (1) of Theorem \ref{thm:mainsurfaces}:

\begin{cor}
\label{cor:orbif}
A normal compact 
complex surface is primitive symplectic if and only if $X$ is the contraction of  an ADE configuration of rational curves on a K3 surface.  
\end{cor}

\proof If $X$ is smooth, then it is primitive symplectic if and only if it is a K3 surface, and there is nothing to prove.

If $X$ is a singular primitive symplectic surface, then it is a symplectic variety, hence it has canonical singularities. Since canonical singularities are rational, by Proposition \ref{prop:contraction} we see that the minimal resolution $f:S\longrightarrow X$ of $X$ is a K3 surface. Moreover, since canonical singularities are ADE singularities on surfaces, the resolution $f$ is the contraction of an ADE configuration of curves on $S$ (see \cite[Section 3 Chapter 3]{BHPV} for the definition and the properties of ADE singularities). 

Conversely, if $X$ is the contraction of an ADE configuration of curves on a K3 surface $S$, then $X$ is a normal compact complex surface that carries a holomorphic symplectic form $\sigma$ on its smooth locus $X^{s}$. If $f:S\longrightarrow X$ is the contraction, then we see that $\sigma$ extends to a holomorphic $2-$form on $S$, so by \cite{B2} we see that $X$ is compact complex symplectic surface. By Corollary \ref{cor:singprim} we are done.\endproof


\begin{remark}
{\rm We recall that an irreducible curve on a K3 surface has negative self intersection if and only if it is smooth, rational and its self-intersection is $-2$. Since the lattices spanned by their
roots are the ADE lattices we obtain that all the contractions of a K3 surface are contractions of ADE configurations of smooth rational curves and are in particular canonical.}
\end{remark}

As a consequence of Corollary \ref{cor:orbif}, it follows that in order to classify all possible primitive irreducible symplectic surfaces one needs to classify all ADE configurations of rational curves on K3 surfaces. This will be the aim of Section \ref{sec_ADEconfig}.

\par\bigskip

We now wish to notice that there are primitive symplectic surfaces which are not irreducible symplectic surfaces: the reason is that while conditions (1) and (2) of Definition \ref{def: THE DEFINTION} are always verified by all contractions of ADE configurations of rational curves on K3 surfaces, condition (3) is not trivially satisfied, as the following well known example shows.

\begin{esem}
{\rm Let $T$ be a $2-$dimensional complex torus and $\iota:T\longrightarrow T$ the involution mapping $p\in T$ to $-p$. Let $X:=T/\iota$, so that $X$ is a singular symplectic surface having 16 singular points, and the resolution of the singularities $f:S\longrightarrow X$ is such that $S$ is a K3 surface, called Kummer surface of $T$. The morphism $f$ is a contraction of 16 disjoint $(-2)-$curves, but the surface $X$ is not irreducible symplectic: the quotient morphism $q:T\longrightarrow X$ is a finite quasi-\'etale covering, but since $H^{0}(T,\Omega_{T}^{1})\neq 0$ we see that the algebra of holomorphic forms on $T$ is not spanned by the pull-back of the symplectic form on $X$. 
}
\end{esem} 

Therefore, in order to classify the irreducible symplectic surfaces one needs not only to classify all possible ADE configurations of rational curves on K3 surfaces, but also to give necessary and sufficient conditions that guarantee that their contraction satisfies condition (3) of \ref{def: THE DEFINTION}.

We will use the following notation: if $S$ is a K3 surface and $B$ is an ADE configuration of $(-2)-$curves on $S$, we let $f_{B}:S\longrightarrow X_{B}$ be the contraction morphism of the curves in $B$, and we let $X_{B}$ be the singular surface obtained by this contraction.

In order to understand under which conditions the surface $X_{B}$ is an irreducible symplectic surface, we need to consider all finite quasi-\'etale coverings $Y\longrightarrow X_{B}$, and to calculate the dimension of $h^{0}(Y,\Omega_{Y}^{[p]})$ for $p=0,1,2$. To do so we first relate the finite quasi-\'etale coverings of $X_{B}$ with some induced coverings of the K3 surface $S$, as follows:

\begin{prop}
\label{prop:qetfp}
Let $S$ be a K3 surface, $B$ an ADE configuration of $(-2)-$curves on $S$, and $f_{B}:S\longrightarrow X_{B}$ the contraction of the curves in $B$. Let $g:Y\longrightarrow X_{B}$ be a finite quasi-\'etale covering of degree $n$, let $Z:=S\times_{X_{B}}Y$ and $\widetilde{g}:Z\longrightarrow S$ be the natural projection. Then $\widetilde{g}$ is a finite covering of degree $n$ whose branch locus is contained in the exceptional locus $B$ of $f_{B}$.
\end{prop}

\proof Let $D\subseteq X_B$ be the branch locus of $g$ and $U:=X_B\setminus D$. As $g$ is quasi-\'etale, we see that $D$ is given by a finite number of points. Moreover, we let $Y_{U}:=g^{-1}(U)$, which is an open subset of $Y$, and $S_{U}:=f_{B}^{-1}(U)$, which is an open subset of $S$. 

Notice that $g_{|Y_{U}}:Y_{U}\longrightarrow U$ is \'etale of degree $n$, and if we let $Z_{U}:=S_{U}\times_{U}Y_{U}$ (which is an open subset of $Z$), then the projection $\widetilde{g}_{|Z_{U}}:Z_{U}\longrightarrow S_{U}$ is \'etale and has the same degree of $g_{|Y_{U}}$. It follows that $\widetilde{g}:Z\longrightarrow S$ is generically finite of degree $n$.

We are only left with showing that $D\subseteq f_B(B)$, where we notice that $f_B(B)$ is the singular locus of $X_B$. To prove this, suppose that $p\in D$ is a smooth point of $X_{B}$. Hence $f_{B}^{-1}(p)$ is a single point of $S$ which is contained in the branch locus of $\widetilde{g}$. But since $S$ is smooth and $Z$ is normal, by purity the branch locus of $\widetilde{g}$ has no isolated points, and we get a contradiction.\endproof

The previous result tells us that any finite quasi-\'etale covering $g:Y\longrightarrow X_{B}$ induces by base change a generically finite covering (of the same degree) $\widetilde{g}:Z\longrightarrow S$ of the K3 surface $S$, whose branch locus is contained in $B$. Next section will be aimed to understand the geometry of $Z$ and of the branch locus of $\widetilde{g}$.

\section{Coverings of singular symplectic surfaces}\label{covering of singular K3}

We described the primitive symplectic surfaces in the previous section. To identify the ones which are also irreducible symplectic or even simple symplectic, we have to consider their finite quasi \'etale covers. Therefore, in this section we consider these covers and in particular we show that we can restrict the study to Galois covers with out loss of generalities.

Moreover, as finite quasi \'etale covers of a symplectic surface $X$ are naturally associated to branched covers of its resolution $S$ as discussed in the previous section, it will be useful to recall some results on branched covers of K3 surfaces.

\subsection{Galois and non Galois covers}\label{subsec:Galois and non Galois}

Let $\pi:Y\longrightarrow X$ be a covering of $X$, which can be Galois or not.

Suppose first that $\pi$ is a Galois covering branched on a subset of $Sing(X)$.
We recall that all the Galois covers are assumed to be normal. 
Since ADE singularities on surfaces are quotient singularities, if $X$ is a symplectic surface, then it is an orbifold. So now we consider Galois covers of orbifolds. 

Let $X$ be an $n$-dimensional orbifold with only isolated singularities, such that every singular point $p$ of $X$ has an open neighborhood isomorphic to a neighborhood of the origin in the quotient $\mathbb{C}^n/G_{p}$ for some subgroup $G_{p}$ of $SL(n,\mathbb{C})$.

\begin{lemma}
\label{lem:orbicov}
Let $X$ be an $n$-dimensional orbifold with only isolated singularities and $\pi:Y\longrightarrow X$ a Galois covering branched only on a subset of $Sing(X)$. Then $Y$ is an $n-$dimensional orbifold with only isolated singularities, and if $q$ is a singular point of $Y$ then $q$ has an open neighborhood which is isomorphic to a neighborhood of the origin in the quotient $\mathbb{C}^n/H_{q}$ where $H_{q}$ is a normal subgroup of $G_{\pi(q)}$. 
\end{lemma} 

\proof Consider the restriction $\pi':Y\setminus\pi^{-1}(Sing(X))\rightarrow  X\setminus Sing(X)$ of $\pi$. Then $\pi'$ is an \'etale Galois covering of the manifold $X\setminus Sing(X)$, so that $Y\setminus\pi^{-1}(Sing(X))$ is smooth. 

Let us analyse locally what happens around a singular point $p\in Sing(X)$. The point $p$ has an open neighborhood $U_p$ which is isomorphic to $\mathbb{C}^n/G_{p}$, for some subgroup $G_{p}$ of $SL(n,\mathbb{C})$. The universal covering of $U_p$ is $\mathbb{C}^n$, so each Galois covering of $U_p$ is a quotient of $\mathbb{C}^n$ by a normal subgroup $H_{p}$ of $G_{p}$. Each connected component of $\pi^{-1}(U_p)$ is of the form $\mathbb{C}^n/H_{p}$, hence each point $q$ such that $q\in\pi^{-1}(p)$ has an open neighborhood which is isomorphic to $\mathbb{C}^{n}/H_{p}$.\endproof

Suppose now that $\pi:Y\ra X$ is not Galois. As before $Y\setminus\pi^{-1}(Sing(X))$ is smooth, and we still look at what happens around a singular point $p\in Sing(X)$. Let $U_{p}$ be an open neighborhood of $p$ in $X$ which is isomorphic to $\mathbb{C}^{n}/G_{p}$. Then we have a commutative diagram
$$\xymatrix{&\mathbb{C}^{n}\ar[dl]^{/H_{p}}\ar[dr]_{/G_{p}}\\V_p\ar[rr]_{\pi}&&U_p}$$
where $V_p$ is one of the connected components of $\pi^{-1}(U_p)$: in this case the group $H_{p}$ is still a subgroup of $G_{p}$, but it is non necessarily normal. 

In the following will be useful to compare a non Galois covering with its Galois closure. Hence we recall the definition and well known properties of the Galois closure.

\begin{defn}
Let $\pi:Y\ra X$ be a finite covering and $L$ the Galois closure of the function fields extension $k(X)\subset k(Y)$. The Galois closure of $\pi$ is given by $Z$, which is the normalization of $Y$ in $L$, and the induced maps $\alpha:Z\ra X$, $\beta:Z\ra Y$. 
\end{defn}

Notice that $\alpha$ and $\beta$ are Galois coverings and that $\alpha$ does not factorize into further Galois covers by the minimality of the extension $L$.

\begin{lemma}\label{lem: Galois closure}
Let $\pi:Y\ra X$ be a finite covering with branch locus $D$ and $Y$ a normal variety, and let $\alpha:Z\ra X$ be its Galois closure. Then the branch locus of $\alpha$ is $D$. 
\end{lemma} 

\proof The Galois closure of $\pi$ is such that  $\alpha:Z\ra X$ is a Galois $G$-covering such that $\alpha=\beta\circ\pi$ for $\beta:Z\ra Y\simeq Z/H$ for a certain subgroup $H\leq G$. So we have the following diagram
$${\xymatrix{Z\ar[dr]^{\beta}\ar[ddr]_{\alpha}\\
		&Y=Z/H\ar[d]^{\pi}\\
		&X=Z/G.}}$$

If $p\in X$ is a point in the branch locus of $\pi$, then its inverse image under $\alpha=\pi\circ \beta$ has cardinality strictly less than $|G|$, therefore it is in the branch locus of $\alpha$.

If $q\in X$ is a point in the branch locus of $\alpha$ but not in that of $\pi$, then each point in $\pi^{-1}(q)$ is in the branch locus of $\beta$. Since $\beta$ is a Galois cover, there exists a subgroup $K$ of $G$ that stabilizes each point in $\beta^{-1}(\pi^{-1}(q))=\alpha^{-1}(q)$. The action of $K$ on $\alpha^{-1}(q)\subset Z$ is then the identity, so every $k\in K$ commutes with any element in $G$: this implies that $K$ is a normal subgroup of $G$. It follows that $\alpha$ factorizes through the quotient $Z/K$, contradicting the minimality of the Galois closure.\endproof 

As a consequence, whenever we have a covering of an orbifold with isolated singularities, we may always suppose that it is a Galois covering without changing the branch locus.

\subsection{Generalities on covers of K3 surfaces}

Let us look at what happens in the case of orbifolds of dimension 2 with ADE singularities. 

\begin{prop}
\label{prop: Kodaira dimension ADE covers} 
Let $X$ be a surface with only ADE singularities and $\pi:Y\rightarrow  X$ a finite covering whose branch locus is contained in $Sing(X)$. Then $Y$ has only ADE singularities and $\kappa(X)=\kappa(Y)$.
\end{prop}

\proof Since $X$ is an orbifold, a singular point $p$ of $X$ is such that there exists an open neighborhood isomorphic to $\mathbb{C}^2/G_{p}$ where $G_{p}$ is a finite subgroup of $GL(2,\mathbb{C})$. The singularity in $p$ is an ADE singularity if and only if $G_{p}\subset SL(2,\mathbb{C})$. 

By Lemma \ref{lem:orbicov} the surface $Y$ is a $2-$dimensional orbifold with only isolated singularities, and if $q$ is a singular point of $Y$, then it has an open neighborhood isomorphic to a neighborhood of the origin in the quotient $\mathbb{C}^{2}/H_{q}$ where $H_{q}$ is a subgroup of $G_{\pi(q)}$: it is then a finite subgroup of $SL(2,\mathbb{C})$, and hence $Y$ is a surface with only ADE singularities. 

Finally, $\pi:Y\ra X$ is a quasi-\'etale between surfaces with canonical singularities, so the Kodaira dimension of $Y$ equals the one of $X$ by \cite[Section 3, Page 51]{Cat}.\endproof

As a consequence we see that finite quasi-\'etale coverings of contractions of ADE configurations of $(-2)-$curves on K3 surfaces may only be of very special type, i.e., their minimal models can only be either K3 surfaces or $2-$dimensional complex tori. More precisely, we have the following:

\begin{cor}
\label{cor: canonical ADE covers}
Let $S$ be a K3 surface, $B$ an ADE configuration of $(-2)-$curves on $S$ and $f_{B}:S\rightarrow X_{B}$ be the contraction of the curves in $B$. Then: 
\begin{enumerate}
    \item the surface $X_B$ is irreducible symplectic if and only if it does not admit any finite quasi-\'etale covering whose minimal model is a 2-dimensional complex torus;
    \item the surface $X_B$ is simple symplectic if and only if it does not admit any non-trivial finite quasi-\'etale covering at all.
\end{enumerate}
\end{cor}


\proof By construction, $X_{B}$ is a primitive symplectic variety and satisfies conditions (1) and (2) of Definition \ref{def: THE DEFINTION}. Moreover, it has only ADE singularities and $\kappa(X_{B})=0$. 

Let $\pi:Y\longrightarrow X_{B}$ be a finite quasi-\'etale cover of $X_{B}$. Then, by Proposition \ref{prop: Kodaira dimension ADE covers}, $Y$ is a possibly singular compact complex surface with only ADE singularities and $\kappa(Y)=0$. By Proposition \ref{prop:contraction} the minimal model of $X_B$ is a K3 surface $S$ and by Proposition \ref{prop:qetfp} the finite quasi-\'etale cover $\pi$ induces a branched cover of $S$, $\pi_S:Y_S\ra S$ whose branch locus is contained in $B$. 

By construction $Y_S$ is birational to $Y$, so let us denote $Y'$ the minimal model of both $Y$ and $Y_S$. In particular $\kappa(Y')=\kappa(Y_{S})=\kappa(Y)=0$ and by \cite[Theorem 1.5 and Remark 1.5.2]{GKKP} we have that $p_{g}(Y_S)=p_g(Y')=h^{[2],0}(Y)$. Since $Y_S$ is a cover of $S$ we have $p_{g}(Y_S)\geq p_g(S)=1$, which, together with $\kappa(Y_S)=0$ implies $h^{[2],0}(Y)=p_g(Y')=p_g(Y_S)=1$. 

By the Enriques-Kodaira classification of compact K\"ahler surfaces, it follows that $Y'$ is either a torus or a K3 surface. If $Y'$ is a torus, then $h^{[1],0}(Y)=h^{1,0}(Y')=h^{1,0}(Y_S)=2\neq 0$; if $Y'$ is a K3 surface, then $h^{[1],0}(Y)=h^{1,0}(Y')=h^{1,0}(Y_S)=0$. It follows that if $Y'$ is a torus, then $h^{[1],0}(Y)=2$ and so the finite quasi-\'etale cover $\pi$ does not satisfy the point (3) of Definition \ref{def: THE DEFINTION}: hence $X_B$ is not irreducible symplectic. 

On the other hand, if there are no covers $Y$ of $X$ such that $Y'$ is a torus, then for every quasi-\'etale finite cover $\pi:Y\longrightarrow X_{B}$ we have that $Y'$ is a K3 surface and hence $h^{[1],0}(Y)=0$. As a consequence $X_B$ is an irreducible symplectic surface because all the finite quasi-\'etale covers satisfy the point (3) of Definition \ref{def: THE DEFINTION}.

As a cover of a smooth surface is either \'etale or branched on a codimension 1 subvariety, it follows that $X_B$ admits a non-trivial finite quasi-\'etale cover if and only if its smooth locus admits a non trivial finite \'etale cover, and hence if and only if it is not simply connected. In particular, $X_B$ is simple if and only if it does not admit any non-trivial finite quasi-\'etale cover at all.
\endproof

By Corollary \ref{cor: canonical ADE covers}, in order to determine if $X_B$ is an irreducible symplectic variety and if it is simple, one has to determine if it admits a finite quasi-\'etale covering whose minimal model is a $2-$dimensional complex torus or a K3 surface, which is the topic of the next section. 

\section{Covers of K3 surfaces and ADE configurations}\label{sec: Covers of K3}
The aim of this Section is to describe in detail the covers of K3 surfaces branched on ADE configurations of rational curves, since these correspond to finite quasi \'etale covers of singular symplectic surfaces. The results are summarized in Theorem \ref{theorem: recap section 4}.

In Section \ref{subsec: covers with tori}, we determine the ADE configurations which are the branch locus of a cover of a K3 surface with a surface birational to a torus, see Theorem \ref{thm:mainfujiki}. The contraction of these configurations are symplectic surfaces, but not irreducible symplectic surface. 

In Section \ref{subsec: covers with K3 surfaces}, we consider the analogous problem for covers with surfaces birational to K3 surfaces, see Proposition \ref{prop: MG}. To this purpose, we have to consider effective divisible classes on K3 surfaces, see Definition \ref{def: divisible}, and describe their shapes, see Proposition \ref{prop: divisible}.

Finally, in Section \ref{subsec: simple symplectic surfaces}, we identify the ADE configurations whose contraction produce simple symplectic surfaces, i.e. which cannot appear as branch locus of any covers, see Theorem \ref{thm:mainsimple}.\\

Before all this, we fix some preliminary notation and results about covers of K3 surfaces.

An ADE configuration of rational curves on a K3 surface $S$ corresponds to a lattice, namely the lattice spanned by the classes of these curves in the N\'eron-Severi lattice of $S$. If a given ADE configuration $B$ exists on a K3 surface $S$, then the associated lattice $\Lambda$ can be embedded (not necessarily primitively) in the N\'eron--Severi of $S$, since the classes of the $(-2)$-curves are the irreducible roots of $\Lambda$ and are contained in $NS(S)$. The converse is also true, as we will see in Lemma \ref{lemma: primitive embdding of Lambda and existence K3} and in Theorem \ref{theorem: existence of ADE configuration}: it will follows that a given ADE configuration $B$ exists on a K3 surface $S$ if and only if the associated lattice $\Lambda$ can be embedded (not necessarily primitively) in the N\'eron--Severi of $S$. 

Let $\Lambda=\oplus_{i=1}^\alpha A_{a_i}\bigoplus\oplus_{j=1}^\delta D_{d_j}\bigoplus \oplus_{h=1}^\epsilon E_{e_h}$ be the lattice associated to an ADE configuration. The lattice $\Lambda$ is a negative definite even lattice of rank $$\rk(\Lambda)=\sum_{i=1}^{\alpha}a_i+\sum_{j=1}^{\delta}d_j+\sum_{h=1}^{\epsilon}e_h.$$ 
\begin{lemma}\label{lemma: primitive embdding of Lambda and existence K3} 
	Let $\Lambda:=\oplus_{i=1}^\alpha A_{a_i}\bigoplus\oplus_{j=1}^\delta D_{d_j}\bigoplus \oplus_{h=1}^\epsilon E_{e_h}$.
	\begin{enumerate}
		\item If there exists a K3 surface $S$ on which there is the ADE configuration of rational curves of type $\oplus_{i=1}^\alpha A_{a_i}\bigoplus\oplus_{j=1}^\delta D_{d_j}\bigoplus \oplus_{h=1}^\epsilon E_{e_h}$, then the lattice $\Lambda$ is embedded in $NS(S)$ and hence in $\Lambda_{K3}$.
		\item If the lattice $\Lambda$ is primitively embedded in $\Lambda_{K3}$, then there exists a K3 surface $S$ with an ADE configuration $\oplus_{i=1}^\alpha A_{a_i}\bigoplus\oplus_{j=1}^\delta D_{d_j}\bigoplus \oplus_{h=1}^\epsilon E_{e_h}$ of smooth rational curves.
	\end{enumerate}
\end{lemma}

\proof If $S$ is a K3 surface on which there exists the required configuration, each curve of this configuration is an algebraic class on $X$. Hence the lattice spanned by these curves is contained in $NS(S)$. 

Conversely, if $\Lambda$ can be primitively embedded in $\Lambda_{K3}$ and $\Lambda$ is a negative definite lattice, by the surjectivity of the period map there exists a K3 surface $S$ (indeed a family of $\Lambda$-polarized K3 surfaces) such that $NS(S)\simeq \Lambda$, and $S$ is a non--projective K3 surface. 

By the Riemann--Roch Theorem, each class with self intersection $-2$ contained in the N\'eron--Severi group of a K3 surface is either effective, or the opposite of an effective class. So, up to a sign, we assume that the generators of $\Lambda$ are effective and represent a possibly reducible curve. We recall that $\Lambda$ is negative definite, so each effective $(-2)$-class in $\Lambda$ is either the class of an irreducible curve or the sum of $(-2)$-classes with prescribed intersection properties. 

Since the ADE lattices are generated by their irreducible roots and $NS(S)\simeq \Lambda$, the irreducible roots of $\Lambda$ correspond (up to a sign) to irreducible effective class on the K3 surface $S$. In particular, each irreducible root is (up to a sign) the class of a smooth irreducible rational curve on $S$, i.e. on $S$ there is a configuration of $(-2)-$curves whose associated lattice is $\Lambda$. 
\endproof

\subsection{Coverings with tori}\label{subsec: covers with tori}

By Proposition \ref{prop:qetfp} and Corollary \ref{cor: canonical ADE covers}, a contraction $X_{B}$ of an ADE configuration $B$ of $(-2)-$curves on a K3 surface $S$ is an irreducible symplectic surface if and only if $S$ has no generically finite covering whose branch locus is contained in $B$ and which is given by a surface birational to a $2-$dimensional complex torus.

The classification of the branch loci of the coverings of K3 surfaces which are birational to $2-$dimensional complex tori may be obtained as a consequence of a result of Fujiki: more precisely, in \cite{Fuj2} the groups $G$ of automorphisms of a 2-dimensional complex torus $A$ such that $A/G$ is birational to a K3 surface $\widetilde{A/G}$ are classified. 

In particular, in \cite[Lemma 3.19]{Fuj2} the points of $A$ with non trivial stabilizer are classified. Some of them are identified by the action of the group. In \cite{NKummer}, \cite{Bertin}, \cite{W}, \cite{G}, \cite{R} the configuration of curves in the branch locus of $A\rightarrow \widetilde{A/G}$ and the attached lattices is given for all the groups (the order two case is studied in \cite{NKummer}, the other cyclic cases in \cite{Bertin}, the remaining ones in \cite{W}, \cite{G} and \cite{R}). These known results are summarized in the following theorem.

\begin{thm}{\rm(\cite{Fuj2}, \cite{NKummer}, \cite{Bertin}, \cite{W}, \cite{G}, \cite{R})}
\label{thm:fujiki} Let $f:A\longrightarrow S$ be a generically finite covering of a K3 surface $S$ where $A$ is a surface whose minimal model is a torus, and let $B$ its branch locus. Then the ADE configuration type of $B$ is one of the following:
\begin{enumerate}
    \item $B_{1}:=A_{1}^{\oplus 16}$;
    \item $B_{2}:=A_{2}^{\oplus 9}$;
    \item $B_{3}:=A_3^{\oplus 4}\oplus A_1^{\oplus 6}$;
    \item $B_{4}:=A_5\oplus A_2^{\oplus 4}\oplus A_1^{\oplus 5}$;
    \item $B_{5}:=D_4^{\oplus 2}\oplus A_3^{\oplus 3}\oplus A_1^{\oplus 2}$;
    \item $B_{6}:=D_4^{\oplus 4}\oplus A_1^{\oplus 3}$;
    \item $B_{7}:=D_5\oplus A_3^{\oplus 3}\oplus A_2^{\oplus 2}\oplus A_1$;
    \item $B_{8}:=E_6\oplus D_4\oplus A_2^{\oplus 4}\oplus A_1$;
    \item $B_{9}:=A_3^{\oplus 6} \oplus  A_1$;
    \item $B_{10}:=A_5 \oplus  A_3^{\oplus 2} \oplus A_2^{\oplus 4}$.
\end{enumerate}
\end{thm}

This allows us to state and prove the following, which is point (2) of Theorem \ref{thm:mainsurfaces}:

\begin{thm}
\label{thm:mainfujiki}
Let $S$ be a K3 surface, $B$ an ADE configuration of smooth rational curves on $S$ and $f_{B}:S\longrightarrow X_{B}$ the contraction morphism. Then $X_{B}$ is an irreducible symplectic surface if and only if $B$ does not contain any of the ADE configurations of Theorem \ref{thm:fujiki}.
\end{thm}

\proof We preliminary observe that $X_B$ admits a finite quasi-\'etale cover with a surface birational to  torus if and only if $S$ admits a cover (branched on a subconfiguration $B'\subseteq B$) with a surface whose minimal model is a torus. Moreover, $S$ admits a cover with surface whose minimal model is a torus if and only if it admits a Galois cover of this type. Indeed, if $\pi:Y\ra S$ is a  cover branched on an ADE configuration such that $h^{1,0}(Y)=2$, its Galois closure $\alpha:Z\ra S$ is a cover branched on the same ADE configuration (see Lemma \ref{lem: Galois closure}). It follows that $h^{1,0}(Z)$ is either 0 or 2. Since the Galois closure factors through the cover, the surface $Z$ is a cover of the surface $Y$, which implies that $h^{1,0}(Z)\geq h^{1,0}(Y)$. So $h^{1,0}(Z)=2$ and $Z$ is birational to a torus. 

Hence, by Corollary \ref{cor: canonical ADE covers}, $X_B$ is not an irreducible symplectic surface if and only if $S$ admits a Galois cover with a surface whose minimal model is a torus, which implies that $S$ is birational to a quotient of a torus, in particular it is the smooth minimal model of the quotient of a torus. This implies that $S$ contains one of the configuration listed in Theorem \ref{thm:fujiki}. 

Conversely, let us suppose that $S$ contains one of the configurations listed in Theorem \ref{thm:fujiki} and let us denote it as $B$. Then it is known that $B$ is the branch locus of a cover: this is proved in \cite{NKummer} for $G=\Z/2\Z$; in \cite{Bertin} for $G=\Z/n\Z$, $n=3,4,6$; in \cite{G} for $G=D_{12}$ (the binary dihedral group of order 12), $G=T_{24}$ (the binary tetrahedral group  of order 24) and for two possible actions of $G=Q_8$ (the quaternion group); in \cite{R} for the remaining actions of $G=Q_8$ and $G=T_{24}$. So if $S$ contains a configuration $B$ of curves among the ones listed in Theorem \ref{thm:fujiki}, there exists a cover of $S$ branched on such a configuration whose minimal model is a torus. This implies that $X_B$ is not an irreducible symplectic surface.\endproof
\begin{remark}\label{rem: curves implies covers for tori}
{\rm
As briefly recalled in the proof of the previous Theorem, if $B$ is one of the configurations of curves in Theorem \ref{thm:fujiki}, then it suffices that $B$ appears as a configuration of rational curves on a K3 surface $S$ to assure that $B$ is the branch locus of a cover of $S$. This is a very special property of the configurations in Theorem \ref{thm:fujiki}. A priori, the existence of a cover branched on a set of curves, depends also on the embedding of the curves in the surface, as discussed in Example \ref{esem: ADE no imples Galois}. }
\end{remark}
The main consequence of Theorem \ref{thm:mainfujiki} for our purposes is that in order to classify all the irreducible singular surfaces, one has to classify all the possible ADE configurations of $(-2)-$curves on K3 surfaces, and then exclude those containing a subconfiguration of the form listed in Theorem \ref{thm:mainfujiki}. 

\subsection{Covers with K3 surfaces}\label{subsec: covers with K3 surfaces}
Having completely characterized primitive and irreducible symplectic surfaces  in Corollary \ref{cor:orbif} and Theorem \ref{thm:mainfujiki} in terms of contractions of curves on K3 surfaces, we are left with finding a similiar characterization for simple symplectic surfaces, i.e., irreducible symplectic surfaces whose smooth locus is simply connected. 

As simple symplectic surfaces are irreducible symplectic, by Theorem \ref{thm:mainfujiki}, they are obtained by contractions of certain ADE configurations of curves on K3 surfaces that do not contain any of the configurations in the list of Theorem \ref{thm:fujiki}. By Corollary \ref{cor: canonical ADE covers} any finite quasi-\'etale covering of such a contraction has to be birational to a K3 surface.

The possible ADE configurations of curves on K3 surfaces which are the branch locus of a finite covering with another K3 surface have been studied and classified by Xiao in \cite{Xi}. The main result in this direction is the following.


\begin{thm}{\rm (Xiao, \cite{Xi})}
\label{thm:xiao }
Let $f:Y\longrightarrow S$ be a generically finite covering of a K3 surface where $Y$ is a surface birational to K3 surface, and let $B$ its branch locus. Then the ADE configuration type of $B$ is one of the 81 configuration in Table \cite[Table 2]{Xi}.
\end{thm}

The proof of the previous result is contained in \cite{Xi}, where the Galois covers of K3 surfaces with a K3 surface are classified. More precisely, Xiao proves that if $W$ is the quotient of a K3 surface by a finite group, the exceptional locus of the desingularization of $W$ forms an ADE configuration among the 81 possibilities listed in \cite[Table 2]{Xi}. Conversely, Xiao proves that for every ADE configuration $B$ contained in \cite[Table 2]{Xi} there is a K3 surface $S$ having $B$ as ADE configuration of $(-2)-$curves, and that has a Galois covering whose branch locus is $B$.

We notice that each ADE configuration in \cite[Table 2]{Xi} corresponds uniquely to a group $G$ such that there exists a $G$-Galois cover branched on that configuration. For this reason we denote the configuration  in \cite[Table 2]{Xi} as $B_G$ (and with $B_G^i$, $i=1,2$ in the cases $G=Q_8, T_{24}$, since in these cases the same group is attached to more than one ADE configuration).  We observe that if $G$ and $G'$ are two different groups, then $B_G\neq B_{G'}$.

\begin{remark}
{\rm It is important to notice that even if for every ADE configuration $B_G$ in \cite[Table 2]{Xi} there is a K3 surface $S$ (actually, a family of K3 surfaces) for which $B_G$ appears as the branch locus of a finite $G$-Galois covering of $S$, it is not necessarily true that if $W$ is a specific K3 surface which has $B_G$ as an ADE configuration of $(-2)-$curves, then there is a $G$-Galois covering of $W$ whose branch locus is $B_G$. Example \ref{esem: ADE no imples Galois} shows this phenomenon. We observe that the configurations in Theorem \ref{thm:fujiki} behave differently, as discussed in Remark \ref{rem: curves implies covers for tori}.}
\end{remark}

In order to classify the possible covers that a K3 surface has, we need to discuss first the presence of the so called divisible classes, which are related to cyclic covers (in Section \ref{subsub cyclic}). Secondly, we present known results on general Galois cover of K3 surfaces with a surface birational to another K3 surface (in Section \ref{subsub Xiao}).
\subsubsection{Cyclic covers of K3 surfaces and divisible classes}\label{subsub cyclic}
 
\begin{defn}\label{def: divisible}
Let $n\in\mathbb{N}$, $n>1$. An (effective) divisor $D$ is an (effective) $n$-divisible class in the Picard group of a surface $S$ if there exists $L\in Pic(S)$ such that $nL\sim D$.
\end{defn}
	
An effective divisible class is associated to a cyclic non trivial cover of the surface branched on it (see e.g. \cite[Chapter 1, Section 17]{BHPV}). A priori on a K3 surface there could be no effective divisible classes, one effective divisible class or more than one. In the latter case one has more than one cyclic cover of the surface, and possibly a Galois cover whose Galois group is not cyclic. 

Since we are interested in covers branched on ADE configuration of curves, the effective $n$-divisible class we consider are given by linear combination with positive coefficients of $(-2)$-curves.

\begin{esem}
\label{esem: ADE no imples Galois}
{\rm Let us consider the first entry of \cite[Table 2]{Xi}: the group $G$ in this case is $\Z/2\Z$ and the configuration type of $B_{\Z/2\Z}$ is $A_1^{\oplus 8}$, which geometrically corresponds to a set of 8 disjoint rational curves.

 There are two different families of (non projective) K3 surfaces admitting this configuration of curves: the first family, denoted $F_{1}$, is given by the K3 surfaces whose N\'eron-Severi lattice contains primitively $A_1^{\oplus 8}$: the generic member of this family is a K3 surface $S$ such that $NS(S)\simeq A_1^{\oplus 8}$ and the transcendental lattice is $T(S)\simeq U^{\oplus 3}\oplus A_{1}^{\oplus 8}$.

The second family, denoted $F_{2}$, is given by the K3 surfaces whose N\'eron-Severi lattice primitively contains the so-called Nikulin lattice $N$, i.e., a negative definite lattice of rank 8 and length 6 spanned by eight classes $N_i$, $i=1,\ldots, 8$ with $N_iN_j=-2\delta_{i,j}$ and by the class $(\sum_iN_i)/2$. This is an overlattice of index 2 of $A_1^{\oplus 8}$ and it is the lattice obtained requiring that the 8 disjoint rational curves form a 2-divisible set of curves. 
	The generic member of this family is such that $NS(S)\simeq N$ and  $T(S)\simeq U^{\oplus 3}\oplus N$.



It follows that among the two families of K3 surfaces admitting $B_{\Z/2\Z}\simeq A_1^{\oplus 8}$ as an ADE configuration of rational curves, the K3 surfaces which are general member of the family $F_{i}$ admit a double cover with a K3 surface branched on $B_{\Z/2\Z}$ if and only if $i=2$.}


\end{esem}

\begin{remark}
{\rm There is a fundamental difference between the list in \cite[Table 2]{Xi} and the one in Theorem \ref{thm:fujiki}. Indeed in the latter, each ADE configuration is necessarily the branch locus of a Galois cover with a torus (see Remark \ref{rem: curves implies covers for tori}). This is not the case for the list in \cite[Table 2]{Xi}, as the previous Example shows. This is due to the fact that the existence of the curves of an ADE configuration in Theorem \ref{thm:fujiki} forces the minimal lattice spanned by these curves to have some divisibility relations, which allows one to construct the required Galois cover and this is not the case for the configuration mentioned in Theorem \ref{thm:xiao }. For example, considering disjoint rational curves, one has that if a K3 surface admits 16 disjoint rational curves, it admits a double cover with a surface birational to a complex torus (see \cite{NKummer}), but if a K3 surface admits 8 disjoint rational curves we don't know if it admits a double cover with a surface birational to a K3 surface. 
}
\end{remark}


By the previous Remark, the existence of a given ADE configuration $B_{G}$ of rational curves on a K3 surface $S$ is not enough to guarantee the existence of a (rational) $G$-cover of $S$ with a K3 surface. Nevertheless, Proposition \ref{prop: MG} shows that there is a lattice theoretic condition characterizing the K3 surfaces admitting a (rational) $G$-cover with another K3 surface: it is not sufficient to consider the lattice associated with the ADE configuration $B_G$, but one has to consider a specific overlattice of it.

To fix the notation needed in what follows, we recall that the second cohomology group of any K3 surface, endowed with the cup product, is isometric to the unique even unimodular lattice of signature $(3,19)$. This lattice will be denoted $\Lambda_{K3}$.

Consider an ADE lattice $\Lambda$ that is not primitively embedded in the lattice $\Lambda_{K3}$. There might be an overlattice $\Lambda'$ of $\Lambda$ of finite index which can be primitively embedded in $\Lambda_{K3}$, and in this case $\Lambda$ is embedded (but not primitively) in $\Lambda_{K3}$. Since $\Lambda'/\Lambda$ is a product of cyclic groups of finite order, see \cite{N}, $\Lambda'$ is generated by the generators of $\Lambda$ and by some linear combinations with rational coefficients of them, which yields divisible classes.

\begin{lemma}\label{lemma: Lambda1 and cyclic cover} Let $\Lambda$ be an ADE configuration of rational curves on $S$. There exists an overlattice $\Lambda'$ of $\Lambda$ of finite index $r>1$ which is embedded in $NS(S)$ if and only it there exists an effective divisible class in $NS(S)$. The latter implies the existence of a cyclic branched cover of $S$.
\end{lemma}

\proof
Assume that a given set of smooth rational curves on a K3 surface $S$ is divisible. Hence there exists a linear combination $D$ (with integer coefficient) of these curves such that $\frac{1}{n}D\in NS(S)$ for some positive integer $n$. This implies that the lattice $\Lambda$ spanned by the set of these smooth rational curves is embedded in $NS(S)$, but not primitively embedded, and that the lattice  $\Lambda':=\langle \Lambda, \frac{1}{n}D\rangle$ is an overlattice of finite index of $\Lambda$ embedded in $NS(S)$.

Conversely, if $\Lambda':=\langle \Lambda, \frac{1}{n}D\rangle$ is embedded in $NS(S)$, then $\frac{1}{n}D$ is a divisible class and we can assume it is effective (otherwise one can add classes of curves in such a way that the coefficients of the irreducible curves appearing in $\frac{1}{n}D$  are all contained in $(0,1)$, i.e. that the coefficients in $D$ are integers in $\{1,2, \ldots, n-1\}$). So, $\frac{1}{n}D$ is an effective divisibile class associated to an $n$-cyclic cover.\endproof

We now describe the shape of the divisible classes and fix some preliminary notation: 
\begin{itemize}
    \item given an integer $n\geq 1$, let $C_1,\cdots,C_n$ be the irreducible curves in an \textit{$A_n$-configuration} on a surface $S$ numbered in such a way that $C_{i}C_{i+1}=1$ for every $i=1,\ldots,n-1$ and $C_{i}C_{j}=0$ if $|i-j|\geq 2$;
    \item given integers $n_{1},\cdots,n_{m}\geq 1$, we let $C_{i}^{(j)}$ be the irreducible curves in the configuration $A_{n_{j}}$, and $C_{1}^{(j)},\cdots,C_{n_{j}}^{(j)}$ are numbered as described in the previous point.
    \item we denote $\sum_{i=1}^{n_{j}}[\alpha_{i}]_{n}C_{i}^{(j)}$ the linear combination such that  $\alpha_i\in\mathbb{Z}$ and $[\alpha_i]_n$ is the class of $\alpha_i$ in $\Z/n\Z$. 
    \item we let $$[V]_{n}^{(j)}:=\sum_{i=1}^{n_{j}}[i]_{n}C_{i}^{(j)},\,\,\,\,\,\,\,\,[kV]_{n}^{(j)}:=\sum_{i=1}^{n_{j}}[ki]_{n}C_{i}^{(j)}.$$  In particular we observe that $[V]_n^{(j)}/n$ spans the discriminant group of $A_{n_{j}}$.
\end{itemize}

\begin{prop}
\label{prop: divisible}
Let $\Lambda$ be the lattice of an ADE configuration of rational curves on a K3 surface $S$ and $\Lambda'$ be an overlattice due to the presence of effective divisible classes. Then there are only a finite number of ways to construct $\Lambda'$, each of them obtained by adding a finite number of divisible classes, and in each case the lattice $\Lambda$ has to contain a particular sublattice $\Gamma$ that depends on the divisible classes one has to add. The lattices $\Gamma$ with minimal rank which admit an overlattice with a prescribed lattice $\Gamma'$ are listed in Table \ref{table:1}: the first column gives the number of divisible classes that are added; the second column gives the divisibility of the added classes; the third column gives the lattice $\Gamma$; the fourth column gives the class that becomes divisible; the fifth column gives the Galois  group of the cover branched on all the divisible classes; the last column has a $T$ if the cover is a torus, and has K3 otherwise.
\end{prop}
\begin{center}
\begin{table}[h!]
\centering
\caption{Possibilities for overlattices of a given ADE lattice}
\label{table:1}
\begin{tabular}{|c|c|c|c|c|c|}
\hline

$m$ & $n_i$ & $\Gamma$ & & $G=\Lambda'/\Lambda$ & \\
\hline\hline

1 & 2 & $A_1^{\oplus 8}$ & $\sum_{j=1}^8[V]_2^{(j)}$ 
& $\Z/2\Z$ & K3\\
\hline

2 & $\begin{array}{c}2\\2\end{array}$ & $A_1^{\oplus 12}$ & $\begin{array}{c} \sum_{j=1}^{8}[V]_{2}^{(j)}\\ \sum_{j=1}^4[V]_2^{(j)}+\sum_{h=9}^{12}[V]_2^{(h)}\end{array}$ 
& $(\Z/2\Z)^2$ & K3\\
\hline

3 & $\begin{array}{c}2\\2\\2\end{array}$ & $A_1^{\oplus 14}$ & $\begin{array}{c}\sum_{j=1}^8[V]_2^{(j)}\\ \sum_{j=1}^4[V]_2^{(j)}+\sum_{h=9}^{12}[V]_2^{(h)}\\ \sum_{j=0}^4\left([V]_2^{(1+4j)}+[V]_{2}^{(4j+2)}\right)\end{array}$ & $(\Z/2\Z)^3$ & K3\\
\hline

4 & $\begin{array}{c}2\\2\\2\\2\end{array}$ & $A_1^{\oplus 15}$ & $\begin{array}{c}\sum_{j=1}^8[V]_2^{(j)}\\ \sum_{j=1}^4[V]_2^{(j)}+\sum_{h=9}^{12}[V]_2^{(h)}\\ \sum_{j=0}^4\left([V]_2^{(1+4j)}+[V]_2^{(4j+2)}\right)\\ \sum_{j=1}^8[V]_2^{(2j-1)}\end{array}$ & $(\Z/2\Z)^4$ & K3\\
\hline

5 & $\begin{array}{c}2\\2\\2\\2\\2\end{array}$ & $A_1^{\oplus 16}$ & $\begin{array}{c}\sum_{j=1}^8[V]_2^{(j)}\\ \sum_{j=1}^4[V]_2^{(j)}+\sum_{h=9}^{12}[V]_2^{(h)}\\ \sum_{j=0}^4\left([V]_2^{(1+4j)}+[V]_{2}^{(4j+2)}\right)\\ \sum_{j=1}^8[V]_2^{(2j-1)}\\ \sum_{j=1}^{16}[V]_2^{(i)}\end{array}$ & $(\Z/2\Z)^5$ & $T$\\
\hline

1 & 3 & $A_2^{\oplus 6}$ & $\sum_i^6[V]_3^{(i)}$ & $\Z/3\Z$ & K3\\
\hline

2 & $\begin{array}{c}3\\3\end{array}$ & $A_2^{\oplus 8}$ & $\begin{array}{c}\sum_{j=1}^6[V]_3^{(j)}\\ \sum_{j=1}^3[2V]_3^{(j)}+\sum_{h=7}^{9}[V]_3^{(h)}\end{array}$ & $(\Z/3\Z)^2$ & K3\\
\hline

3 & $\begin{array}{c}3\\3\\3\end{array}$ & $A_2^{\oplus 9}$ & $\begin{array}{c}\sum_{j=1}^6[V]_3^{(j)}\\ \sum_{j=1}^3[2V]_3^{(j)}+\sum_{h=7}^{9}[V]_3^{(h)}\\ \sum_{j=1}^4\left([V]_3^{(2j-1)}+[2V]_3^{(2j)}\right)+[V]_3^{(9)}\end{array}$ & $(\Z/3\Z)^3$ & $T$\\
\hline

1 & 5 & $A_4^{\oplus 4}$ & $\sum_{i=1}^2[V]_5^{(2i-1)}+[2V]_5^{(2i)}$ & $\Z/5\Z$ & K3\\
\hline

1 & 7 & $A_6^{\oplus 3}$ & $[V]_7^{(1)}+[2V]_7^{(2)}+[3V]_7^{(3)}$ & $\Z/7\Z$ & K3\\
\hline

1 & 4 & $A_3^{\oplus 4}\oplus A_1^{\oplus 2}$ & $\sum_{i=1}^4[V]_4^{(i)}+\sum_j=5^6 [V]_2^{(j)}$ & $\Z/4\Z$ & K3\\
\hline

2 & $\begin{array}{c}4\\2\end{array}$ & $A_3^{\oplus 4}\oplus A_1^{\oplus 4}$ & $\begin{array}{c}\sum_{i=1}^4[V]_4^{(i)}+\sum_{j=5}^6 [V]_2^{(j)}\\ \sum_{i=1}^2[2V]_4^{(i)}+\sum_{j=5}^8 [V]_2^{(j)}\\ \end{array}$ & $\Z/4\Z\times \Z/2\Z$ & K3\\
\hline

2 & $\begin{array}{c}4\\4\end{array}$ & $A_3^{\oplus 6}$ & $\begin{array}{cc} \sum_{i=1}^4[V]_4^{(i)}+ [2V]_4^{(5)}\\ \ [2V]_4^{(2)}+[3V]_4^{(2)}+ \sum_{i=3}^6[V]_4^{(i)}\end{array}$ & $(\Z/4\Z)^2$ & K3\\
\hline

3 & $\begin{array}{c}4\\2\\2\end{array}$ & $A_3^{\oplus 4}\oplus A_1^{\oplus 6}$ & $\begin{array}{cc} \sum_{i=1}^4[V]_4^{(i)}+ \sum_{j=5}^6[V]_2^{(j)}\\ \sum_{i=1}^2[2V]_4^{(i)}+\sum_{j=5}^8[V]_2^{(j)}\\ \sum_{i=2}^3[2V]_4^{(i)}+\sum_{j=5}^6[V]_2^{(j)}+\sum_{h=9}^{10}[V]_2^{(h)}\end{array}$ & $(\Z/4\Z)\times(\Z/2\Z)^2$ & $T$\\
\hline

1 & 8 & $A_7^{\oplus 2}\oplus A_3\oplus A_1$ &  $[V]_8^{(1)}+[3V]_8^{(2)}+[V]_4^{(3)}+[V]_2^{(4)}$ & $\Z/8\Z$ & K3\\
\hline
\end{tabular}
\end{table}
\end{center}
\proof As an example we discuss the cases related with $2-$divisible classes, which are essentially contained in \cite{NKummer} and \cite{Nsympl}. 

In \cite{NKummer} it is proved that the $2-$divisible classes of disjoint rational curves consists either of 8 or of 16 curves. In the first case the double cover surface is birational to a K3 surface, in the latter to a $2-$dimensional complex torus. 
	
Moreover, it is also proved that on a K3 surface there are at most 16 disjoint smooth rational curves and if there are 16 disjoint smooth rational curves, then they form a divisible set. In this case there are also other divisible classes, which are the one described in the case 5) in Table \ref{table:1}.
	
Assume now that there are less than 16 disjoint smooth rational curves. Each divisible classed is of the sum of 8 curves, and two different divisible classes share some curves. Moreover, If $v_1$ and $v_2$ are two $2-$divisible classes sharing $k$ curves, then $v_1+v_2$ gives a divisible class which is the sum of $16-2k$ curves, so that $k=4$. 	
It follows that all the $2-$divisible classes are given by the sum of 8 disjoint rational curves, and that two classes have exactly 4 curves in common. This forces the choice of the classes $v_1$, $v_2$, $v_3$ and $v_4$ to be as described. 
	
All the other cases in which the cover is a K3 surface are discussed in \cite{Nsympl}, the ones in which the cover is a torus are discussed in \cite{Bertin}.\endproof

The existence of a divisible class of order $pq$ with $(p,q)=1$ is equivalent to the existence of both a $p$-divisible class and of a $q$-divisible class. Since by \cite{Nsympl} the order of the divisible class is at most 8, combining the previous results one obtains the cases which there are covers with Galois groups $\Z/6\Z$ (composition of cyclic covers by $\Z/3\Z$ and $\Z/2\Z$) and $\Z/6\Z\times \Z/2\Z$ (composition of the Galois cover with Galois groups $(\Z/2)^2$ and $\Z/3\Z$). 

\subsubsection{Galois covers of K3 surfaces with surfaces birational to a K3 surface}\label{subsub Xiao}
Let $B_G$ be a configuration of curves which arises by the desingularization of the quotient of a K3 surface by the action of a group of symplectic automorphisms (e.g. 8 disjoint rational curves). We saw in Example \ref{esem: ADE no imples Galois} that the presence of $B_G$ on a K3 surface does not imply that the surface admits a Galois cover (e.g. the 8 disjoint rational curves may or may not be a divisible set). In particular, the existence of the cover is not only related with the presence of the lattice $\Lambda\simeq B_G$  inside the N\'eron--Severi group of the K3 surface, but also with the existence of a certain overlattice $\Lambda'$ of $\Lambda\simeq B_G$ embedded in the N\'eron--Severi group. In the following proposition we discuss the properties of such overlattices $\Lambda'$, which will be denoted by $M_G$ for historical reasons.

\begin{prop}
\label{prop: MG} 
Let $G$ be a group acting symplectically on a K3 surface and $B_G$ (resp. $B^{i}_{G}$ for $i=1,2$ for $G=Q_{8},T_{24}$) the corresponding ADE configuration in  \cite[Table 2]{Xi}.
\begin{enumerate}
    \item There is a finite index overlattice $M_G$ (resp. $M_G^i$) of $B_G$ (resp. $B_G^i$) characterizing the K3 surfaces $S$ admitting a (rational) $G$-cover with another K3 surface, i.e., if there is a K3 surface $Y$ such that $G\subset\Aut(Y)$ and $S$ is birational to $Y/G$ then $M_G$ (resp. $M^{i}_{G}$) is primitively embedded in $NS(S)$.
    \item If $M_G$ (resp. $M^{i}_{G}$) is embedded in $NS(S)$ for a K3 surface $S$, then the embedding is primitive and $S$ admits a $G$-Galois  
    cover with a K3 surface.
    \item The quotient $M_G/B_G$ (resp. $M^{i}_{G}/B^{i}_{G}$) is not trivial if and only if $G\neq \mathfrak{A}_5, \mathfrak{A}_6,\ L_2(7) (=\textrm{PSL}(2,7)), M_{20}$. 
\end{enumerate} 
\end{prop}

\proof The complete list of the groups acting symplectically on a K3 surface is known after \cite{M2}, \cite{K}, \cite{Xi}. In particular, in \cite{Xi} the set of points with non trivial stabilizer for the action of $G$ on a K3 surface are given: this provides the list of the lattices $B_{G}$ contained in \cite[Table 2]{Xi} of the ADE configurations resolving the singular points of the quotient. 
The lattice $M_G$ is by definition the minimal primitive sublattice, of the N\'eron--Severi group of a K3 surface which is the minimal resolution of the quotient of a K3 by a symplectic automorphism, containing the curves arising by the desingularization of the quotient, which are the curves in $B_G$. Hence it is primitively embedded in the N\'eron--Severi of the surface.

One can directly check on \cite[Table 2]{Xi} that there are no different groups attached to the same ADE configuration.

Moreover, in the seventh column of \cite[Table 2]{Xi}, there is a description of the divisible sets of curves appearing in the lattice $M_G$. With the exceptions mentioned in the statement, for all the groups $G$ appearing, $M_G/B_G$ is not trivial, i.e. there are divisible effective classes. Since every divisible class is associated with a cyclic branched cover, we obtain that if $G$ appears in the list \cite[Table 2]{Xi} and $G\neq \mathfrak{A}_5, \mathfrak{A}_6,\ L_2(7), M_{20}$, then $S$ admits a cyclic (and in particular Galois) cover.  
The uniqueness of $M_G$ is proved in \cite[Lemma 6]{Xi} and the fact that the embedding of $M_G$ in the N\'eron--Severi group of a K3 surface implies the presence of a branched $G$-Galois cover is proved in \cite{Nsympl} for the cyclic group $G$ and in \cite{Wi} for all the other groups $G\not\in\{\mathfrak{A}_5, \mathfrak{A}_6,\ L_2(7), M_{20}\}$. 

The cases of the ADE configuration $B_G$ where $G\in\{\mathfrak{A}_5, \mathfrak{A}_6,\ L_2(7), M_{20}\}$ are discussed in \cite{Wi}. More precisely, if $G\in\{\mathfrak{A}_5, \mathfrak{A}_6,\ L_2(7), M_{20}\}$, in \cite[Question 3.1]{Wi} and in the discussion below, the author proves that if $M_G=B_G$ admits a unique (up to isometry) embedding in $\Lambda_{K3}$, then all the K3 surfaces with the ADE configurations of curves given by $B_G$ are $G$-covered by a K3 surface (essentially because the fundamental group of the surface obtained by removing the ADE configuration $B_G$ is necessarily $G$). 

If $G=\mathfrak{A}_5, M_{20}$, then the embedding of $B_G=M_G$ in $\Lambda_{K3}$ is unique, while if $G=\mathfrak{A}_6, L_{2}(7)$ there are two different embeddings of $B_G=M_G$ in $\Lambda_{K3}$  (see \cite{Wi}). Nevertheless, if $G=\mathfrak{A}_{6}, L_{2}(7)$ there are also two different symplectic actions of $G$ on a K3 surface. Indeed, a symplectic action of $G$ on a K3 surface determines an action of $G$ on $\Lambda_{K3}$ and hence two lattices $\Gamma_G:=\Lambda_{K3}^{G}$ and $\Omega_G:=\left(\Lambda_{K3}^{G}\right)^{\perp}$. In \cite{Ha} the lattices $\Gamma_G$ are determined for each group $G$ acting symplectically on a K3 surface. In particular, if $G\in\{\mathfrak{A}_6, L_{2}(7)\}$ there are two different lattices $\Gamma_G^1$ and $\Gamma_{G}^{2}$, meaning that $G$ has two different actions on $\Lambda_{K3}$, both related with a symplectic action on a K3 surface. 

Let now $Y_i$, $i=1,2$ be a K3 surface whose transcendental lattice $T_{Y_i}$ is isometric to $\Gamma_G^i$. Then $Y_i$ admits a symplectic action of $G$, and $Y_1$ and $Y_{2}$ are general members of two different families of K3 surfaces with a symplectic action of $G$. Hence $Y_i/G$ is a singular surfaces whose resolution $S_i$ is a K3 surface containing an ADE configuration $B_G$, and $S_1$ and $S_2$ are members of different families of K3 surfaces. Each of them corresponds to one of the two different embeddings of $B_G=M_G$ in $\Lambda_{K3}$. We conclude that each embedding of $B_G=M_G$ with $G\in\{\mathfrak{A}_6, L_{2}(7)\}$ corresponds to a K3 surface admitting a $G$-Galois cover.\endproof

\begin{remark}
{\rm By a direct inspection, there are no lattices in the list of Theorem \ref{thm:fujiki} which are contained in the list in \cite[Table 2]{Xi}: this means that the lattices associated with the existence of a cover of a K3 surface $S$ with another K3 surface do not contain any lattice associated with the existence of a cover of $S$ with a torus.}
\end{remark}

\begin{remark}
\label{rem: smooth quasi etale cover}
{\rm Let $Y$ be a K3 surface admitting a symplectic action of the group $G$. Let $S$ be the resolution of $Y/G$ and $B_G$ the lattices generated by the curves resolving the singularities. Let $X_{B_G}$ the contraction of $B_G$. Then $X_{B_G}=Y/G$ and the quasi-\'etale cover $Y\ra Y/G=X_{B_{G}}$ is smooth. Vice versa, if a Galois quasi-\'etale cover of $X_B$ (a contraction of a K3 surface $S$) is smooth, then it is induced by a quotient of a group $G$ acting symplectically on the cover surface; in particular it is branched on all the singular points of $X_B$ and $B$ is either a $B_G$ in \cite[Table 2]{Xi} or in the list of Theorem \ref{thm:fujiki}.}
\end{remark}

\subsection{Simple symplectic surfaces}\label{subsec: simple symplectic surfaces}

We are now in the position to prove a characterization of simple symplectic surfaces in terms of contractions of curves on a K3 surface, which is as follows:

\begin{thm}
\label{thm:mainsimple}
Let $S$ be a K3 surface, $B$ an ADE configuration of smooth rational curves on $S$, $f_{B}:S\longrightarrow X_{B}$ the contraction morphism, and suppose that $X_{B}$ is an irreducible symplectic surface. The following are equivalent:
\begin{enumerate}
    \item the surface $X_{B}$ is simple symplectic;
    \item for every subconfiguration $B'$ of $B$ appearing in \cite[Table 2]{Xi}, there is no Galois covering of $S$ branched along $B'$;
    \item either there are no subconfigurations of $B$ which are of the form $B_G$ for a $G$ as in \cite[Table 2]{Xi} or for every subconfiguration of $B$ of the form $B_G$ there is no an embedding of the corresponding $M_G$ in $NS(S)$.
\end{enumerate}
\end{thm}

\proof Suppose first that $X_{B}$ is simple, and let $B'$ be a subconfiguration of $B$ that appears in \cite[Table 2]{Xi}. Suppose that there is a Galois covering $\pi:W\longrightarrow S$ branched along $B'$. Let $g:W\longrightarrow W'$ be the contraction of the curves in $\pi^{-1}(B')$: hence there is a finite quasi-\'etale covering $\pi':W'\longrightarrow X_{B}$ whose degree is the one of $\pi$. This implies the existence of a non-trivial element of the fundamental group of the smooth locus of $X_{B}$, contradicting the fact that $X_{B}$ is simple. 

Conversely, suppose that for every subconfiguration $B'$ of $B$ appearing in \cite[Table 2]{Xi}, there is no Galois covering of $S$ branched along $B'$. If $X_{B}$ was not simple, there would be a non-trivial finite quasi-\'etale covering $\pi:Y\longrightarrow X_{B}$. Since $X_{B}$ is an irreducible symplectic surface, the surface $Z:=Y\times_{X_{B}}S$ has a K3 surface as a minimal model, and it is a finite covering of $S$. As in the proof of Theorem \ref{thm:mainfujiki}, this implies that $S$ has a Galois covering by a surface that is birational to a K3 surface, hence by Theorem \ref{thm:xiao } the branch locus $B'$ is in \cite[Table 2]{Xi}. Since $B'$ is contained in the exceptional locus of $f_{B}$, which is $B$, we get a contradiction. 

By Proposition \ref{prop: MG}, the existence of a Galois cover branched on $B'$ implies that $B'$ is one of the ADE configurations $B_G$ appearing in \cite[Table 2]{Xi} and that the overlattice $M_G$ of $B_G$ is embedded in the N\'eron--Severi group of $S$.
\endproof

As a particular case we get the following:

\begin{cor}
Let $S$ be a K3 surface and $B$ an ADE configuration of $(-2)-$curves on $S$ that does not contain any configuration of \cite[Table 2]{Xi}. Then the contraction $X_{B}$ of the curves in $B$ is a simple symplectic surface.
\end{cor}

\proof All the configurations listed in Theorem \ref{thm:fujiki} contain at least one configuration in \cite[Table 2]{Xi} by Theorem \ref{thm:mainfujiki} it follows that $X_{B}$ is an irreducible symplectic surface. Now apply Theorem \ref{thm:mainsimple}, since $B$ has no subconfigurations appearing in \cite[Table 2]{Xi}.\endproof

\begin{esem}
{\rm We may now consider again the contraction of the ADE configuration $B_{\Z/2\Z}\simeq A_{1}^{\oplus 8}$ for the K3 surfaces of the families $F_{1}$ and $F_{2}$ that were described in Example \ref{esem: ADE no imples Galois}.}

{\rm By Theorem \ref{thm:mainsimple} and Example \ref{esem: ADE no imples Galois}, and using the fact that there is no subconfiguration of $B_{\Z/2\Z}$ appearing in \cite[Table 2]{Xi} beside $B_{\Z/2\Z}$ itself, we then get that the contraction of the curves of $B_{\Z/2\Z}$ on a K3 surface in the family $F_{1}$ is a simple symplectic surface, while the contraction of the curves of $B_{\Z/2\Z}$ on a K3 surface in the family $F_{2}$ is a non-simple irreducible symplectic surface.}
\end{esem}

To get a complete characterization of simple symplectic surfaces by Theorem \ref{thm:mainsimple}, it remains to consider the lattices mentioned in point (3) of Proposition \ref{prop: MG}. 
\begin{remark}
\label{rem: subset S}
{\rm By point (3) of Proposition \ref{prop: MG} if $G\in\{\mathfrak{A}_5, \mathfrak{A}_6,\ L_2(7), M_{20}\}$, then $B_G=M_G$. By \cite[Table 2]{Xi} we have the following:
\begin{itemize}
    \item if $G=\mathfrak{A}_5$, then $B_G=A_4^{\oplus 2}\oplus A_2^{\oplus 3} \oplus A_1^{\oplus 4}$;
    \item if $G=\mathfrak{A}_6$, then $B_G=A_4^{\oplus 2}\oplus A_3^{\oplus 2}\oplus A_2^{\oplus 2}\oplus A_1$;
    \item if $G=L_2(7)$, then $B_G=A_6\oplus A_3^{\oplus 2}\oplus A_2^{\oplus 3}\oplus A_1$;
    \item if $G=M_{20}$, then $B_G=D_4\oplus A_4^{\oplus 2}\oplus A_2^{\oplus 3}\oplus A_1$. 
\end{itemize} 
If we let $\mathcal{S}:=\{B_{\mathfrak{A}_{5}},B_{\mathfrak{A}_{6}},B_{L_{2}(7)},B_{M_{20}}\}$, 
we have that if $M_G\simeq B_G\in \mathcal{S}$ the existence of an embedding of $B_G\in\mathcal{S}$ in the N\'eron--Severi group of a K3 surface $S$ guarantees that $S$ admits a $G$-Galois cover even if there are no divisible classes in $B_G$, as explained in the proof of Proposition \ref{prop: MG}. }
\end{remark}

Hence we have the following, which  is point (3) of Theorem \ref{thm:mainsurfaces}:

\begin{cor}
\label{cor: no lattices in S}
Let $S$ be a K3 surface, $B$ an ADE configuration of $(-2)-$curves on $S$, and suppose that the contraction $X_{B}$ of $B$ is an irreducible symplectic surface. Then $X_B$ is a simple symplectic surface if and only if the embedding of $B$ in $NS(S)$ is primitive and $B$ does not contain any configuration in the set $\mathcal{S}$ defined in Remark \ref{rem: subset S}.
\end{cor}

\proof Suppose first that $B$ is primitively embedded in $NS(S)$ and $B$ does not contain any configuration in the set $\mathcal{S}$. If $X_B$ admits a finite quasi-\'etale cover, then $S$ admits a $G$-Galois cover branched on a subconfiguration of $B$. This implies that $B$ contains a configuration $B_G$ of curves and that $NS(S)$ contains primitively $M_G$. Since $B$ does not contain any of the configuration of the set $\mathcal{S}$, we see that $G\neq \mathfrak{A}_5,\mathfrak{A}_6, L_2(7), M_{20}$ (see Remark \ref{rem: subset S}), so by Proposition \ref{prop: MG} the lattice $M_G$ is a proper finite index overlattice of $B_G$, so $B_G$ is not primitively embedded in $NS(S)$, getting a contradiction. It follows that $X_{B}$ is simple. 

Conversely, suppose that $X_{B}$ is simple and that $B$ is embedded in $NS(S)$, but either the embedding is not primitive or $B$ contains a configuration in $\mathcal{S}$. This implies that either there is a divisible class in $NS(S)$, 
or $S$ admits a $G$-cover with $G\in\{\mathfrak{A}_5,\mathfrak{A}_6, L_2(7), M_{20}\}$. In the first case the divisible class can be chosen to be effective, since the $(-2)$-classes in $B_G$ can be chosen, up to a sign, to be the classes of smooth rational curves, and the divisible class can be always chosen with positive coefficients (as in proof of Lemma \ref{lemma: Lambda1 and cyclic cover} the coefficients of the divisible class are contained in $(0,1)$). So, if there is a divisible class, there is an effective divisible class and so a cyclic cover branched on these curves. 

In both the cases, $S$ admits a Galois cover branched on some curves in $B$: this implies that $X_B$ admits a quasi-\'etale Galois cover, contradicting the fact that $X_{B}$ is simple. \endproof

As a conclusion, we summarize all the previous results in the following:

\begin{theorem}\label{theorem: recap section 4}
Let $S$ be a K3 surface, $B$ an ADE configuration of $(-2)-$curves on $S$ and $X_{B}$ the contraction of the curves in $B$.
\begin{itemize}
    \item If $B$ does not contain any configuration in the list \ref{thm:mainfujiki}, then: 
    \begin{enumerate}
        \item the surface $X_{B}$ is irreducible symplectic;
        \item the surface $X_{B}$ is simple symplectic if and only if $B$ is primitively embedded in $NS(S)$ and does not contain any configuration in $\mathcal{S}$.
    \end{enumerate}
    \item If $B$ contains a configuration in the list in \ref{thm:mainfujiki}, then $X_B$ is primitive symplectic but not irreducible symplectic.
\end{itemize}
\end{theorem}

\begin{remark}
{\rm In the proof of Theorem \ref{thm: simple} we will consider the ADE configurations which can appear on a K3 surface and which contain an ADE configuration in $\mathcal{S}$. In particular we obtain that there are exactly nine ADE configurations which are excluded by the hypothesis ``$B$ does not contain any configuration in $\mathcal{S}$" of Theorem \ref{theorem: recap section 4}.}
\end{remark}

\section{The ADE configurations on K3 surfaces}\label{sec_ADEconfig}

As we saw in the previous section, in order to classify all primitive symplectic surfaces, we need to provide a complete classification of all possible ADE configurations of $(-2)-$curves on K3 surfaces: the primitive symplectic surfaces are the contractions of these ADE configurations; the irreducible symplectic surfaces are those obtained by contracting the configurations which do not contain any ADE configuration in the list of Theorem \ref{thm:mainfujiki}, and the simple ones are obtained by contracting a configuration which is primitively embedded in the N\'eron--Severi group and does not contain moreover the configurations in the set $\mathcal{S}$, see Theorem \ref{thm: simple}.

In this section we explain the procedure which allows one to find all the ADE configurations of rational curves on a K3 surface and how to find the ones which correspond to irreducible symplectic surfaces and to simple symplectic surfaces: in Section \ref{subsec: finidinf lists} we present the procedure, which can be implemented in a computer algorithm described in Appendix \ref{sec: the algortihm GAP}; in Section \ref{subsec: proof} we recall all the results which are previously shown in order to give the explicit proof of Theorem \ref{thm:classification}.

A synthetic description of the possible ADE configurations of $(-2)-$curves on K3 surfaces is provided in \cite{Sh}, where some conditions to find a complete list of these configurations are given, even if there is neither an explicit algorithm nor a complete list. The purpose of this section is to describe an explicit and different algorithm which allows to produce a complete list and to distinguish among the ADE configurations associated to specific type of primitive symplectic varieties. Notice that also the papers \cite{cat1}, \cite{cat2} contain an algorithmic approach to the problem of classifying singularities on K3 surfaces.\\

Given a lattice $L$, the length $\ell(L)$ of $L$ is the minimal number of generators of the discriminant group $A_L:=L^{\vee}/L$. Recall moreover that the discriminant group $A_{L}$ of a lattice $L$ is naturally endowed with a quadratic form $q_L$ induced by the bilinear form on $L$ and called discriminant form. 

The discriminant group and form of the ADE lattices are as follows: $$\left(A_{L}, q_L\right)=\left\{\begin{array}{ll}
	\left(\Z/(i+1)\Z,\left[-\frac{i}{i+1}\right]\right) & L=A_{i}\\ \left(\Z/4\Z,\left[\frac{-(2n+1)}{4}\right]\right) & L=D_{2n+1}\\ \left((\Z/2\Z)^2,\left[\begin{array}{cc}-\frac{n}{2}&-\frac{1}{2}\\-\frac{1}{2}&-1\end{array}\right]\right) & L=D_{2n}\\ \left(\{1\},-\right) & L=E_{8}\\ (\Z/2\Z,\left[-\frac{3}{2}\right]) & L=E_{7}\\ \left(\Z/3\Z,\left[-\frac{4}{3}\right]\right) & L=E_{6}.\end{array}\right.$$

Finally, we recall that the discriminant group of a direct sum of lattices is the product of the discriminant groups of each summands and the discriminant form is the one induced, i.e., two different summands are orthogonal with respect to the discriminant form and on each summand the discriminant form is the one of the associated lattice. 

This allows to compute the discriminant group, the discriminant form and the length of each possible ADE configuration.

\begin{lemma}\label{lem:twocond} If a negative definite lattice $L$ is primitively embedded in $\Lambda_{K3}$, then $\rk(L)\leq 19$ and $\ell(L)\leq min(\rk(L), 22-\rk(L))$.\end{lemma}

\proof The lattice $\Lambda_{K3}$ is the unique unimodular even lattice of signature $(3,19)$ and the lattice $\Lambda$ is an even negative definite lattice, so the result follows by \cite[Theorem 1.12.2]{N}. For the reader convenience we write an explicit proof.

Since $L$ is negative definite, the rank of $L$ has to be less than or equal to the one of the negative part of $\Lambda_{K3}$. Moreover, for every lattice $M$, $\ell(M)\leq \rk(M)$, which clearly implies that $\ell(L)\leq \rk(L)$. Since $\Lambda_{K3}$ is unimodular and $L$ is primitively embedded in $\Lambda_{K3}$, its orthogonal complement $L^{\perp_{\Lambda_{K3}}}$ has the same discriminant group of $L$ (see \cite[Proposition 1.6.1]{N}) and in particular the same length. So $\ell(L)=\ell(L^{\perp_{\Lambda_{K3}}})\leq \rk(L^{\perp_{\Lambda_{K3}}})=22-\rk(L)$, which concludes the proof. \endproof

\begin{example}
\label{example: 8A2}
{\rm Let us consider the ADE configuration $\Lambda=A_2^{\oplus 8}$. Then $\rk(\Lambda)=16$ and $\ell(\Lambda)=8$, hence by Lemma \ref{lem:twocond}, there is no primitive embedding of $\Lambda$ in $\Lambda_{K3}$. Nevertheless, by Proposition \ref{prop: divisible} there are two finite index overlattices of $\Lambda$, one of index 3 and one of index $3^2$ (see Table \ref{table:1}, $m=1$, $n_i=(3)$ and $m=2$, $n_i=(3,3)$). Denote $\Lambda'_1$ and $\Lambda'_2$ these overlattices. We observe that the rank does not change considering finite index overlattices, i.e. $\rk(\Lambda)=\rk(\Lambda_1')=\rk(\Lambda_2')=16$, but the length drops, in particular $\ell(\Lambda_1')=6$ and $\ell(\Lambda_2')=4$. So the lattices $\Lambda_i'$ could be primitively embedded in $\Lambda_{K3}$, and if so we have a chain of embeddings (not all primitive) $\Lambda\subset\Lambda_i'\subset\Lambda_{K3}$ which suffices to guarantee that there is a K3 surface with the configuration $\Lambda=A_2^{\oplus 8}$ of curves.}
\end{example}

Example \ref{example: 8A2} shows that in order to understand if a given ADE configuration $\Lambda$ appears on a K3 surface, it is necessary to consider also its overlattices $\Lambda'$, which can be primitively embedded in $\Lambda_{K3}$. Our aim is to produce such an overlattice $\Lambda'$ for a given lattice $\Lambda$. To do so, we will add to $\Lambda$ some divisible classes, which are the generators of $\Lambda'/\Lambda$ and which have to satisfy the hypotheses of Proposition \ref{prop: divisible}, i.e. they have to be contained in Table \ref{table:1}.

From a lattice theoretic point of view the situation is the following. Let $L$ be a lattice and $n\in\mathbb{N}$.  Let $v$ be a linear combination of the elements of a basis of $L$ with coefficients in $\frac{1}{n}\mathbb{Z}$. If $\langle L,v\rangle$ is an even lattice, the linear combination $v$ has to satisfy some conditions: first of all, in order to get that $\langle L,v\rangle$ is a lattice (with the natural pairing induced by $\mathbb{Q}-$linear extension by the one of $L$), we need to have $\lambda\cdot v\in\mathbb{Z}$ for every $\lambda\in\Lambda$, which implies $v\in L^{\vee}$. Moreover, in order to have that $\langle L,v\rangle$  is even, we need $v\cdot v\in 2\Z$.  The two conditions we mentioned are well known and are equivalent to require that $v\in A_{L}$ spans an isotropic subspace of $A_{L}$ with respect to the discriminant form $q_{L}$, see \cite[Section 4$^o$]{N}.

In our context, $L$ is an ADE lattice $\Lambda$ which is embedded in the N\'eron--Severi group of a K3 surface $S$, so we need moreover that the $(-2)$-classes generating $\Lambda$ represent (up to a sign) smooth irreducible curves on $S$. This imposes several further restrictions on $v$, which in particular can not represent a $-2$ class, see e.g. \cite{G} and \cite{Sc}. 

So, to construct an overlattice $\Lambda'$ of finite index of $\Lambda$, we have to add a class $v$ such that $v\in L^{\vee}$, $v\cdot v\in 2\Z$ and $v\cdot v< -2$. By construction, $nv\in\Lambda$ and it is divisible since $v\in NS(S)$. Moreover, it can be chosen to be effective, as in proof of Lemma \ref{lemma: Lambda1 and cyclic cover}.

This lattice theoretic setting implies Proposition \ref{prop: divisible} for the lattices with minimal rank among those admitting a given set of divisible classes. Of course, the configurations listed in Proposition \ref{prop: divisible} may appear as sublattices of other ADE configurations of higher rank. For example, the lattice $A_7$ contains the lattice $A_3^{\oplus 2}$, and so if one is looking for e.g. the configuration $A_3^{\oplus 4}\oplus A_1^{\oplus 2}$, one can find it embedded in $A_7^{\oplus 2}\oplus A_1^{\oplus 2}$. Similar phenomena appear considering the lattices $D_n$, $E_6$ and $E_7$. For example, the lattice $D_5$ contains $A_3\oplus A_1$ and the generator of $A_{D_5}$ contains classes both in the copy of $A_3$ and in the copy of $A_1$ embedded in $D_5$.  If one is looking for e.g. the configuration $A_3^{\oplus 4}\oplus A_1^{\oplus 2}$, one can find it embedded in $D_5^{\oplus 2}\oplus A_3^{\oplus 2}$ (with each of the two copies of $A_3\oplus A_1$ embedded in on $D_5$). We don't describe all these kinds of embeddings, since this would produce a long and tedious list, but the program described in Section \ref{sec: the algortihm GAP} recognizes and uses all of them. 

To give just an idea of the possible embeddings of ADE lattices in ADE lattices, we observe that $A_m$ can be embedded in $A_{m'}$ if $m\leq m'$, in $D_n$ if $m\leq n-1$ and in $E_h$ if $m\leq h$; the lattices $D_n$ cannot be embedded in $A_m$, but can be embedded in $D_{n'}$ if $n\leq n'$ and in $E_h$ if $n\leq h-2$; the lattice $E_h$ can be embedded only in the lattices $E_{h'}$ with $h\leq h'$. The interested reader may refer to \cite{Nish} for the proofs of these and similar results. The situation becomes more complicated when one considers direct sums of ADE lattices, for which a computer aided calculation is needed.

All the previous results lead to the main one of this section:

\begin{theorem}
\label{theorem: existence of ADE configuration}
A K3 surface $S$ contains an ADE configuration of smooth rational curves associated to a lattice $\Lambda$ if and only if either $\Lambda$ or an overlattice of finite index $\Lambda'$ of $\Lambda$ obtained as in Proposition \ref{prop: divisible} is primitively embedded in $\Lambda_{K3}$.
\end{theorem}

\proof The result is a consequence of Lemma \ref{lemma: primitive embdding of Lambda and existence K3} and of the fact that when one adds the divisible classes described in Proposition \ref{prop: divisible}, one is not adding $(-2)-$classes. This guarantees that the $(-2)$-classes generating $\Lambda$ are (up to a sign) classes of smooth irreducible rational curves, see e.g. \cite[Proposition 3.2]{G} and \cite[Corollary 1.2]{Sc}.
\endproof

Hence, in order to establish if a K3 surface with a given ADE configuration $B$ of rational curves exists, one just needs to understand if the lattice $\Lambda$ associated to $B$, or an overlattice $\Lambda'$ obtained from $\Lambda$ by adding divisible classes as in Proposition \ref{prop: divisible} is primitively embedded in $\Lambda_{K3}$. This is a completely lattice theoretical problem, and will be treated in the next section.

\subsection{Finding the list of all the ADE configurations of rational curves}\label{subsec: finidinf lists}

We start by recalling some of the results of Nikulin about the existence of an embedding of a lattice into a unimodular one and  we observe that $\Lambda$ (and hence its overlattice $\Lambda'$) is, in our context, a negative definite lattice, since the ADE-lattices are negative definite.
The main result is \cite[Theorem 1.12.2]{N}, that we state here under the assumption that the unimodular lattice is $\Lambda_{K3}$ and the lattice $\Lambda$ is negative definite, which is the case of our interest. In the following, if $\Lambda$ is a lattice, we will let $\Lambda_p$ be the $p$-adic lattice, see \cite[Section 7$^{o}$]{N}.

\begin{theorem}
\label{thm: criterion of existence embedding Nikulin}
Let $M$ be a negative definite even lattice with discriminant form $q_M$. The following properties are equivalent:
\begin{itemize}
    \item There is a primitive embedding of $M$ into $\Lambda_{K3}$ with rank $\rk(M)$, signature $(0,\rk(M))$ and discriminant form $q_M$.
    \item  The following conditions are simultaneously satisfied:
    \begin{enumerate} 
	\item $\rk(M)\leq 19$, $22-\rk(M)\geq \ell(M)$;
	\item if $p$ is an odd prime such that $22-\rk(M)=\ell(M_p)$, then we have that $-|A_{M}|=d(q_{A_{M_p}})\mod (\Z_p^*)^2$;
	\item if $22-\rk(M)=\ell(M_2)$ then either $q_M=q_{\langle 2\rangle}\oplus q_{M'}$ for a certain lattice $M'$ or $|A_{M}|=\pm d(q_{A_{M_p}})\mod (\Z_2^*)^2$.
    \end{enumerate}
\end{itemize}
\end{theorem}

The previous Theorem has an important consequence, which is \cite[Corollary 1.12.3]{N}, that again we state here under the assumption that the unimodular lattice is $\Lambda_{K3}$ and $\Lambda$ is negative definite.

\begin{corollary}
	\label{cor: criterion of existence embedding Nikulin}
	If $M$ is a lattice such that $\rk(M)\leq 19$ and $\ell(M)\leq 21-\rk(M)$, then there exists a primitive embedding of $M$ in $\Lambda_{K3}$.
\end{corollary}

\begin{example}
{\rm In Example \ref{example: 8A2} we constructed two overlattices of $\Lambda\simeq A_2^{\oplus 8}$. The lattice $\Lambda'_2$ satisfy the condition of Corollary \ref{cor: criterion of existence embedding Nikulin}. Therefore, $\Lambda_2'$ admits a primitive embedding in $\Lambda_{K3}$ and by Theorem \ref{theorem: existence of ADE configuration} $A_2^{\oplus 8}$ appears as configuration of rational curves on a K3 surface.

The overlattice $\Lambda_1'$ does not satisfy the condition of Corollary \ref{cor: criterion of existence embedding Nikulin}; nevertheless, it satisfies Theorem \ref{thm: criterion of existence embedding Nikulin} (in particular for $p=3$ one has to check the condition (2) in Theorem \ref{thm: criterion of existence embedding Nikulin}). So, there exists also a primitive embedding of $\Lambda_1'$ in $\Lambda_{K3}$, which provides a different proof of the fact that $A_2^{\oplus 8}$ is a configuration which appears on K3 surfaces.}
\end{example}

We first construct all the possible ADE configurations whose corresponding lattice $\Lambda$ is such that $\rk(\Lambda)\leq 19$, and we let $\mathcal{L}_{tot}$ be the list of all the corresponding lattices. 

Starting from $\mathcal{L}_{tot}$ we extract the sublist $\mathcal{L}_{real}$ of the ADE configurations that are realized as ADE configurations of $(-2)-$curves on some K3 surfaces. To do so we consider the following procedure which has to be applied to each $\Lambda\in \mathcal{L}_{tot}$.  The algorithm is explained in details in the Appendix \ref{sec: the algortihm GAP}, here we explain the mathematical issues on which it is based. Our algorithm produces, in a more explicit way and with more information, the same list which can be calculated by Shimada's algorithm.

First we observe that a finite index overlattice $\Lambda'$ of $\Lambda$ has the same rank of $\Lambda$, but the length of $\Lambda'$ is strictly less than the length of $\Lambda$ if the index of the inclusion is $r>1$. Moreover, we recall that if a lattice $\Lambda$ is embedded in $\Lambda_{K3}$, then its saturation in $\Lambda_{K3}$ is by definition primitively embedded in it. So, if $\Lambda$ admits an embedding in $\Lambda_{K3}$, then there exists an overlattice of $\Lambda$ of finite index $r$ (possibly $r=1$) of $\Lambda$ (i.e. its saturation) which is primitively embedded in $\Lambda_{K3}$. In particular, in the following points (ii) and (iii) we use that if there are no overlattices of $\Lambda$ primitively embedded in $\Lambda_{K3}$, then $\Lambda$ can not be embedded in $\Lambda_{K3}$. 

{\bf The procedure.}
\begin{enumerate}
    \item[(i)] First, one has to consider the overlattices $\Lambda'$ of $\Lambda$ of maximal finite index $r$ (possibly $r=1$) as in Proposition \ref{prop: divisible}. The overlattice $\Lambda'$ of a prescribed index $r$ of $\Lambda$ is not necessarily unique. If at least one overlattice $\Lambda'$ with maximal index $r$ satisfies Corollary \ref{cor: criterion of existence embedding Nikulin} then $\Lambda'$ can be primitively embedded in $\Lambda_{K3}$: in this case $\Lambda\in\mathcal{L}_{real}$ and the procedure ends. 
    \item[(ii)] If $\Lambda$ is not as in (i), and $\ell(\Lambda')> 22-\rk(\Lambda)=22-\rk(\Lambda')$ for every possible overlattice $\Lambda'$ of maximal index $r$, then $\Lambda$ is not in $\mathcal{L}_{real}$ and the procedure ends.
    \item[(iii)] If $\Lambda$ is not as in (i) and (ii), then $\ell(\Lambda')= 22-\rk(\Lambda)=22-\rk(\Lambda')$; one has to check if conditions (2) and (3) of Theorem \ref{thm: criterion of existence embedding Nikulin} are satisfied for at least one overlattice $\Lambda'$: $\Lambda$ is contained in $\mathcal{L}_{real}$ if and only if both conditions are satisfied (for every prime $p$) for at least one overlattice $\Lambda'$ and the procedure ends.
\end{enumerate}

By Corollary \ref{cor:orbif}, the list $\mathcal{L}_{real}$ gives the list of all possible primitive symplectic surfaces. The calculations to compile the list $\mathcal{L}_{real}$ according to the previous ideas will be performed using some GAP4 and SAGE programs, and the complete algorithm is described in the Appendix \ref{sec: the algortihm GAP}.

\begin{remark}\label{rem two overlattices}
{\rm As mentioned above there is a delicate point: the overlattice of maximal index $r$ is not necessarily uniquely determined. The knowledge of $r$ suffices to check if either the conditions in Corollary \ref{cor: criterion of existence embedding Nikulin} or $\ell(\Lambda')> 22-\rk(\Lambda)=22-\rk(\Lambda')$ are satisfied. So, to apply the points (i) and (ii) of the previous procedure, it suffices to compute $r$. But, if the procedure arrives at the point (iii), i.e. if $\ell(\Lambda')=22-rk(\Lambda')$, then one has to check the discriminant form of the overlattice $\Lambda'$ in order to determine if it admits a primitive embedding in $\Lambda_{K3}$ (indeed one has to check the conditions (2) and (3) in Theorem \ref{thm: criterion of existence embedding Nikulin}). To do that one really needs the explicit form of the discriminant form.
It depends on the divisibile classes we are adding and hence one really needs to compute the lattice $\Lambda'$: it is possible that there exists two overlattices $\Lambda'_1$ and $\Lambda'_2$ of $\Lambda$ such that $r=\left|\frac{\Lambda'_1}{\Lambda}\right|=\left|\frac{\Lambda'_2}{\Lambda}\right|$ 
but the discriminant forms of $\Lambda_1'$ and of $\Lambda_2'$ are different: this implies that the generators of $\frac{\Lambda'_1}{\Lambda}$ are not mapped to the ones of $\frac{\Lambda'_2}{\Lambda}$ by an isometry of $\Lambda$.

As a consequence it is possible that only one between $\Lambda_1'$ and $\Lambda_2'$ admits a primitive embedding in $\Lambda_{K3}$. As remarked in (iii) of the previous procedure, if at least one of them admits such an embedding, $\Lambda$ is an admissible ADE configuration, i.e., $\Lambda\in\mathcal{L}_{real}$.

This implies that if a lattice $\Lambda$ does not fall into the steps (i) and (ii) of the previous procedure, in order to establish if $\Lambda\in\mathcal{L}_{real}$ one has to check the conditions (2) and (3) of Theorem \ref{thm: criterion of existence embedding Nikulin} for all the possible overlattices of $\Lambda$ such that $\ell(\Lambda')=22-rk(\Lambda')$ (and not just for one). If the conditions are satisfied for at least one overlattice, then $\Lambda\in\mathcal{L}_{real}$.

For computational reasons our algorithm produces only one overlattice $\Lambda'$ and not all of them. Hence a priori it is possible that the algorithm at the first level discards lattices which have to be considered because the choice of the divisible classes is not the right one. However, the algorithm produces a list of the lattices which are discarded because of the condition (iii) and one is able to apply to these cases the criterion given by Shimada, as explained in Remark \ref{rem: comparing with Shimada} in order to possibly modify the list $\mathcal{L}_{real}$ reintroducing the lattices which were wrongly excluded, see Remark \ref{rem: comparing with Shimada}. }
\end{remark}	

\begin{example}{\rm We provide an explicit example of the "delicate" situation mentioned above, i.e. of the existence of two different overlattices $\Lambda_1'$ and $\Lambda_2'$ of the same lattice $\Lambda$ with the same index, such that one of them cannot be primitively embedded in $\Lambda_{K3}$ and the other can. Let $\Lambda:=D_6^{\oplus 2}\oplus D_5\oplus A_1^{\oplus 2}$. We denote $d_i^{(j)}$  the basis of the $j$-th copy of a $D_h$ lattice in $\Lambda$ and $a_1^{h}$ the generator of the $h$-th copy of $A_1$. The intersection properties of the $d_i$'s are the following $d_i^2=-2$, $d_id_{i+1}=1$ if $i=2,3,4, (5)$ and  $d_1d_3=1$, all the others intersections are zero. The discriminant group of $\Lambda$ is $(\Z/2\Z)^6\times \Z/4\Z$ and $\Lambda$ contains twelve disjoint rational curves (but not 14), so one can add two 2-divisible classes and the maximal index of an overlattice of $\Lambda$ constructed as in Proposition \ref{prop: divisible} is $2^2$. 

Let $\mathcal{F}:=\{d_2^1,d_4^1,d_6^1, d_2^2,d_4^2,d_6^3,a_1^1, a_1^2\}$, then $\mathcal{F}$ contains eight disjoint rational curves and we put $v_{\mathcal{F}}:=\sum_{C\in\mathcal{F}}C/2$ to obtain one divisible class. 

Now let 
$$\mathcal{D}_1=\{d_1^1,d_2^1,d_1^2,d_2^2,d_1^3,d_2^3,a_1^1, a_1^2\}\mbox{ and }\mathcal{D}_2=\{d_1^1,d_4^1,d_6^1, d_1^2,d_2^2,d_1^3,d_2^3,a_1^1\}.$$ Both the sets $\mathcal{D}_i$'s contain eight disjoint rational curves and have exactly 4 curves in common with $\mathcal{F}$, hence, denoted $v_j:=:=\sum_{C\in\mathcal{D}_j}C/2$, one obtains two different overlattices of $\Lambda$:
$$\Lambda_1:=\langle \Lambda, v_{\mathcal{F}}, v_1\rangle \mbox{ and }\Lambda_2:=\langle \Lambda, v_{\mathcal{F}}, v_2\rangle.$$ 
We observe that $\rk(\Lambda)=\rk(\Lambda_1')=\rk(\Lambda_2')=19$, $\ell(\Lambda)=7$ and $\ell(\Lambda_1')=\ell(\Lambda_2')=3$. The discriminant group of $\Lambda_1'$ is generated by $\{(d_1^1+d_4^1+d_6^1+d_1^2+d_4^2+d_6^2)/2,\ (2d_1^2+2d_4^2+2d_6^2+d_1^3+d_2^3+2d_4^3)/4, (d_2^1+d_4^1+d_6^1+a_1^1)/2\}$ and the one of $\Lambda'_2$ by $\{(d_1^1+d_4^1+d_6^1+d_1^2+d_4^2+d_6^2)/2,\ (2d_1^2+2d_4^2+2d_6^2+d_1^3+d_2^3+2d_4^3)/4, (d_1^1+d_2^1+a_1^1)/2\}$ (where the last elements of the two bases are different). Computing the discriminant form on these bases one can check that the one of $\Lambda_1'$ satisfies (3) of Theorem \ref{thm: criterion of existence embedding Nikulin} and the one of $\Lambda_2'$ does not. This implies that only $\Lambda_1'$ admits a primitive embedding in $\Lambda_{K3}$, which in any case suffices to prove that $\Lambda$ admits a (non primitive) embedding in $\Lambda_{K3}$.}
\end{example}
\subsection{The proof of Theorem \ref{thm:classification}}\label{subsec: proof}
In this section we give the complete proof of Theorem \ref{thm:classification}, by using the output of the algorithm in the Appendix \ref{sec: the algortihm GAP}.

As explained in the previous section we have a list $\mathcal{L}_{real}$ of all the ADE configurations which effectively appear on at least one K3 surface.

We provide a partition of $\mathcal{L}_{real}$ into two subsets $\mathcal{L}_{Kum}$ and $\mathcal{L}_{irr}$ where $\Lambda$ is in $\mathcal{L}_{Kum}$ if and only if $\Lambda$ contains one of the lattices listed in Theorem \ref{thm:fujiki}: as a consequence, by Theorem \ref{thm:mainfujiki} the list $\mathcal{L}_{irr}$ provides the list of all possible irreducible symplectic surfaces, while the list $\mathcal{L}_{Kum}$ provides the list of all primitive symplectic surfaces having a finite quasi-\'etale covering by a $2-$dimensional complex torus. 

The next proposition allows us to identify the lattices in $\mathcal{L}_{Kum}$.
\begin{prop}
\label{prop: abelian covers}
There are exactly 10 ADE configurations on K3 surfaces whose contraction is a singular K3 surface which admits a quasi-\'etale cover by a complex torus and in particular is not a singular irreducible symplectic surface. These ADE configurations are exactly those listed in Theorem \ref{thm:fujiki}.
\end{prop}

\proof By Theorem \ref{thm:mainfujiki}, the ADE configurations whose contraction is a singular surface having a finite quasi-\'etale cover that is a torus have to contain one of the ADE configurations $B_{1},\cdots,B_{10}$ listed in Theorem \ref{thm:fujiki}. 

As the ADE configurations on a K3 surface have rank at most 19, and as the configurations $B_{5},\cdots,B_{10}$ have rank 19, it follows that if $B$ is an ADE configuration of $(-2)-$curves on a K3 surface and contains $B_{i}$ for $i=5,\cdots,10$, then $B=B_{i}$. 

Suppose now that $B$ is an ADE configuration of $(-2)-$curves on a K3 surface $S$ that contains $B_{1}=A_{1}^{\oplus 16}$. The unique (up to isometries) finite index overlattice of $A_1^{\oplus 16}$ which is primitively embedded in $\Lambda_{K3}$ contains five $2-$divisible classes, and there are no overlattices of an ADE configuration on K3 surfaces which contain more than five $2-$divisible classes (see \cite{NKummer}). 

If $B$ contains properly $B_{1}$, its rank has to be $r=16+h$ for some $h>0$, and its length has to be $\ell=16+x$ for some $x\geq 0$. The finite index overlattices of the lattice $\Lambda$ associated to $B$ which are due to the presence of 2-divisible classes in $\Lambda$ have rank $r$ and length $\ell'=6+x$. The condition $\ell'\leq\min(r,22-r)$ is $6+x\leq 6-h$ therefore $x\leq -h$ which is impossible, so $B=B_{1}$.

Let us now suppose that $B$ contains $B_{2}=A_{2}^{\oplus 9}$: the finite index overlattice of $A_2^{\oplus 9}$ which is contained in $\Lambda_{K3}$ contains three $3-$divisible classes and there are no overlattices of ADE configurations on K3 surfaces which contain more than three 3-divisible classes. The argument for this is similar to the one of the previous case: if there were other divisibile classes, the minimal model of the cover branched on the associated curves should be a surface with trivial canonical bundle and Euler characteristic different from 0 and 24, which is impossible. 

The same proof of the previous case implies that $B=B_{2}$ or $B=B_{2}\oplus A_{1}$. In the latter case the rank is 19, and the length of the finite index overlattice obtained adding the $3-$divisible classes is 3. Hence one has to check the discriminant form and in particular condition (2) in Theorem \ref{thm: criterion of existence embedding Nikulin}. A direct computation shows that this lattice does not admit a primitive embedding in $\Lambda_{K3}$, so $B=B_{2}$.

Similar arguments apply to the cases $(3)$ and $(4)$ of Theorem \ref{thm:fujiki}. \endproof

\begin{remark}
{\rm 
There is a more geometric argument to show that an ADE  configuration $\Lambda$ which properly contains a configuration $B_i$ in Theorem \ref{thm:fujiki} does not appear on a K3 surface, at least when  $\Lambda=B_i\oplus R$ for an ADE configuration $R$. Indeed, let $S$ be the K3 surface containing the ADE configuration $\Lambda$, $f:A\ra S$ be the cover of $S$ branched on $B_i$ and $\beta:A\ra A'$ be the contraction of $A$ to its minimal model $A'$. The inverse image $f^{-1}(R)$ of $R$ consists of rational curves, which are not $(-1)$-curves and are disjoint from the curves contracted by $\beta$ (the curves contracted by $\beta$ lie on $f^{-1}(B_i)$). So $A'$ would be a torus containing rational curves, which is impossible. The argument is more delicate if $R$ is not a direct summand of $\Lambda$, since in this case the rational curves in $A$ intersect some $(-1)$-curves and one should exclude that they are contracted by $\beta$.}
\end{remark}

We are now ready to state and prove the following, which provides the proof of points (1) and (2) of Theorem \ref{thm:classification}:

\begin{theorem}\label{theorem: output program}
The list $\mathcal{L}_{real}$ contains 5836 lattices and for all $\Lambda\in\mathcal{L}_{real}$ the exists a $20-\rk(\Lambda)$ dimensional family of K3 surfaces whose generic element contains an ADE-configuration $B$ of smooth rational curves associated $\Lambda$.
		
The list $\mathcal{L}_{Kum}$ contains 10 lattices. For each $\Lambda\in\mathcal{L}_{Kum}$, the contraction $X_B$ of a configuration $B$ associated to $\Lambda$ is a primitive symplectic surface which is not irreducible symplectic.
		
The list $\mathcal{L}_{irr}$ contains 5826 lattices. For each $\Lambda\in\mathcal{L}_{irr}$, the contraction $X_B$ of a configuration $B$ associated to $\Lambda$ is an irreducible symplectic surface.
\end{theorem}

\proof The number of lattices contained in each list is the output of the algorithm implemented on SAGE and described in details in Appendix \ref{sec_ADEconfig}. By Theorem \ref{theorem: existence of ADE configuration}, if $\Lambda\in\mathcal{L}_{real}$ and $\Lambda'$ is constructed as in the previous procedure, the family of the K3 surfaces such that $\Lambda'$ is primitively embedded in their N\'eron--Severi group generically contains K3 surfaces which admits the ADE configuration described by $\Lambda$ as set of rational curves. Since $\Lambda'$ is an overlattice of finite index of $\Lambda$, $\rk(\Lambda)=\rk(\Lambda')$ and the family of the K3 surfaces such that $\Lambda'$ is primitively embedded in their N\'eron--Severi group has dimension $20-\rk(\Lambda')$. The results on the sublist $\mathcal{L}_{Kum}$ and $\mathcal{L}_{irr}$ are a consequence of Theorem \ref{thm:mainfujiki}.\endproof
	
The procedure described above (and more precisely the algorithm given in Remark \ref{rem: appendix simple}) can be used in order to obtain a more refined information, that is the list of the ADE configurations whose contraction produce a simple symplectic surface. This proves the point (3) of Theorem \ref{thm:classification}.

\begin{defn}\label{def: mathcalT mathcalS}
We consider the following two sets of lattices:     
\begin{eqnarray*}\begin{array}{ll}\mathcal{T}=&\{D_{4}\oplus A_{4}^{\oplus 2}\oplus A_{2}^{\oplus 3}\oplus A_{1},
D_{5}\oplus A_{4}^{\oplus 2}\oplus A_{2}^{\oplus 2}\oplus A_{1}^{\oplus 2},
D_{7}\oplus A_{4}\oplus A_{2}^{\oplus 3}\oplus A_{1}^{\oplus 2},\\
&E_{6}\oplus A_{4}^{\oplus 2}\oplus A_{2}\oplus A_{1}^{\oplus 3},
E_{8}\oplus A_{4}\oplus A_{2}^{\oplus 2}\oplus A_{1}^{\oplus 3}
\};\\
\mathcal{S}=&\{B_{\mathfrak{A}_{5}},B_{\mathfrak{A}_{6}},B_{L_{2}(7)},B_{M_{20}}\}.\end{array}\end{eqnarray*}\end{defn}
\begin{theorem}
\label{thm: simple}
There are 4697 ADE configurations which correspond to simple symplectic surfaces. They are the ADE configurations whose lattice $\Lambda$ is primitively embedded in $\Lambda_{K3}$ (without adding further divisible classes) with the exception of the lattices contained in $\mathcal{S}\cup \mathcal{T}$.
\end{theorem}

\proof By Corollary \ref{cor: no lattices in S}, we have to consider the ADE lattices $\Lambda\in\mathcal{L}_{\rm irr}$ such that $\Lambda$ admits a primitive embedding in $\Lambda_{K3}$, and in particular it is not necessary to add a divisible class (since each divisible class gives rise to a cyclic cover). Moreover, the lattices in $\mathcal{S}$ correspond to Galois covers even if they do not contain divisible classes, by \cite[Table 2]{Xi}.

The lattices $B_{\mathfrak{A}_{6}},B_{L_{2}(7)},B_{M_{20}}\in\mathcal{S}$ have rank 19, so there are no higher rank lattices containing one of them and which is contained in $\mathcal{L}_{irr}$. 

The lattice $B_{\mathfrak{A}_{5}}=A_{4}^{\oplus 2}\oplus A_{2}^{\oplus 3}\oplus A_{1}^{\oplus 4}$ has rank 18.
If a rank 19 lattice $\Lambda$ contains primitively $B_{\mathfrak{A}_5}$, then there exists a Galois $\mathfrak{A}_5$-cover of the K3 surfaces containing the ADE configuration corresponding to $\Lambda$, because these K3 surfaces are at the boundary of the family of K3 surfaces whose N\'eron--Severi group is $B_{\mathfrak{A}_{5}}$ and the K3 surfaces in this family admit the required cover by \cite{Xi}.
So we look for the rank 19 lattices containing $B_{\mathfrak{A}_{5}}$.
These may be obtained in four ways: by adding to $B_{\mathfrak{A}_{5}}$ a copy of $A_1$; by adding a root to one of the $A_j$ lattices appearing as direct summands of $B_{\mathfrak{A}_{5}}$; by connecting two direct summands $A_j$ and $A_i$ of $B_{\mathfrak{A}_{5}}$ with an extra root; by connecting three direct summands $A_{i}$, $A_{j}$ and $A_{k}$ of $B_{\mathfrak{A}_{5}}$ with an extra root. The lattices one obtains in this way are: 
\begin{itemize}
    \item $A_{4}^{\oplus 2}\oplus A_{2}^{\oplus 3}\oplus A_{1}^{\oplus 5}$,
    \item $A_{4}^{\oplus 2}\oplus A_{2}^{\oplus 4}\oplus A_{1}^{\oplus 3}$,
    \item $A_{4}^{\oplus 2}\oplus A_{3}\oplus A_{2}^{\oplus 2}\oplus A_{1}^{\oplus 4}$,
    \item $A_{5}\oplus A_{4}\oplus A_{2}^{\oplus 3}\oplus A_{1}^{\oplus 4}$,
    \item $D_{5}\oplus A_{4}\oplus A_{2}^{\oplus 3}\oplus A_{1}^{\oplus 4}$,
    \item $A_{4}^{\oplus 2}\oplus A_{3}\oplus A_{2}^{\oplus 3}\oplus A_{1}^{\oplus 2}$,
    \item $A_{4}^{\oplus 3}\oplus A_{2}^{\oplus 2}\oplus A_{1}^{\oplus 3}$,
    \item $A_{6}\oplus A_{4}\oplus A_{2}^{\oplus 3}\oplus A_{1}^{\oplus 3}$,
    \item $E_{6}\oplus A_{4}\oplus A_{2}^{\oplus 3}\oplus A_{1}^{\oplus 3}$,
    \item $A_{5}\oplus A_{4}^{\oplus 2}\oplus A_{2}\oplus A_{1}^{\oplus 4}$,
    \item $A_{7}\oplus A_{4}\oplus A_{2}^{\oplus 2}\oplus A_{1}^{\oplus 4}$,
    \item $E_{7}\oplus A_{4}\oplus A_{2}^{\oplus 2}\oplus A_{1}^{\oplus 4}$,
    \item $A_{9}\oplus A_{2}^{\oplus 3}\oplus A_{1}^{\oplus 4}$,
    \item $D_{4}\oplus A_{4}^{\oplus 2}\oplus A_{2}^{\oplus 3}\oplus A_{1}$,
    \item $D_{5}\oplus A_{4}^{\oplus 2}\oplus A_{2}^{\oplus 2}\oplus A_{1}^{\oplus 2}$,
    \item $D_{7}\oplus A_{4}\oplus A_{2}^{\oplus 3}\oplus A_{1}^{\oplus 2},$
    \item $E_{6}\oplus A_{4}^{\oplus 2}\oplus A_{2}\oplus A_{1}^{\oplus 3}$,
    \item $E_{8}\oplus A_{4}\oplus A_{2}^{\oplus 2}\oplus A_{1}^{\oplus 3}.$
\end{itemize}
A direct inspection on the output of the program described in Appendix \ref{sec: the algortihm GAP}, shows that exactly the last five lattices are primitively embedded in $\Lambda_{K3}$. These are the lattices contained in the set $\mathcal{T}$ defined above.
   
By \cite[Table 2]{Xi} and the discussion above, the lattices in $\mathcal{S}\cup \mathcal{T}$ correspond to the branch locus of a cover of a K3 surface by a surface birational to a K3, hence they do not correspond to simple symplectic surfaces. 
\endproof
For each $\Lambda\in\mathcal{L}_{real}$ we now consider again point (i) of the previous algorithm: the lattice $\Lambda'$ is the finite index overlattice of $\Lambda$ (constructed as in Proposition \ref{prop: divisible}) with maximal index $r$. If $r=1$, $\Lambda=\Lambda'$ and $\Lambda$ is necessarily primitively embedded in $\Lambda_{K3}$: this determines at least one family of K3 surfaces whose N\'eron--Severi group contains primitively $\Lambda$. We don't known a priori if the primitive embedding of $\Lambda$ in $\Lambda_{K3}$ is unique up to isometries, so we don't know if this determines just one family or more than one.
	
If $r>1$, $\Lambda'$ is primitively embedded in $\Lambda_{K3}$ and this determines at least one family of K3 surfaces whose generic member admits an ADE configuration associated to $\Lambda$. But it is possible that there exists another lattice $\Lambda''$ such that $\Lambda\hookrightarrow \Lambda''\hookrightarrow \Lambda'$, where all the inclusions are of finite non trivial index, and $\Lambda''$ admits a primitive embedding in $\Lambda_{K3}$. In this case there is also a family of K3 surfaces associated to the primitive embedding $\Lambda''$ in $\Lambda_{K3}$, and the generic member of this family admits as well the ADE configuration associated to $\Lambda$. As $\Lambda'\neq \Lambda''$, we conclude that there are at least two different families of K3 surfaces whose generic member admits the same ADE configuration of curves, described by the lattice $\Lambda$. This is exactly the phenomenon described in Example \ref{esem: ADE no imples Galois} with $\Lambda=\Lambda''=A_1^{\oplus 8}$ and $\Lambda'=N$; see also Theorem \ref{theo: A1 singularities} below for the the cases $\Lambda=A_1^{\oplus k}$, $k=1,\ldots ,15$.
	
In conclusion, the same ADE configuration of rational curves may correspond to more than one family essentially for two different reasons: either there are several different finite index overlattices of $\Lambda$ which can be primitively embedded in $\Lambda_ {K3}$; or there is an overlattice of finite index (possibly 1) of $\Lambda$ having more than one primitive embedding in $\Lambda_{K3}$. This second phenomenon appears for example for the lattices $B_{\mathfrak{A}_6}$ and  $B_{L_2(7)}$, see \cite[Discussion below Question 3.1]{Wi} and the proof of Proposition \ref{prop: MG}.\\

As a matter of example we focus our attention on irreducible symplectic surfaces with the simplest possible singularities, i.e., $A_1$ singularities. In this case the same configuration of curves may corresponds to more than one family, and the first phenomenon described above appear. As above we denote $\Lambda'$ the overlattice of maximal finite index of an ADE configuration $\Lambda$ obtained as in Proposition \ref{prop: divisible}.
	
\begin{theorem}
\label{theo: A1 singularities}
Let $X$ be a symplectic surface with $k$ singular points, all of type $A_1$. Then $k\leq 16$ and $X$ is obtained by contracting an ADE configuration associated to the lattice $\Lambda\simeq A_1^{\oplus k}$ on a K3 surface $S$. 

If moreover $X$ is an irreducible symplectic surface, then $k\leq 15$ and:
\begin{enumerate}
    \item if $1\leq k\leq 7$, then $\Lambda=\Lambda'$ and $X$ is simple;
    \item if $k=8,9,10,11$ then $\Lambda'/\Lambda=\Z/2\Z$ and there are two possibilities: either $\Lambda$ is primitively embedded in $\Lambda_{K3}$ and $X$ is simple or $\Lambda'$ is primitively embedded in $\Lambda_{K3}$ and $X$ is not simple;
    \item if $k\geq 12$, then $X$ is not simple and the following possibilities appears: 
    \begin{itemize}
        \item if $k=12$ then $\Lambda'/\Lambda=(\Z/2\Z)^2$ and either $\Lambda'$ is primitively embedded in $NS(S)$ or $\Lambda''$ is primitively embedded in $NS(S)$, where $\Lambda''$ is an overlattice of index 2 of $\Lambda$; 
        \item if $k=13,14,15$ then $\Lambda'/\Lambda=(\Z/2\Z)^{k-11}$ and the minimal primitive overlattice of $\Lambda$ primitively embedded in $NS(S)$ is necessarily $\Lambda'$.
    \end{itemize}
    \item there are 21 irreducible families of symplectic surfaces with singularities of type $A_1$, among them 20 correspond to irreducible symplectic surfaces and among them 11 correspond to simple surfaces. The general member is non-projective.
\end{enumerate}
\end{theorem}
\proof The Theorem follows directly by checking the overlattices of $A_1^{\oplus k}$ constructed in Proposition \ref{prop: divisible} which admit a primitive embedding in $\Lambda_{K3}$ (applying modification of the algorithm as described in Remark \ref{rem: appendix simple}). See also \cite{NKummer} and \cite[Theorem 8.6, Corollary 8.9, Remark 8.10]{GS} for more detailed computations in case $k=16$ and $k\leq 15$ respectively.	\endproof
By the previous theorem (as by Example \ref{esem: ADE no imples Galois}) we see that the same ADE configuration could correspond to different families of irreducible symplectic surfaces, and it could happen that one of them is a family of simple symplectic surfaces, the others are not.

\section{Further remarks on singular irreducible surfaces}

Before moving to the study of Hilbert schemes of points on singular irreducible symplectic surfaces, we wish to add some further observations about irreducible symplectic surfaces.

We first notice that if $B$ is an ADE configuration of smooth rational curves that may be realized on a K3 surface, the general member of the family of K3 surfaces admitting $B$ as an ADE configuration of rational curves is non-projective, but there is at least one projective K3 surface in this family. This is made precise by the following, that was already used by Shimada in \cite{Sh}.

\begin{lemma}
\label{lem:projADE}
Let $B$ be an ADE configuration appearing on a K3 surface, and let $\mathcal{F}$ be the family of K3 surfaces admitting $B$ as an ADE configuration of rational curves.
\begin{enumerate}
    \item The general member of the family $\mathcal{F}$ is non-projective.
    \item There are K3 surfaces in the family $\mathcal{F}$ that are projective and codimension one subfamilies of $\mathcal{F}$ which correspond to families of projective K3 surfaces admitting the ADE configuration $B$ of rational curves.
\end{enumerate}
\end{lemma}

\proof The ADE configuration $B$ corresponds to a negative definite lattice $\Lambda$. Since $B$ appears as an ADE configuration of rational curves on a K3 surface, we see that $\Lambda$ is embedded in the N\'eron-Severi group of some K3 surface. It follows that $\Lambda$ is embedded in $\Lambda_{K3}$: its saturation $\Lambda^{s}$ in $\Lambda_{K3}$ is then a lattice which is primitively embedded in $\Lambda_{K3}$ and which has $\Lambda$ as a finite index sublattice. Notice that a priori $\Lambda^{s}$ is not uniquely determined by $\Lambda$, since it depends on the embedding of $\Lambda$ in $\Lambda_{K3}$. 

By the surjectivity of the period map, there is a K3 surface $S$ such that $NS(S)=\Lambda^{s}$ and the transcendental lattice $\Gamma$ is the lattice orthogonal to $\Lambda^{s}$ in $\Lambda_{K3}$. The existence of $S$ is due to the surjectivity of the period map. Since $NS(S)\simeq\Lambda^{s}$ is negative definite, we see that $S$ is non-projective. Moreover, since $\Lambda\subset NS(S)$, the surface $S$ contains an ADE configuration $B$ of rational curves. 

Now, the transcendental lattice $\Gamma$ of $S$ has signature $(3,19-\rk(\Lambda))$. Let $A\in\Gamma$ be a class with positive self-intersection, and let $\Delta$ be the lattice that is orthogonal to $A$ in $\Gamma$. Then $\Delta$ is primitively embedded in $\Gamma$, and hence in $\Lambda_{K3}$, and its signature is $(2,19-\rk(\Lambda))$. 

By the surjectivity of the period map there is a K3 surface $S'$ whose transcendental lattice is $\Delta$. By construction, the N\'eron--Severi of $S'$ contains the positive class $A$ and the lattice $\Lambda$. Hence $A$ is effective and we may assume that the ample cone of $S'$ is such that $A$ is in its the closure (this can be done up to reflections in $(-2)$-classes). In particular, $A$ has a non negative intersection with all the irreducible smooth rational curves, i.e. it is a nef divisor. Moreover, since $A^2>0$, $A$ its pseudoample and its orthogonal complement in $NS(S')$ is an overlattice of finite index (possibly 1) of $\Lambda$. By \cite[Proposition 3.2]{G}, the $(-2)$ classes of $\Lambda$ (up to sign) correspond to smooth rational curves on $S'$. It follows that $S'$ is a projective K3 surface admitting $B$ as an ADE configuration of rational curves. The families of the $NS(S')$-polarizd K3 surfaces have codimension 1 in $\mathcal{F}$ since $\rk(NS(S'))=\rk(\Lambda)+1$ and satisfy by construction the statement.
\endproof

\begin{remark} \label{rem: no 4k+2 polarization for Kummer}
{\rm The self-intersection of the class $A$ in the proof of the previous Lemma may verify some further conditions. This means that each ADE configuration which is admissible for a non-projective K3 surface is admissible also for a projective one, but \textit{a priori} the degree of the polarizations of the projective model may assume only particular values. For example, it is known that there exist projective K3 surfaces admitting sixteen disjoint rational curves and an orthogonal polarization of degree $d$ if and only if $d$ is a multiple of 4 (see \cite{NKummer})}.
\end{remark}

As a consequence of the previous Lemma we get the following result, showing that all examples of singular irreducible symplectic surface may be realized by projective surfaces. More precisely, we have the following:

\begin{theorem}
\label{theorem:projsiss}
For every ADE configuration $B$ in the list $\mathcal{L}_{irr}$ there is a projective irreducible symplectic surface whose singularities are the contraction of $B$.
\end{theorem}

\proof By Lemma \ref{lem:projADE} there is a projective K3 surface $S'$ which admits the ADE configuration $B$ of rational curves, whose N\'eron--Severi group is an overlattice of finite index of $\Z A\oplus \Lambda$, where $A$ is a class with a positive self intersection, as in the proof of Lemma \ref{lem:projADE}. Since $B\in\mathcal{L}_{irr}$, by Theorem \ref{thm:mainfujiki} the contraction of the curves in $B$ produces a singular irreducible symplectic surface, and the map $\varphi_{|A|}:S'\ra \mathbb{P}^n$ provide a projective model of $X_B$ which is projective since $S$, since it contracts exactly the curves in $B$. In particular $X_B$ is a projective singular irreducible symplectic surface as required.
\endproof
\begin{remark}{\rm In view of Theorem \ref{theorem:projsiss}, one may ask if, given an ADE configuration $B$, there exists a singular irreducible symplectic surface $X_B'$ (with singularities obtained by the contraction of $B$)  which has a certain model in a projective space.

For example, on can ask if $X_B'$ can be described as a double cover of $\mathbb{P}^2$ for a chosen ADE configuration $B$. This is equivalent to asking if there exists a K3 surface $S'$ whose N\'eron--Severi group is an overlattice of finite index of $\langle 2\rangle\oplus \Lambda_B$; indeed, denoted $A$ the generator of $\langle 2\rangle$, the map $\varphi_{|A|}$ contracts the curves in $B$ and provides a double cover $X_B\rightarrow \mathbb{P}^2$ (as $A$ is pseudoample and defines necessarily a $2:1$ morphism, see \cite{SD}). 

Since $A$ is obtained as in the proof of Lemma \ref{lem:projADE}, the existence of $A$ is equivalent to the existence of a vector with self-intersection 2 in the transcendental lattice $T_{S}$ of a K3 surface $S$, whose N\'eron--Severi group is a finite index overlattice of $\Lambda_B$. For example, by using this observation and Theorem \ref{theo: A1 singularities} and by computing the transcendental lattice of the surfaces appearing in that Theorem, one obtains that an irreducible symplectic surface whose singularities are all of type $A_1$ admits a model that is a double cover of $\mathbb{P}^2$ (which is branched on a possibly reducible sextic with simple normal crossing singularities). Nevertheless, as recalled in Remark \ref{rem: no 4k+2 polarization for Kummer}, a symplectic surface with $A_1$-singularities which is not an irreducible symplectic surface does not admit such a model.

One can ask a similar question, for example, for projective model of $X_B'$ as singular quartic in $\mathbb{P}^3$. In this case one has to replace the condition $A^2=2$ with the condition $A^2=4$, but now one has an extra problem: it is still true that $\varphi_{|A|}$ defines a map $X_B\rightarrow \mathbb{P}^3$, but this map can be $1:1$ to a quartic or $2:1$ to a quadric, according to conditions described in \cite[Proposition 5.7]{SD}. So one has to analyze the properties of the overlattice $NS(S')$ of $\Z A\oplus \Lambda_B$ in order to determine if the map to $\mathbb{P}^3$ is $1:1$ or $2:1$.}
\end{remark}

Another interesting property is that the transcendental lattice of every singular irreducible symplectic surface may be realized as the transcendental lattice of a K3 surface, as the following proves:

\begin{prop}
Let $X$ be an irreducible symplectic surface. The transcendental lattice $T(X)$ of $X$ and its N\'eron--Severi group $NS(X)$ are primitively embedded in $\Lambda_{K3}$, and there is a K3 surface $S$ such that $T(X)\simeq T(S)$. Moreover, if $X$ is non simple we have that $\rk(T(X))\leq 14$.
\end{prop}

\proof By Corollary \ref{cor:orbif} the surface $X$ is the contraction of an ADE configuration $B$ of rational curves on a K3 surface $S$. Let $f:S\ra X$ be the contraction map and let $\Lambda$ be the lattice corresponding to the ADE configuration $B$. Then $H^2(S,\Q)\simeq (H^2(X,\Z)\oplus\Lambda)\otimes \Q$ and $\Lambda\subset NS(S)$. It follows that $T(X)\simeq T(S)$ and that $NS(X)$ is isometric to the orthogonal lattice of $\Lambda$ in $NS(S)$. Since $T(S)$ and $NS(S)$ admit primitive embeddings in $\Lambda_{K3}$, the same holds for $T(X)$ and $NS(X)$. 

Finally, notice that if $X$ is non-simple, by Proposition \ref{prop:qetfp} there is a Galois cover of $S$: by \cite[Table 2]{Xi}, the rank of $NS(S)$ is at least 8, concluding the proof.\endproof

We now look at the automorphism group of a singular irreducible symplectic surface:

\begin{prop}
\label{prop:automorphisms}
Let $X$ be a primitive symplectic surface obtained as a contraction of a K3 surface $S$. Then $Aut(X)\subseteq Aut(S)$.
\end{prop}

\proof Let $B$ be the ADE configuration of rational curves on $S$ whose contraction gives $X$, and let $f:S\longrightarrow X$ the contraction morphism. For all $\alpha\in Aut(X)$,  $\tilde{\alpha}:=f^{-1}\circ\alpha\circ f$ is a birational transformation of the K3 surface $S$. It follows that $\tilde{\alpha}\in Aut(S)$.\endproof
%
%



\begin{remark}
{\rm Let $X$ be a primitive symplectic surface and $\alpha\in Aut(X)$ a finite order automorphism. The quotient $X/\alpha$ is a primitive symplectic surface if and only if $\alpha$ is symplectic, i.e., it preserves the symplectic form defined on the smooth locus of $X$. Moreover, $X/\alpha$ is obtained by contracting an ADE configuration of rational curves on $S/\widetilde{\alpha}$ (see the proof of Proposition \ref{prop:automorphisms} for the definition of $\widetilde{\alpha}$). If $W$ is the minimal model of $S/\widetilde{\alpha}$, then $X/\alpha$ is obtained contracting an ADE configuration on the K3 surface $W$, because  $S/\widetilde{\alpha}$ is obtained contracting an ADE configuration on $W$.}
\end{remark}
It is interesting to remark how the fundamental group of the smooth part of an irreducible symplectic surface is related with the lattice properties of the ADE configuration considered. So we recall the following result, due to Xiao  \cite[Lemma 2]{Xi} but re-written in our context

\begin{prop}
\label{prop: fundamentale group XB}

Let $S$ be a K3 surface which contains and ADE configuration $B$. Let $\Lambda$ be the lattice spanned by the curves of $B$ and $\Lambda'$ the minimal primitive overlattice of $\Lambda$ contained in $NS(S)$. Let $S\ra X_B$ be the contraction of $B$, $X_{B}^{s}$ the smooth locus of $X_{B}$ and $G:=\pi_{1}(X_B^{s})\simeq \pi_{1}(S-B)$. Then $$\left(\frac{G}{[G,G]}\right)^{*}=\frac{\Lambda'}{\Lambda},$$where $(G/[G,G])^{*}$ is the dual of $G/[G,G]$.
\end{prop}

We observe that $G$ is necessarily one of the group acting symplectically on a K3 surface and that $\frac{G}{[G,G]}$ is an Abelian group acting on a K3 surface. So the Abelianization of the fundamental group of the smooth part of an irreducible symplectic surface is trivial or contained in $\{\Z/n\Z,(\Z/a\Z)^2,\Z/2\Z^i, \Z/j\Z\times \Z/2\Z\}$ where $n=2,\ldots, 8$, $a=2,3,4$, $i=3,4$, $j=4,6$. For each group $G$, the datum $\Lambda'/\Lambda$ is given in \cite[Table 2]{Xi}.

\section{Hilbert schemes of points on singular symplectic surfaces}

Let $X$ be a primitive symplectic surface, and let $S$ be the K3 surface which is the smooth minimal model of $X$. The aim of this section is to discuss the following questions:

\begin{enumerate}
    \item[Q1)] Is $Hilb^n(X)$ an irreducible (or primitive) symplectic variety (or orbifold)?
    \item[Q2)] Does $Hilb^n(X)$ admit a symplectic resolution $\widetilde{Hilb^n(X)}$? In this case, is $\widetilde{Hilb^n(X)}$ an irreducible symplectic manifold?
    \item[Q3)] If yes, is $\widetilde{Hilb^n(X)}$ deformation equivalent to $Hilb^n(S)$?
\end{enumerate}

Before going more into the discussion of these questions, we first recall the definition of the Hilbert scheme of points on a singular surface, and some basic properties. 

If $X$ is a normal surface (complex or algebraic) and $n\in\mathbb{N}$, we let $X^{n}$ be the product of $n$ copies of $X$ and $Sym^{n}(X)$ the quotient of $X^{n}$ by the natural action of the group $\mathfrak{S}_{n}$ of permutations of $\{1,\cdots,n\}$. 

We define the Hilbert scheme of $n-$points on $X$ as the scheme $Hilb^{n}(X)$ parameterizing the functor mapping a scheme $\mathcal{S}$ the set of $\mathcal{S}-$flat families of $0-$dimensional subschemes on $\mathcal{S}\times X$ of length $n$ (see \cite{Ber}). It comes with a morphism $$\phi_{n,X}:Hilb^{n}(X)\longrightarrow Sym^{n}(X),$$called the Hilbert-Chow morphism of $X$.

If $X$ is smooth, by \cite{Fog} we know that $\phi_{n,X}:Hilb^{n}(X)\longrightarrow Sym^{n}(X)$ is smooth, and by \cite{Cha}, \cite{Hai}, \cite{ES} it is the blow-up of $Sym^{n}(X)$ along the diagonal $\Delta$ with reduced scheme structure, where if $p:X^{n}\longrightarrow Sym^{n}(X)$ is the quotient morphism, we let $\Delta$ be the image under $p$ of $$\{(x_{1},\cdots,x_{n})\in X^{n}\,|\,x_{i}=x_{j}\,\,{\rm for}\,\,{\rm some}\,\,i\neq j\}.$$

The first result we will need in what follows is the following:

\begin{prop}
\label{prop:genhilb2}
Let $X$ be a projective symplectic surface and $n\in\mathbb{N}$.
\begin{enumerate}
    \item The Hilbert scheme $Hilb^{n}(X)$ is a projective symplectic variety.
    \item If $n=2$, the Hilbert scheme $Hilb^{2}(X)$ has quotient singularities, and it is the blow-up of $Sym^{2}(X)$ along the diagonal $\Delta$.
\end{enumerate} 
\end{prop}

\proof The surface $X$ is covered by open affine subspaces of the form $\mathbb{A}^{2}/G$ where $G$ is a subgroup of $SL(2,\mathbb{Z})$. By Lemma 2.6 of \cite{Ber} we have that $Hilb^{n}(X)$ is covered by open subspaces of the form $Hilb^{n}(\mathbb{A}^{2}/G))$, and the restriction of the Hilbert-Chow morphism $\phi_{n,X}$ to $Hilb^{n}(\mathbb{A}^{2}/G))$ is $\phi_{n,\mathbb{A}^{2}/G}$. 

The fact that $Hilb^{n}(X)$ is covered by open subsets of the form $Hilb^{n}(\mathbb{A}^{2}/G)$ implies by Corollary 6.6 of \cite{CY} that $Hilb^{n}(X)$ is normal, irreducible, reduced and has symplectic singularities. Moreover, by Proposition 2.13 of \cite{Ber} we know that $Hilb^{n}(X)$ is projective.

Now, Proposition 2.1 of \cite{Yam} gives that $\phi_{2,\mathbb{A}^{2}/G}$ is the blow-up of $Sym^{2}(\mathbb{A}^{2}/G)$ along the diagonal: as the restriction of $\phi_{2,X}$ to $Hilb^{n}(\mathbb{A}^{2}/G)$ is $\phi_{2,\mathbb{A}^{2}/G}$, it follows that $\phi_{2,X}$ is the blow-up of $Sym^{2}(X)$ along the diagonal. 

The fact that $Hilb^{2}(X)$ has quotient singularities is proved in \cite{Yam}.\endproof

Because of this, from now on in this section we will only consider singular symplectic surfaces that are projective.

\begin{remark}
\label{remark:kahlercase}
{\rm If $X$ is a K\"ahler singular symplectic surface which is not projective, and if we define $Hilb^{n}(X)$ as the blow-up of $Sym^{n}(X)$ along the diagonal, then by Corollary 3.2.1 and Proposition 1.2.1(v) of \cite{Var} we would have that $Hilb^{n}(X)$ is K\"ahler.}
\end{remark}





We first relate the Hilbert scheme of points on a primitive symplectic surface and the Hilbert scheme of points on its smooth minimal model.

\begin{lemma}
\label{lemma: biration Hilb^n} 
Let $X$ be a projective primitive symplectic surface and $S$ its smooth minimal model. Then any resolution of the singularities of $Hilb^n(X)$ is birational to $Hilb^n(S)$.
\end{lemma}

\proof The surfaces $S$ and $X$ are birational by definition. Hence $S^n$ and $X^n$ are birational. The symmetric product $Sym^{n}(S)$ is the quotient of $S^n$ by the natural action of the group $\mathfrak{S}_n$, and the action of $\mathfrak{S}_n$ on $S^n$ induces an action of $\mathfrak{S}_n$ on $X^n$, whose quotient is $Sym^{n}(X)$.

Let $B$ be the ADE configuration of rational curves on $S$ whose contraction is $X$, and $X^{s}$ be the smooth locus of $X$. Then $X^{s}$ and $S\setminus B$ are isomorphic, and the action of $\mathfrak{S}_{n}$ on the open subset $(S\setminus B)^{n}$ of $S^{n}$ and on the open subset $(X^{s})^{n}$ of $X^{n}$ coincide. 

It follows that $Sym^{n}(S)$ and $Sym^{n}(X)$ are birational. Since the Hilbert-Chow morphism $\phi_{n,Y}:Hilb^n(Y)\longrightarrow Sym^{n}(Y)$ is birational for any variety $Y$ (see Theorem 2.16 of \cite{Ber}), we obtain the statement.
\endproof


In order to see if $Hilb^{n}(X)$ is a primitive symplectic variety when $X$ is a primitive symplectic surface, we have to calculate the space of closed holomorphic $2-$forms on the smooth locus of $Hilb^{n}(X)$, or equivalently the space of reflexive $2-$forms on $Hilb^{n}(X)$. 

To do so, we recall that if $Y$ is a normal compact K\"ahler space, we let $h^{[p],0}(Y)$ be the dimension of the space of reflexive $p-$forms on $Y$, i.e., $$h^{[p],0}(Y)=\dim H^{0}(Y,\Omega_{Y}^{[p]}).$$ These numbers are birational invariants for varieties with canonical singularities, as the following shows:

\begin{lem}
\label{lem:hp0}
Let $X$ and $Y$ be two normal compact K\"ahler spaces with canonical singularities. If $X$ and $Y$ are birational, then $h^{[p],0}(X)=h^{[p],0}(Y)$ for every $p\geq 0$.
\end{lem}
\proof Let $\widetilde{X}$ (resp. $\widetilde{Y}$) be a resolution of the singularities of $X$ (resp. of $Y$). Then by Corollary 1.8 of \cite{KS} we have that $$h^{[p],0}(X)=h^{p,0}(\widetilde{X}),\,\,\,\,\,\,\,\,\,h^{[p],0}(Y)=h^{p,0}(\widetilde{Y}).$$Since $X$ and $Y$ are birational, it follows that $\widetilde{X}$ and $\widetilde{Y}$ are birational, and since they are compact K\"ahler manifolds we have $h^{p,0}(\widetilde{X})=h^{p,0}(\widetilde{Y})$, so the result follows.\endproof

The birational invariance of the numbers $h^{[p],0}$ allows us to calculate these numbers for the Hilbert scheme of points on a primitive symplectic surface.

\begin{prop}
\label{prop:hp0hilb2}
Let $X$ be a projective primitive symplectic surface.
For every $p\geq 0$ we have that $$h^{[p],0}(Hilb^{n}(X))=\left\{\begin{array}{ll} 1, & p\equiv 0\mod 2, p\leq 2n\\ 0, & \mbox{otherwise}\end{array}\right.$$
\end{prop} 
\proof 
Recall that from Theorem  \ref{thm:mainsurfaces} there is a K3 surface $S$ such that $X$ is a contraction of an ADE configuration of smooth rational curves on $S$. By Lemma \ref{lemma: biration Hilb^n} we then have that $Hilb^{n}(X)$ and $Hilb^{n}(S)$ are two compact K\"ahler spaces with canonical singularities which are birational: by Lemma \ref{lem:hp0} we then have $$h^{[p],0}(Hilb^{n}(X))=h^{p,0}(Hilb^{n}(S)),$$and the statement follows since $Hilb^{n}(S)$ is an irreducible symplectic manifold.\endproof 

\subsection{Existence of a symplectic resolution}

The problem of the existence of a symplectic resolution for $Hilb^{n}(X)$ where $X$ is a symplectic surface is considered in \cite{CY} in the case of the quotient of an affine plane by the action of a finite subgroup of $SL(2,\mathbb{C})$. This allows us to prove the following result, answering to the questions Q1, Q2 and Q3. The following proposition and subsequent remark prove point (1) of Theorem \ref{thm:hilb2thm}.

\begin{prop}
\label{prop: Craw}
Let $X$ be a primitive symplectic surface that is projective and $n\in\mathbb{N}$. 
\begin{enumerate}
    \item The Hilbert scheme $Hilb^{n}(X)$ is a projective primitive symplectic variety of dimension $2n$ that has a unique projective symplectic resolution $\widetilde{H}$.
    \item The symplectic resolution $\widetilde{H}$ of $Hilb^{n}(X)$ is an irreducible symplectic manifold deformation equivalent to $Hilb^{n}(S)$, where $S$ is the minimal model of $X$.
\end{enumerate}\end{prop}
\proof Since $X$ is a primitive symplectic surface that is projective, then $Hilb^{n}(X)$ is a normal, projective symplectic variety by Proposition \ref{prop:genhilb2}, and by Proposition \ref{prop:hp0hilb2} we have $$h^{[2],0}(Hilb^{n}(X))=1.$$ 
Moreover, if $Y$ is a resolution of the singularities of $Hilb^{n}(X)$, as the singularities of $Hilb^{n}(X)$ are rational the Leray spectral sequence implies that $$h^{1}(Hilb^{n}(X),\mathcal{O}_{Hilb^{n}(X)})=h^{1}(Y,\mathcal{O}_{Y})=h^{1,0}(Y)=h^{1,0}(Hilb^{n}(S))=0$$where the last two equalities come from the fact that $Hilb^{n}(S)$ is birational to $Hilb^{n}(X)$ and hence to $Y$, and the fact that $Hilb^{n}(S)$ is an irreducible symplectic manifold. As a consequence, we see that $Hilb^{n}(X)$ is a primitive symplectic variety.

The remaining part of the proof is entirely due to M. Mauri: we warmly thank him for having suggested it to us.

Let $Z$ be a $\mathbb{Q}-$factorial terminalization of $Hilb^{n}(X)$, whose existence is guaranteed by \cite{BCHM}. We notice that $Z$ is a variety having terminal singularities and a holomorphic symplectic form on its smooth locus, so $Z$ is a symplectic variety. Moreover, since $Z$ is birational to $Hilb^{n}(X)$, by Proposition \ref{prop:hp0hilb2} we have that $$h^{[2],0}(Z)=h^{[2],0}(Hilb^{n}(X))=1.$$Finally, if $\widetilde{Z}$ is a resolution of the singularities of $Z$, then $\widetilde{Z}$ is a resolution of the singularities of $Hilb^{n}(X)$ as well, so the Leray spectral sequence gives $$h^{1}(Z,\mathcal{O}_{Z})=h^{1}(\widetilde{Z},\mathcal{O}_{\widetilde{Z}})=h^{1,0}(\widetilde{Z})=h^{1,0}(Hilb^{n}(S))=0.$$It follows that $Z$ is a primitive symplectic variety.

Now, notice that $Z$ and $Hilb^{n}(S)$ are two $\mathbb{Q}-$factorial projective varieties with terminal singularities which are birational: this implies that $Z$ and $Hilb^{n}(S)$ are isomorphic in codimension 1, and the birational morphism induces an isomorphism between $Pic(Z)\otimes\mathbb{Q}$ and $Pic(Hilb^{n}(S))\otimes\mathbb{Q}$. Since both $Z$ and $Hilb^{n}(S)$ are projective primitive symplectic varieties, by Theorem 6.16 of \cite{BL} we get that $Z$ and $Hilb^{n}(S)$ are locally trivially deformation equivalent. Since $Hilb^{n}(S)$ is smooth, it then follows that $Z$ is smooth as well, and hence it is a symplectic resolution of $Hilb^{n}(X)$.\endproof

\begin{remark}
{\rm 
We notice that $Hilb^{n}(X)$ and $Hilb^n(S)$ are deformation equivalent but the deformation is not locally trivial: indeed $Hilb^{n}(X)$ and $Hilb^{n}(S)$ belong to different locally trivial deformation classes, the first being singular while the second is smooth.}

\end{remark}

The study of the symplectic resolution of $Hilb^2(X)$ is considered when it is birational equivalent to the variety of lines on a certain singular cubic fourfold in \cite{BCL, LiB, Y1}.

\begin{remark}
{\rm If $X$ is a singular symplectic surface, not necessarily projective, by \cite[Corollary 6.6]{CY} there exists a unique local symplectic resolution of each singularity of $Hilb^{n}(X)$. The local symplectic resolutions of the singularities of $Hilb^{n}(X)$ glue together to give a symplectic resolution $Y$ of $Hilb^{n}(X)$. Nevertheless, this does not imply a priori that the K\"ahler metrics on the local resolutions glue together to give a global K\"ahler metric on $Y$. If $X$ is projective, Proposition \ref{prop: Craw} implies that $Y$ is K\"ahler.}

{\rm As soon as we know that $Y$ is K\"ahler, then even Proposition \ref{prop: Craw} (2) holds. Indeed if $S$ is the smooth minimal model of $X$, then $S$ is a K3 surface and $Y$ and $Hilb^{n}(S)$ are birational by Lemma \ref{lemma: biration Hilb^n}. Recall that if $Z$ is a compact K\"ahler holomorphic symplectic manifold which is birational to an irreducible symplectic manifold, then $Z$ is an irreducible symplectic manifold as well (see the proof of Corollary 6.2.7 of \cite{HL}): it follows that $Y$ is an irreducible symplectic manifold, hence by \cite[Theorem 4.6]{H} we have that $Y$ and $Hilb^{n}(S)$ are deformation equivalent.}


{\rm Finally, we notice that by Proposition 4.3 of \cite{BL1} the symplectic resolution of $Hilb^{n}(X)$ is K\"ahler for a general $X$ in a given locally trivial deformation class.} 
\end{remark}

\subsection{Hilbert schemes of two points on simple symplectic surfaces}

We have already seen that the Hilbert scheme of points on a primitive symplectic surface is a primitive symplectic variety.

It is natural to ask if more can be said if $X$ is an irreducible symplectic surface or a simple symplectic surface. The result we will prove is that if $X$ is simple, then $Hilb^{2}(X)$ is an irreducible symplectic orbifold.

To do so, we first need a result about the fundamental group of a Hilbert scheme of points on a smooth surface, not necessarily compact. The following is a classical result for compact complex surfaces (see for example \cite{Lehn}), and the proof for the non-compact case is identical: we write a proof for the convenience of the reader.

\begin{lem}
	\label{lem:pi1hilb}
	Let $S$ be a (non necessarily compact) smooth complex surface. Then there is an isomorphism $$\pi_{1}(Hilb^{2}(S))\simeq\frac{\pi_{1}(S)}{[\pi_{1}(S),\pi_{1}(S)]}.$$
\end{lem}

\proof Consider a path $\alpha:[0,1]\longrightarrow S$ and let $p_{0}:=\alpha(0)$ and $p_{1}:=\alpha(1)$. We suppose that $p_{1}\neq p_{0}$, and consider the two continuous paths $$\alpha_{1}:[0,1]\longrightarrow S\times S,\,\,\,\,\,\,\,\,\alpha_{1}(t):=(p_{0},\alpha(t)),$$
$$\alpha_{2}:[0,1]\longrightarrow S\times S,\,\,\,\,\,\,\,\,\alpha_{2}(t):=(\alpha(t),p_{0}).$$Clearly we have $\alpha_{1}(0)=(p_{0},p_{0})$, $\alpha_{1}(1)=(p_{0},p_{1})$, $\alpha_{2}(0)=(p_{0},p_{0})$ and $\alpha_{2}(1)=(p_{1},p_{0})$.

It then follows that the path $\beta:=\alpha_{1}^{-1}*\alpha_{2}$ is a path in $S\times S$ starting at $(p_{0},p_{1})$ and ending at $(p_{1},p_{0})$.

Let now $$\sigma:S\times S\longrightarrow S\times S,\,\,\,\,\,\,\,\,\,\sigma(p,q):=(q,p),$$and notice that $\beta=\sigma\circ\beta^{-1}$: indeed, for every $t\in[0,1/2]$ we have $$\beta(t)=\alpha_{1}^{-1}(2t)=\alpha_{1}(1-2t)=(p_{0},\alpha(1-2t))$$ and $$\sigma\circ\beta^{-1}(t)=\sigma(\beta(1-t))=\sigma(\alpha_{2}(2(1-t)-1))=$$ $$=(\sigma(\alpha(1-2t),p_{0}))=(p_{0},\alpha(1-2t)),$$while for every $t\in[1/2,1]$ we have $$\beta(t)=\alpha_{2}(2t-1)=(\alpha(2t-1),p_{0}),$$and $$\sigma\circ\beta^{-1}(t)=\sigma(\beta(1-t))=\sigma(\alpha_{1}(1-2(1-t)))=$$ $$=\sigma(p_{0},\alpha(2t-1))=(\alpha(2t-1),p_{0}).$$

Let now $\Delta_{S}$ be the diagonal of $S\times S$, and notice that $\Delta_{S}$ has real codimension 4 in $S\times S$: it follows that there is a path $\gamma$ in $S\times S$ which is path-homotopic to $\beta$ and such that for every $t\in[0,1]$ we have $\gamma(t)\notin\Delta_{S}$.

Notice that since $\gamma$ and $\beta$ are path-homotopic, we have that $$\sigma\circ\gamma^{-1}\simeq(\sigma\circ\beta^{-1})=\beta\simeq\gamma,$$where $\simeq$ denotes the path-homotopy equivalence. Notice that both $\gamma$ and $\sigma\circ\gamma^{-1}$ do not intersect $\Delta_{S}$, and as $\Delta_{S}$ has codimension 4 in $S\times S$, the path-homotopy equivalence between $\gamma$ and $\sigma\circ\gamma^{-1}$ may be realized in $S\times S\setminus\Delta_{S}$. It follows that $\gamma\simeq\sigma\circ\gamma^{-1}$ in $S\times S\setminus\Delta_{S}$.

Let now $Z$ be the blow-up of $S\times S$ along $\Delta_{S}$ with reduced structure and let $\rho':Z\longrightarrow S\times S$ be the blow-up morphism. If we let $E'$ be the exceptional divisor of $\rho'$, we have that $\rho'$ gives an isomorphism between $Z\setminus E'$ and $S\times S\setminus\Delta_{S}$: as consequence, the path $\gamma$ defines a path $\gamma'$ in $Z\setminus E'$ whose starting point is $z_{0}:=(\rho')^{-1}(p_{0},p_{1})$ and whose ending point is $z_{1}:=(\rho')^{-1}(p_{1},p_{0})$. Moreover, if we let $\sigma':Z\longrightarrow Z$ be the lifting of $\sigma$, then the path-homotopy between $\gamma$ and $\sigma\circ\gamma^{-1}$ gives a path-homotopy between $\gamma'$ and $\sigma'\circ(\gamma')^{-1}$ in $Z\setminus E'$.

We now let $\pi':Z\longrightarrow Hilb^{2}(S)$ be the morphism such that the following diagram 
$$\begin{CD}
	Z @ >\rho'>> S\times S\\
	@V \pi' VV @VV \pi V\\
	Hilb^{2}(S) @ >\rho>> Sym^{2}(S)
\end{CD}$$
is commutative, where $\pi:S\times S\longrightarrow Sym^{2}(S)$ is the quotient morphism for the action of $\mathfrak{S}_{2}$, and $\rho:Hilb^{2}(S)\longrightarrow Sym^{2}(S)$ is the Hilbert-Chow morphism.

Let $E$ be the exceptional divisor of $\rho$, and $\overline{\gamma}:=\pi'\circ\gamma'$: then $\overline{\gamma}$ is a loop in $Hilb^{2}(S)\setminus E$ whose base point is $y_{0}:=\pi'(z_{0})$ (i.e., the subscheme of $S$ given by the two distinc points $p_{0}$ and $p_{1}$), and we have $$\overline{\gamma}^{-1}=\pi'\circ(\gamma')^{-1}=\pi'\circ\sigma'\circ(\gamma')^{-1}\simeq\pi'\circ\gamma'=\overline{\gamma},$$and hence $\tau:=[\overline{\gamma}]\in\pi_{1}(Hilb^{2}(S)\setminus E,y_{0})$ is an element of order 2.

Now, notice that since $\Delta_{S}$ has real codimension 4 in $S\times S$, we have the following chain of group isomorphisms $$\pi_{1}(Z\setminus E',z_{0})\simeq\pi_{1}(S\times S\setminus\Delta_{S},(p_{0},p_{1}))\simeq\pi_{1}(S\times S,(p_{0},p_{1}))\simeq$$ $$\simeq\pi_{1}(S\times S,(p_{0},p_{0}))\simeq\pi_{1}(S,p_{0})\times\pi_{1}(S,p_{0}).$$

But now notice moreover that $\pi':Z\setminus E'\longrightarrow Hilb^{2}(S)\setminus E$ is a topological covering of degree 2, so there is an exact sequence $$1\longrightarrow\pi_{1}(Z\setminus E',z_{0})\longrightarrow\pi_{1}(Hilb^{2}(S)\setminus E,y_{0})\longrightarrow\mathbb{Z}/2\mathbb{Z}\longrightarrow 1.$$ Since the class $\tau\in\pi_{1}(Hilb^{2}(S)\setminus E,y_{0})$ is the class of the path $\overline{\gamma}$ whose lift to $Z\setminus E'$ is the path $\gamma'$, which is not a loop, we see that $\tau$ maps to the generator of $\mathbb{Z}/2\mathbb{Z}$. 

Moreover, since $\tau$ has order 2, we see that the previous exact sequence splits, and hence we get an isomorphism $$\pi_{1}(Hilb^{2}(S)\setminus E,y_{0})\simeq\pi_{1}(Z\setminus E',z_{0})\rtimes\mathbb{Z}/2\mathbb{Z},$$and $\tau$ acts on the normal subgroup as $$\tau\cdot[\delta]:=[\gamma'*(\sigma'\circ\delta)*(\gamma')^{-1}].$$Under the isomorphism $\pi_{1}(Z\setminus E',z_{0})\simeq\pi_{1}(S,p_{0})\times\pi_{1}(S,p_{0})$, this action becomes the action that flips the two factors, i.e., we have an isomorphism $$\pi_{1}(Hilb^{2}(S)\setminus E,y_{0})\simeq\pi_{1}(S,p_{0})\times\pi_{1}(S,p_{0})\rtimes\mathbb{Z}/2\mathbb{Z}$$ where the action of $\tau$ on the normal subgroup is $$\tau\cdot([a],[b]):=([b],[a]).$$

Now, consider a tubular neighborhood $U$ of $E$ in $Hilb^{2}(S)$, and let $r:U\longrightarrow E$ be the retraction. Then $r$ is a homotopy equivalence, and its restriction $r_{0}:U\setminus E\longrightarrow E$ is a topological fibration whose fibers are homotopicaly equivalent to $S^{1}$. We now may choose $y_{0}\in U\setminus E$ and $\tau$ to be a generator of the fundamental group of $F:=r_{0}^{-1}(r(y_{0}))$.

Using the Seifert-van Kampen Theorem we see that $\pi_{1}(Hilb^{2}(S),y_{0})$ is the push-out of the commutative diagram
$$\begin{CD}
	\pi_{1}(U\setminus E,y_{0}) @ >>> \pi_{1}(Hilb^{2}(S)\setminus E,y_{0})\\
	@VVV @VVV\\
	\pi_{1}(U,y_{0}) @ >>> \pi_{1}(Hilb^{2}(S),y_{0})
\end{CD}$$
and notice that $\pi_{1}(U,y_{0})\simeq\pi_{1}(E,r(y_{0}))$.

Now, the topological fibration $r_{0}$ gives an exact sequence $$\pi_{1}(F,y_{0})\longrightarrow\pi_{1}(U\setminus E,y_{0})\longrightarrow\pi_{1}(E,r(y_{0}))\longrightarrow 1,$$and $\pi_{1}(F,y_{0})=\mathbb{Z}\tau$. The image of $\tau$ under the morphism $\pi_{1}(F,y_{0})\longrightarrow\pi_{1}(U\setminus E,y_{0})$ is simply $\tau$ itself: we then conclude that $$\pi_{1}(Hilb^{2}(S),y_{0})\simeq\pi_{1}(Hilb^{2}(S)\setminus E,y_{0})/\langle\langle\tau\rangle\rangle\simeq\pi_{1}(S,p_{0})\times\pi_{1}(S,p_{0})\rtimes\langle\tau\rangle/\langle\langle\tau\rangle\rangle.$$

Finally, consider two elements $a,b\in\pi_{1}(S,p_{0})$. The action of $\tau$ gives $\tau(a,1)\tau=(1,a)$, and $(a,1)$ and $(1,b)$ commute. It follows that when one takes the quotient under the relation $\tau=1$, we get an identification of $(a,1)$ and $(1,a)$, and the resulting class commutes with the class of $(b,1)$. As a consequence, the morphism $$\pi_{1}(S,p_{0})\times\pi_{1}(S,p_{0})\times\langle\tau\rangle/\langle\langle\tau\rangle\rangle\longrightarrow\pi_{1}(S,p_{0})/[\pi_{1}(S,p_{0}),\pi_{1}(S,p_{0})]$$mapping $(a,b,\tau)$ to $\overline{ab}$ is an isomorphism.\endproof


\begin{prop}
\label{prop:hilb2}
Let $X$ be primitive symplectic surface and $U$ its smooth locus. Let $Hilb^{2}(X)^{s}$ be the smooth locus of $Hilb^{2}(X)$.
We have a group isomorphism $$\pi_{1}(U)/[\pi_{1}(U),\pi_{1}(U)]\ra \pi_{1}(Hilb^{2}(X)^{s}).$$
\end{prop} 

\proof Let $Z:=\{p_{1},\ldots,p_{n}\}$ be the singular locus of $X$, i.e., $Z=X\setminus U$. Moreover let $\pi:X\times X\longrightarrow Sym^{2}(X)$ be the quotient under the action of $\mathfrak{S}_{2}$ interchanging the two factors and $\rho:Hilb^{2}(X)\longrightarrow Sym^{2}(X)$ be the Hilbert-Chow morphism. Let $\Delta_{X}$ be the diagonal of $X\times X$ and $D_{X}:=\pi(\Delta_{X})$. So $\rho:Hilb^{2}(X)\longrightarrow Sym^{2}(X)$ is the blow-up of $Sym^{2}(X)$ along $D_{X}$.

Notice that $U\times U\subseteq X\times X$ is the smooth locus of $X\times X$, and the image of $U\times U$ in $Sym^{2}(X)$ under $\pi$ is $Sym^{2}(U)$. Let $Z_{i}:=\pi(\{p_{i}\}\times X)=\pi(X\times\{p_{i}\})$, and notice that $Sym^{2}(X)$ is singular along $D_{X}\cup Z_{1}\cup\ldots\cup Z_{n}$. Notice that $D_{X}$ and $Z_{i}$ are isomorphic to $X$ for every $i=1,\ldots,n$, hence the singular locus of $Sym^{2}(X)$ has codimension 2 in $Sym^{2}(X)$. 

Moreover, $Sym^{2}(U)=Sym^{2}(X)\setminus(Z_{1}\cup\ldots\cup Z_{n})$. We let $\Delta_{U}$ be the diagonal of $U$ and $D_{U}:=\pi_{|U\times U}(\Delta_{U})$, then $D_{U}=D_{X}\cap Sym^{2}(U)\subset Sym^2(X)$, therefore $Sym^{2}(U)$ is singular along $D_{U}=Sym^{2}(U)\cap D_{X}$.

It follows that $\rho^{-1}(Sym^{2}(U))$ is the blow of $Sym^{2}(U)$ along $D_{U}$: we then have that $\rho^{-1}(Sym^{2}(U))=Hilb^{2}(U)$, and since $Sym^{2}(U)$ is open in $Sym^{2}(X)$ it follows that $Hilb^{2}(U)\subseteq Hilb^{2}(X)$ is an open subset.

Moreover, notice that $Hilb^{2}(U)$ is smooth, since $U$ is a smooth surface, and hence $Hilb^{2}(U)$ is contained in the smooth locus $Hilb^{2}(X)^{s}$ of $Hilb^{2}(X)$. Notice that $2p_i:=\pi(\{p_{i}\}\times \{p_i\})\in Z_{i}$ for every $i=1,\ldots,n$, and let $W_{i}:=\rho^{-1}(Z_{i}\setminus\{2p_{i}\})$. We see that the singular locus of $Hilb^{2}(X)$ is $\overline{W}_{1}\cup\ldots\cup\overline{W}_{n}$. 

Notice that $\rho^{-1}(2p_{i})$ is isomorphic to the projectivization of the Zariski tangent space $T_{p_{i}}X$ of $X$ at $p_{i}$ (see for example \cite{LiB2}, Prop. 2.19). As $p_{i}$ is an ADE singularity, and as ADE singularities on surfaces are hypersurface singularities, the dimension of the Zariski tangent space of $X$ at $p_{i}$ is 3, so that $\rho^{-1}(2p_{i})\simeq\mathbb{P}^{2}$.

Then $\overline{W}_{i}\cap\rho^{-1}(2p_{i})$ is isomorphic to a closed subset of $\mathbb{P}^{2}$, and hence $V:=Hilb^{2}(X)^{s}\setminus Hilb^{2}(U)$ is given by a finite number of dense open subsets of $\mathbb{P}^{2}$: it follows that $V$ has real codimension 4 in $Hilb^{2}(X)^{s}$, and hence there is a group isomorphism $\pi_{1}(Hilb^{2}(U))\longrightarrow\pi_{1}(Hilb^{2}(X)^{s}).$ (\cite[Theorem 2.3, Chapter X]{God}).
As $U$ is a smooth surface, by Lemma \ref{lem:pi1hilb} we get the statement.
\endproof 

\begin{remark}
{\rm By Propositions \ref{prop: fundamentale group XB} and \ref{prop:hilb2}, if $X$ is an irreducible symplectic surface then $\pi_{1}(Hilb^2(X)^s)\in\{\{1\},\Z/n\Z,(\Z/a\Z)^2,\Z/2\Z^i, \Z/j\Z\times \Z/2\Z\}$ where $n=2,\ldots, 8$, $a=2,3,4$, $i=3,4$ and $j=4,6$.}\end{remark}

As a corollary we get the following, which completes the proof of Theorem \ref{thm:hilb2thm}:

\begin{cor}
\label{cor:simpleisv}
Let $X$ be a primitive simple symplectic surface, then $Hilb^{2}(X)$ is an irreducible symplectic orbifold of dimension 4. 
\end{cor}

\proof The Hilbert scheme $Hilb^{2}(X)$ has quotient singularities by Proposition \ref{prop:genhilb2}, hence it is an orbifold. Moreover, since a symplectic surface is primitive symplectic, by Proposition \ref{prop: Craw} we know that $Hilb^{2}(X)$ is primitive symplectic variety of dimension 4. 



Now, as $X$ is simple, the smooth locus $U$ of $X$ is simply connected. It follows that $\pi_{1}(U)=\{1\}$, and by point (1) of Proposition \ref{prop:hilb2} it follows that $Hilb^{2}(X)^{s}$ is simply connected as well. Therefore $Hilb^{2}(X)$ is an irreducible symplectic orbifold.\endproof 


We now conclude with the proof of Theorem \ref{thm:b2qualunque}, for which we still need an important ingredient about the Hodge numbers and the Betti numbers of quotients of projective orbifolds by the action of a finite group: 

\begin{lem}
\label{lem:invariant}
Let $Z$ be a complex projective variety that is an orbifold, $G$ a finite group action on $Z$ and $Y:=Z/G$. Then for every $n,p,q\geq 0$ we have $$H^{n}(Y,\mathbb{Q})\simeq H^{n}(Z,\mathbb{Q})^{G},\,\,\,\,\,\,\,\,H^{p,q}(Y)\simeq H^{p,q}(Z)^{G},$$where if $k$ is a field and $V$ is a $k-$vector space on which $G$ acts, we let $V^{G}$ be the $G-$invariant subspace of $V$.
\end{lem}

\proof This result is well-known when $Z$ is smooth. In the singular case notice that $Z$ is a CW-complex, and as it is an orbifold we have Poincar\'e duality. It follows that there is the Cartan-Leray spectral sequence in cohomology $$E_{2}^{p,q}=H^{p}(G,H^{q}(Z,\mathbb{Q}))\Rightarrow E^{p+q}=H^{p+q}(Y,\mathbb{Q}).$$As $G$ is finite and $H^{q}(Z,\mathbb{Q})$ is a $\mathbb{Q}-$vector space, we have that $H^{p}(G,H^{q}(Z,\mathbb{Q}))=0$ for every $p\geq 1$, so we get an isomorphism $$H^{n}(Y,\mathbb{Q})\simeq H^{0}(G,H^{n}(Z,\mathbb{Q}))\simeq H^{n}(Z,\mathbb{Q})^{G}$$for every $n\geq 0$.

Finally, recall that as $Z$ and $Y$ are both orbifolds, we have a Hodge decomposition on $H^{n}(Z,\mathbb{Q})$ and $H^{n}(Y,\mathbb{Q})$ for every $n\geq 0$, so for every $p,q\geq 0$ such that $p+q=n$ we have $$H^{p,q}(Y)\simeq H^{p,q}(Z)^{G},$$concluding the proof.\endproof

This allows us to prove that the second Betti number of the Hilbert scheme of 2 points on a singular irreducible symplectic surface $X$ is one more than the second Betti number of $X$, generalizing a well-known result on K3 surfaces (see for example \cite{B}). More precisely, we have:
 
\begin{prop}
\label{prop:h2hilb}
Let $X$ be a projective primitive symplectic surface. Then $$H^{2}(Hilb^{2}(X),\mathbb{Q})\simeq H^{2}(X,\mathbb{Q})\oplus\mathbb{Q}.$$
\end{prop}
	
\proof Let $b:Z\longrightarrow X\times X$ be the blow-up of $X\times X$ along the diagonal, and let $E$ be the exceptional divisor. Moreover, let $S$ be the minimal model of $X$ and $f:S\longrightarrow X$ be the resolution, so that we have $f\times f:S\times S\longrightarrow X\times X$. Let $\widetilde{b}:\widetilde{Z}\longrightarrow S\times S$ be the blow-up along the diagonal, and $\widetilde{E}$ the exceptional divisor. By the universal property of the blow-up, there is a unique morphism $g:\widetilde{Z}\longrightarrow Z$ such that $(f\times f)\circ\widetilde{b}=b\circ g$. So the following diagram is commutative:
$$\begin{CD}
	\widetilde{Z} @ >\widetilde{b}>> S\times S\\
	@VgVV @VVf\times fV\\
	Z @ >b>> X\times X
\end{CD}$$

Let $p:X\times X\longrightarrow Sym^{2}(X)$ be the quotient map and $\phi:Hilb^{2}(X)\longrightarrow Sym^{2}(X)$ be the blow-up of the diagonal. By the universal property of the blow-up, there is a unique morphism $\pi:Z\longrightarrow Hilb^{2}(X)$ such that $p\circ b=\phi\circ\pi$, that is the following diagram is commutative:  
$$\begin{CD}
	Z @ >b>> X\times X\\
	@V\pi VV @VVpV\\
	Hilb^{2}(X) @>\phi>> Sym^{2}(X)
\end{CD}$$

We divide the proof in two steps: first we prove that $g^{*}:H^{2}(Z,\mathbb{Q})\longrightarrow H^{2}(\widetilde{Z},\mathbb{Q})$ is injective; in the second step, we prove that $$H^{2}(Z,\mathbb{Q})=(p_{1,X}\circ b)^{*}H^{2}(X,\mathbb{Q})\oplus(p_{2,X}\circ b)^{*}H^{2}(X,\mathbb{Q})\oplus \mathbb{Q}\cdot E$$ and then we observe that this implies the statement.

\textbf{Step 1}. Consider the normalization $\nu:Z^{\nu}\longrightarrow Z$, and let $$g^{\nu}:\widetilde{Z}\longrightarrow Z^{\nu}$$such that $\nu\circ g^{\nu}=g$, whose existence is granted by the universal property of the normalization since $\widetilde{Z}$ is smooth, and hence normal.

Notice that $\nu$ is a birational morphism, so $b\circ\nu:Z^{\nu}\longrightarrow X\times X$ is a proper birational morphism between normal varieties. We then have that $$K_{Z^{\nu}}=(b\circ\nu)^{*}K_{X\times X}+a_{1}E_{1}+\cdots+a_{n}E_{n}$$where $E_{1},\cdots,E_{n}$ are the irreducible components of codimension 1 of the exceptional locus of $b\circ\nu$ and $a_{1},\cdots,a_{n}\in\mathbb{Q}$. As $X\times X$ is symplectic, we have $K_{X\times X}=0$, and if we let $E^{\nu}:=a_{1}E_{1}+\cdots+a_{n}E_{n}$, we get $K_{Z^{\nu}}=E^{\nu}$.

The divisor $E^{\nu}$ being a $\mathbb{Q}-$divisor, we see that $K_{Z^{\nu}}$ is $\mathbb{Q}-$Cartier. In particular, there is a positive integer $m$ such that $mK_{Z^{\nu}}=mE^{\nu}$ is a Cartier divisor.

As $S$ is a K3 surface, we have that $K_{S\times S}=0$ and 
so $$K_{\widetilde{Z}}=\widetilde{b}^{*}K_{S\times S}+\widetilde{E}=\widetilde{E}.$$
Now, notice that $g^{\nu}:\widetilde{Z}\longrightarrow Z^{\nu}$ is a resolution of the singularities of $Z^{\nu}$, and we have $$g^{\nu}_{*}(mK_{\widetilde{Z}})=g^{\nu}_{*}(m\widetilde{E})=mE^{\nu}=mK_{Z^{\nu}}.$$

As $Z^{\nu}$ is normal and $K_{Z^{\nu}}$ is $\mathbb{Q}-$Cartier, the fact that $g^{\nu}_{*}(mK_{\widetilde{Z}})=mK_{Z^{\nu}}$ gives that $Z$ has canonical singularities. By Elkik's Theorem (see \cite{El}), it then has rational singularities. In particular $R^{1}g^{\nu}_{*}\mathbb{Z}=0$, and hence $$(g^{\nu})^{*}:H^{2}(Z^{\nu},\mathbb{Z})\longrightarrow H^{2}(\widetilde{Z},\mathbb{Z})$$is injective by the Leray spectral sequence.

Now, notice that $\nu:Z^{\nu}\longrightarrow Z$ is a finite morphism, so $R^{1}\nu_{*}\mathbb{Z}=0$. Again by the Leray spectral sequence we get that $$\nu^{*}:H^{2}(Z,\mathbb{Z})\longrightarrow H^{2}(Z^{\nu},\mathbb{Z})$$is injective. It follows that $g^{*}=(g^{\nu})^{*}\circ\nu^{*}:H^{2}(Z,\mathbb{Z})\longrightarrow H^{2}(\widetilde{Z},\mathbb{Z})$ is injective, and hence that $g^{*}:H^{2}(Z,\mathbb{Q})\longrightarrow H^{2}(\widetilde{Z},\mathbb{Q})$ is injective.

\textbf{Step 2}. We now prove that $$H^{2}(Z,\mathbb{Q})=(p_{1,X}\circ b)^{*}H^{2}(X,\mathbb{Q})\oplus(p_{2,X}\circ b)^{*}H^{2}(X,\mathbb{Q})\oplus\mathbb{Q}\cdot E,$$where $p_{1,X},p_{2,X}:X\times X\longrightarrow X$ are the two projections.

To do so, let $P_{1},\cdots,P_{n}\in X$ be the singular points of $X$, $D_{1,i}:=\{P_{i}\}\times X$ and $D_{2,i}:=X\times\{P_{i}\}$, so that the singular locus of $X\times X$ is $$\Sigma_{X}=\bigcup_{j=1}^{2}\bigcup_{i=1}^{n}D_{ji}.$$Notice that each $D_{ji}$ is isomorphic to $X$, and it is a closed subvariety of $X\times X$ of codimension 2.

Let $C_{1},\cdots,C_{r}$ be the $(-2)-$curves on $S$ obtained by blowing up the singular points of $X$: the classes of these curves span a lattice $\Lambda$ of rank $r$ that is primitively embedded in $NS(S)$, and we know that $$H^{2}(S,\mathbb{Q})=f^{*}(H^{2}(X,\mathbb{Q}))\oplus\bigoplus_{i=1}^{r}\mathbb{Q}\cdot C_{i}.$$As $f^{*}:H^{2}(X,\mathbb{Q})\longrightarrow H^{2}(S,\mathbb{Q})$ is injective, we see that $b_{2}(X)=22-r$.
Notice that as $X$ is simply connected we have $$H^{2}(X\times X,\mathbb{Q})=p_{1,X}^{*}H^{2}(X,\mathbb{Q})\oplus p_{2,X}^{*}H^{2}(X,\mathbb{Q}),$$so that $b_{2}(X\times X)=44-2r$.

Let now $\Gamma_{1,i}:=C_{i}\times S$ and $\Gamma_{2,i}:=S\times C_{i}$. We then have that $$H^{2}(S\times S)=p_{1,S}^{*}H^{2}(S,\mathbb{Q})\oplus p_{2,S}^{*}H^{2}(S,\mathbb{Q})=$$ $$=(p_{1,X}\circ (f\times f))^{*}H^{2}(X,\mathbb{Q})\oplus(p_{2,X}\circ(f\times f))^{*}H^{2}(X,\mathbb{Q})\oplus\bigoplus_{j=1}^{2}\bigoplus_{i=1}^{r}\mathbb{Q}\cdot\Gamma_{ji}.$$
As $S\times S$ is smooth and the diagonal is a smooth subvariety of $S\times S$, we have $$H^{2}(\widetilde{Z},\mathbb{Q})=\widetilde{b}^{*}(H^{2}(S\times S,\mathbb{Q})\oplus\mathbb{Q}\cdot\widetilde{E},$$where $\widetilde{E}$ is the exceptional divisor. Using the description of $H^{2}(S\times S,\mathbb{Q})$ we get $$H^{2}(\widetilde{Z},\mathbb{Q})=q_{1}^{*}H^{2}(X,\mathbb{Q})\oplus q_{2}^{*}H^{2}(X,\mathbb{Q})\oplus\bigoplus_{j=1}^{2}\bigoplus_{i=1}^{r}\mathbb{Q}\cdot\widetilde{\Gamma}_{ji}\oplus\mathbb{Q}\cdot\widetilde{E},$$where $q_{i}:=p_{i,X}\circ(f\times f)\circ\widetilde{b}=p_{i,X}\circ b\circ g$ for $i=1,2$, and $\widetilde{\Gamma}_{ji}$ is the proper transform of $\Gamma_{ji}$ under $\widetilde{b}$.

Now, recall that $(f\times f)^{*}:H^{2}(X\times X,\mathbb{Q})\longrightarrow H^{2}(S\times S,\mathbb{Q})$ is injective and that $\widetilde{b}^{*}:H^{2}(S\times S,\mathbb{Q})\longrightarrow H^{2}(\widetilde{Z},\mathbb{Q})$ is injective. As $(f\times f)\circ\widetilde{b}=b\circ g$ it follows that $$b^{*}:H^{2}(X\times X,\mathbb{Q})\longrightarrow H^{2}(Z,\mathbb{Q})$$is injective, and we have $$H^{2}(Z,\mathbb{Q})=(p_{1,X}\circ b)^{*}H^{2}(X,\mathbb{Q})\oplus(p_{2,X}\circ b)^{*}H^{2}(X,\mathbb{Q})\oplus V,$$where $V$ is a $\mathbb{Q}-$vector space of dimension $p\geq 1$ that contains the class of $E$. In particular, $b_{2}(Z)=44+p-2r\geq 45-2r$.

Let us now relate $H^{2}(Z,\mathbb{Q})$ and $H^{2}(\widetilde{Z},\mathbb{Q})$. Notice that for every $j=1,2$ and $i=1,\cdots,r$ there is $1\leq k\leq n$ such that $b(g(\widetilde{\Gamma}_{ji}))=(f\times f)(\widetilde{b}(\widetilde{\Gamma}_{ji}))=D_{jk}$. It follows that $g(\widetilde{\Gamma}_{ji})$ is contained in the proper transform of $D_{jk}$ under $b$, so $g(\widetilde{\Gamma}_{ji})$ has codimension at least 2 in $Z$. As a consequence, we have that $\widetilde{\Gamma}_{ji}$ is not in the image of $g^{*}$. Hence, we have $$g^{*}(H^{2}(Z,\mathbb{Q}))\subseteq q_{1}^{*}H^{2}(X,\mathbb{Q})\oplus q_{2}^{*}H^{2}(X,\mathbb{Q})\oplus\mathbb{Q}\cdot\widetilde{E}.$$

By Step 1 we know that $g^{*}$ is injective, so $$45-2r\leq b_{2}(Z)=\dim(g^{*}H^{2}(Z,\mathbb{Q}))\leq 45-2r,$$which implies that $b_{2}(Z)=45-2r$. It follows that $V$ is a $1-$dimensional $\mathbb{Q}-$vector space that contains the class of $E$, and hence $V=\mathbb{Q}\cdot E$. It then follows that $$H^{2}(Z,\mathbb{Q})=(p_{1,X}\circ b)^{*}H^{2}(X,\mathbb{Q})\oplus(p_{2,X}\circ b)^{*}H^{2}(X,\mathbb{Q})\oplus \mathbb{Q}\cdot E.$$
We are now able to proceed with the proof of the statement. To do so, notice that the action of $\mathfrak{S}_{2}$ on $X\times X$ extends to an action on $Z$, and the morphism $\pi$ is the quotient morphism under this action. Since $$H^{2}(Z,\mathbb{Q})=(p_{1,X}\circ b)^{*}H^{2}(X,\mathbb{Q})\oplus(p_{2,X}\circ b)^{*}H^{2}(X,\mathbb{Q})\oplus \mathbb{Q}\cdot E,$$by Lemma \ref{lem:invariant} the statement follows.\endproof
	
Now, let $X$ be a primitive symplectic surface obtained by contracting an ADE configuration $B$ of rational curves on a K3 surface $S$ (the case where $X$ is smooth is the case where $B$ is empty and $X=S$). Let $\Lambda$ be the lattice associated to $B$. Then we have $$b_{2}(X)=b_{2}(S)-\rk(\Lambda)=22-\rk(\Lambda).$$
\begin{cor}\label{cor: b2}
For every $3\leq b\leq 23$ there is a $4-$dimensional irreducible symplectic orbifold $X$ such that $b_{2}(X)=b$.
\end{cor}
	
\proof The classification of simple symplectic surfaces we provided shows that for every $0\leq r\leq 19$ there is a simple symplectic surface $X$ obtained by contracting an ADE configuration $B$ of rational curves on a K3 surface such that the rank of lattice $\Lambda$ associated to $B$ is $r$ (the case $r=0$ corresponds to the case where $B$ is empty, i.e., the irreducible symplectic surface is smooth). Moreover, by Theorem \ref{theorem:projsiss} we may suppose that $X$ is projective.

By Corollary \ref{cor:simpleisv} and Proposition \ref{prop:h2hilb} then $Hilb^{2}(X)$ is a $4-$dimensional irreducible symplectic variety whose second Betti number is $$b_{2}(Hilb^{2}(X))=b_{2}(X)+1=23-r.$$Since $0\leq r\leq 19$ we see that $4\leq b_{2}(Hilb^{2}(X))\leq 23$, and $b_{2}(Hilb^{2}(X))=23$ if and only if $X$ is smooth.

This provides an example of a $4-$dimensional irreducible symplectic orbifold whose second Betti number is any $4\leq b\leq 23$. The case $b_{2}=3$ may be realized as in \cite{FuMen}: by \cite{Mon} there is an irreducible symplectic manifold $X$ of dimension 4 which carries a symplectic automorphism $\sigma$ of order 11. The quotient $X/\sigma$ is an irreducible symplectic variety by Proposition 2.15 of \cite{Per}, and as in \cite{FuMen} we have that $b_{2}(X/\sigma)=3$.\endproof

    We recall that $b_2(X)=23$ implies that $X$ is smooth, by \cite{FuMen}.
\section{Appendix: The algorithm and the GAP4 and SAGE Scripts}\label{sec: the algortihm GAP}

In this section we explain the GAP4 program \verb|ADE_K3_Vectors_for_SAGE.gap| and the SAGE program \verb|ADE_X_SAGE_v2.sage| we used to determine the singular symplectic surfaces whose quasi-\'etale covers are birational only to K3 surfaces, and therefore yielding irreducible symplectic surfaces. 

As explained in Section \ref{subsec: finidinf lists}, the first step is to present the list $\mathcal{L}_{tot}$ of all possible ADE configurations according to the Picard rank $rk$ of their associated lattices for a potential primitive symplectic surfaces. 

By Lemma \ref{lem:twocond} we have that $1\leq rk \leq 19$. Recall that, geometrically, lattices of ADE type occur as configurations of exceptional divisors of minimal resolutions of rational double points. Therefore they come with a natural intersection matrix (see e.g. \cite[Chapter 14]{Huy}). 

The list \verb|All_Possible_ADE_List| presents the list $\mathcal{L}_{tot}$ of all of ADE type lattices, and it is provided by the function \verb|CompileADEList|, which takes as input a positive integer $rk$ (the rank of the lattice) and returns the list of all possible ADE configurations whose associated lattice has rank $rk$. 

This function is purely combinatoric: it takes the rank of ADE singularities and adds them up until it reaches the rank. As a matter of example, the rank of the configuration $A_2\oplus D_4\oplus E_8$ is $14$. These data are then stored into a record \verb|dato| of the form:
\begin{center}
	\verb|dato:=rec(SingTypeA=[0,1], SingTypeD=[1], SingTypeE=[0,0,1])|.
\end{center} 

Of course this does not mean that there exists a K3 surface that is the resolution of this configuration of singularities: as explained in Section \ref{sec_ADEconfig} only few of these lattices are really attached to a K3 surface, and one needs a more delicate analysis of each lattice. 

The second step is to attach to each ADE configuration its intersection matrix. This is done by the functions \verb|BM_ADE|, that take as input the record \verb|dato| and presents as output the enriched record \verb|dato| (with the same name)  with the intersection matrix of the configuration. 

We are now in the position to check if the configuration yields a primitive symplectic surface. 

To do so we first make use of the following algorithm applying Lemma \ref{prop: divisible}:

\begin{algorithm}
	\begin{algorithmic}
		\REQUIRE List $L:= $\verb|All_Possible_ADE_List| of data, where each $dato$ is a record containing the singularities and the intersection matrix
		\FOR{$dato \in L$} 
		\IF{We can add a two divisible class to $dato$}
		\STATE Add a two divisible class to $dato$  updating with it the record;
		\ENDIF
		\ENDFOR
		\RETURN (\verb|All_Possible_ADE_List|)
	\end{algorithmic}
\end{algorithm}

After running the previous algorithm, the list  \verb|All_Possible_ADE_List| contains the data that are updated, i.e if one can add a $2-$divisible class it has been added to the record. Then we perform the very same algorithm in the search of a second $2-$divisible class, and we do this four times adding up to four 2-divisible independent vectors. 

Afterwards, instead of considering $2-$divisible classes we look at $3-$divisible classes and add all of them. We perform this new algorithm two times. 

Finally we look for $4-$divisible classes and we perform the algorithm only once. 

The function \verb|Compile_ADE_List_with_vectors| of the program, which has as input the list \verb|All_Possible_ADE_List|, outputs a printed list of lattices providing their intersection matrices, and with the list of divisible classes that we have been added. 

This last list deserves now a more careful analysis since we have to check the conditions $\ell\leq 22-rk$, and if we have equality we have to check the conditions given in Theorem \ref{thm: criterion of existence embedding Nikulin} \cite[Theorem 1.12.2]{N}. 

To do so we transfer the list to \verb|SAGE|, where \verb|associated_dynkin_lattice| returns the intersection matrix of the overlattice of each given lattice adding the divisible classes, whenever there are divisible classes that can be added.

Now we calculate the length $\ell$ of the discriminant group of all these overlattices and check if $\ell<22-rk$. If this is the case then the lattice is added to the list \verb|GoodList| of lattices that correspond to a K3 surface and also to an irreducible symplectic surface (see the Remark \ref{rem No Tori}). The elements in \verb|GoodList| are obtained as in step (i) of the Procedure in Section \ref{subsec: finidinf lists}. 

If $\ell>22-rk$, the data is added to the list \verb|DistList| of discarded cases. The elements in \verb|DistList| are obtained as in step (ii) of the Procedure in Section \ref{subsec: finidinf lists}.

Finally, if $\ell=22-rk$ the lattice is added to the list \verb|NikList| of the cases that deserve the study of the genus of their bilinear intersection form. 

At this point \verb|SAGE|'s implemented routines calculate the genus of each lattice (or intersection matrix) and we can read from it the conditions of Theorem \ref{thm: criterion of existence embedding Nikulin} in a case by case analysis: this is done by the function \verb|Nikulin_Test_Final|. This subroutine produces a list \verb|SAGEout_Nik_Good| of admissible ADE configurations on a K3 surface deduced as in (iii) of the Procedure in Section \ref{subsec: finidinf lists}.  The lattices which are in \verb|NikList| that not satisfy the conditions of Theorem \ref{thm: criterion of existence embedding Nikulin} are saved in the list given by the function \verb|Nik_Bad_List_function|. We do the final test explained in Remark \ref{rem: comparing with Shimada} to this very last list. This gives a partition of the list among lattices: some of them have to be discarded, the others are moved in the list \verb|SAGEout_Nik_Good|.

The final product of this program is the required list of configurations that correspond to irreducible symplectic surfaces $\mathcal{L}_{irr}$, which is the union of \verb|GoodList| and \verb|SAGEout_Nik_Good|.
All the scripts of the programs and the lists (SAGE input and output) can e found at the following web site.
\begin{center}
\url{https://github.com/TeoGini/Singular-K3-Surfaces}
\end{center}

\begin{remark}
\label{rem: comparing with Shimada}
{\rm As explained in Remark \ref{rem two overlattices}, there could be more than one maximal rank overlattice of a given lattice $\Lambda$. The algorithm considers only one of them, this is the reason for which the list \verb|Nik_Bad_List_function| could contains even good lattices. To identify these good lattices, we do a meticulous, case-by-case examination by using the algorithm suggested by Shimada. Indeed, in \cite{Sh} he stated that an ADE configuration appear on a K3 surface if and only if it does not contain a lattice which is contained in a finite explicit list.  All the initially different cases are pointed out with the label \verb|MODIFICATO A MANO| in the \verb|github| files \verb|Input_ADE_SAGE.zip|. Moreover, our algorithm does not account for 8-divisible classes since it corresponds just to one lattice, that we added in the end to \verb|SAGEout_Nik_Good|.}
\end{remark}

\begin{remark}
\label{rem No Tori}
{\rm Warnings: the algorithm we described discharges the configurations listed in Theorem \ref{thm:fujiki}, for example $A_1^{\oplus 16}$ that corresponds to the Kummer K3 surface. This is because we perform the algorithm adding $2-$divisible classes only four times, and in order to get the configuration $A_1^{\oplus 16}$ we would have needed to run the algorithm a fifth time to find the missing fifth $2-$divisible class one needs to add. Recall that the cover given by this divisible class is necessarily birational to a torus. As we are not interested in this case, since we already know that its contraction gives a primitive symplectic surface which is not irreducible symplectic, we decided to remove it from the analysis. In the very same way all the configurations given in Theorem \ref{thm:fujiki} are not included in our good lists. This means that the list we provide is not the list of all possible ADE configurations of rational curves on K3 surfaces, but the list of all possibile ADE configurations of rational curves on K3 surfaces whose contraction is an irreducible symplectic surface. By Proposition \ref{prop: abelian covers}, the complete list of the ADE configurations is the output of the program plus the ten ADE configurations in the list of Theorem \ref{thm:fujiki}. }
\end{remark}

\begin{remark}\label{rem: appendix simple}{\rm The algorithm always add divisible classes to the ADE configuration $\Lambda$, if possible, in order to obtain the maximal index overlattice as in (i) of the Procedure of Section \ref{subsec: finidinf lists}. It is possible to run the algorithm without adding any class, simply considering the list $\mathcal{L}_{tot}$ and perform to this list the checks in the SAGE part of the program. Doing that one obtains lists analogue to the previous ones, in particular a list of the good cases which admits a primitive embedding in $\Lambda_{K3}$. This list contains 4706 elements, which are $\mathcal{S}\cup \mathcal{T}\cup \mathcal{L}_{simp}$, where $\mathcal{S}$ and $\mathcal{T}$ are as in Definition \ref{def: mathcalT mathcalS} and $\mathcal{L}_{simp}$ are the ADE configurations $B$ appearing on a K3 surface $S$ in such a way that the contraction of $B$ on $S$ produces a simple irreducible symplectic surface $X_B$, see Theorem \ref{thm: simple}.

Similarly one can run the checks in the SAGE part of the program to a list obtained adding divisible classes of a prescribed number and type, for example to find results analogue to the ones in Theorem \ref{theo: A1 singularities} for different types of singularities.}\end{remark}

	\bigskip
	
	Alice Garbagnati,  Universit\`a degli Studi di Milano,  Dipartimento di Matematica \emph{''Federigo Enriques"}, I-20133 Milano, Italy. \emph{E-mail} \verb|alice.garbagnati@unimi.it|
	
	\bigskip
	
	Matteo Penegini, Universit\`a degli Studi di Pavia,  Dipartimento di Matematica \emph{''Felice Casorati"}, I-27100 Pavia, Italy.
	\emph{E-mail} \verb|matteo.penegini@unipv.it|

\bigskip
	
	Arvid Perego, Universit\`a degli Studi di Genova, DIMA Dipartimento di Matematica, I-16146 Genova, Italy.
	\emph{E-mail} \verb|perego@dima.unige.it|

\end{document}